\documentclass[11pt,reqno]{amsart}
\usepackage{amsthm, amsfonts, amssymb, color}
\usepackage{mathrsfs}
\usepackage{caption}
\usepackage{enumerate}
\usepackage{cases}
\usepackage{amsmath}
\usepackage[T1]{fontenc}
\usepackage{microtype}
\usepackage{footnote}
\usepackage{amstext, amsxtra}
\usepackage{enumitem}
\usepackage{booktabs}
\usepackage{tabularx}
\usepackage{tikz}
\usepackage{enumitem}
\usepackage{xfp} 
\usepackage{dsfont}
\usepackage{newtxtext}
\usepackage[varvw]{newtxmath}
\usepackage{bm}

\usepackage{multirow}

\usepackage[colorlinks, linkcolor=black, citecolor=blue,   hypertexnames=false]{hyperref}
\usepackage{graphicx}
\usepackage{subfigure}
\usepackage[symbol]{footmisc}
\allowdisplaybreaks
\def\R {\mathbb{R}}

\def\N {\mathbb{N}}

\def\Z {\mathbb{Z}}

\def\Re {\mathfrak{Re\,}}
\def\Im {\mathfrak{Im\,}}
\def\eps{\varepsilon}

\def \and {{\qquad\text{and}\qquad}}

\newtheorem{theorem}{Theorem}[section]
\newtheorem{proposition}[theorem]{Proposition}

\newtheorem{lemma}[theorem]{Lemma}

\newtheorem{remark}[theorem]{Remark}

\newcommand{\1}{\mathbf{1}}
\newcommand{\ess}{\operatorname{ess}}

\numberwithin{equation}{section}
\theoremstyle{definition}

\title[CONTROL FOR GRUSHIN-TYPE OPERATOR]
{Observability and stabilizability for Grushin-Type parabolic equations in the whole space}

\author{Shanlin Huang \quad Gengsheng Wang \quad Jiaqi Yu}

\address{Shanlin Huang,  School of Mathematics (Zhuhai), Sun Yat-sen University, Zhuhai 519082, Guangdong, China}
\email{huangshlin6@mail.sysu.edu.cn}

\address{Gengsheng Wang, Hetao Institute of Mathematics and Interdisciplinary Studies (Shenzhen),  518017, China}
\email{wanggs62@yeah.net}

\address{Jiaqi Yu,  School of Mathematics and Statistics, Huazhong University of Science and Technology,  Wuhan,  430074,  P.R. China}
\email{jiaqi\_yu@hust.edu.cn}


\begin{document}

\begin{abstract}
We study observability and stabilization for the whole-space Grushin heat equation
$$
\partial_t u-\Delta_\gamma u=0,\qquad
\Delta_\gamma:=\Delta_x+|x|^{2\gamma}\Delta_y,\qquad
(t,x,y)\in(0,T)\times\mathbb{R}^n\times\mathbb{R}^m,\quad \gamma>0.
$$
For product  sets $E=E_1\times E_2$, we identify
three observability regimes in which thickness, distance from
the degeneracy set $\{x=0\}$, and local density near that set
play different roles.

\textup{(i) In the subcritical regime  $0<\gamma<1$,}
the observability inequality holds for every $T>0$ if and only if both $E_1$ and $E_2$ are thick; hence $T_{\min}(E)=0$.

\textup{(ii) At the critical exponent  $\gamma=1$,}
if $E_2=\mathbb{R}^m$ and $E_1\subset\mathbb{R}^n$ is thick, we determine the exact minimal observable time:
$$
T_{\min}(E)=\frac{a_{\mathrm{ess}}(E_1)^2}{2n},\qquad a_{\mathrm{ess}}(E_1):=\operatorname*{ess\,inf}_{x\in E_1}|x|.
$$
Thus, the threshold depends on  the essential
distance of $E_1$ from the degeneracy set.

\textup{(iii) In the supercritical regime $\gamma>1$,}
we consider $n=1$, $E_2=\mathbb{R}^m$ and thick sets $E_1\subset\mathbb{R}$.  If  $a_{\mathrm{ess}}(E_1)>0$, observability  fails for every $T>0$, so that  $T_{\min}(E)=+\infty$. If $a_{\mathrm{ess}}(E_1)=0$, observability depends on the density of $E_1$ near $0$, measured by
$$
D_{E_1}(r):=\log\!\left(\frac{2r}{|E_1\cap(-r,r)|}\right),\qquad
L_{E_1}:=\limsup_{r\downarrow0}r^2D_{E_1}(r).
$$
If $L_{E_1}=0$, then $T_{\min}(E)=0$; if $0<L_{E_1}<\infty$, a finite sufficient observation time exists. Moreover, observability at time $T$ requires
$$
\limsup_{r\downarrow0}r^{2(1+\gamma)/(\gamma-1)}D_{E_1}(r)\le C_{T,\gamma}.
$$

In contrast, for every $\gamma>0$, exponential and complete stabilizability hold if and only if  both $E_1$ and $E_2$ are thick.  The stabilization criterion is therefore independent of $\gamma$.
\end{abstract}
	




	\maketitle
	

	\section{Introduction and main results}\label{section1}
	
	\subsection{Background and Motivation}
	We study null controllability and stabilizability for the Grushin heat equation on the whole space.
	Consider the following control system:
	\begin{equation}\label{Equ-Gru-Cont}
		\begin{cases}
			\partial_t u-\Delta_\gamma u=\1_E h,
			\qquad\text{in }(0,T)\times\R^n\times\R^m,\\
			u(0,\cdot)=u_0\in L^2(\R^{n+m}),
		\end{cases}
	\end{equation}
	where
	\begin{equation}\label{Gru-op}
		\Delta_\gamma=\Delta_x+|x|^{2\gamma}\Delta_y,
		\qquad \gamma>0.
	\end{equation}
	Given a measurable control region $E\subset\R^n\times\R^m$, we ask how the geometry of $E$ and the degeneracy parameter $\gamma$ affect (i) null controllability, (ii) the {minimal observable time}, and (iii) stabilizability.
	
	
	\begin{enumerate}[label=\textbf{(D$_{\arabic*}$)},leftmargin=2.4em,itemsep=0.45em]
		
		\item \textbf{Null controllability.}
		System \eqref{Equ-Gru-Cont} is {\it null controllable} from $E$ in time $T>0$ if, for every $u_0\in L^2(\R^{n+m})$, there exists a control $h\in L^2((0,T)\times\R^{n+m})$ such that the corresponding solution satisfies $u(T)=0$.
		
		By the Hilbert Uniqueness Method, null controllability from $E$ in time $T$ is equivalent to the following {\it observability inequality}:
		\begin{equation}\label{equ-obs}
			\|e^{T\Delta_\gamma}g\|_{L^2(\R^{n+m})}^2
			\le C_{\text{obs}}\int_0^T
			\|e^{t\Delta_\gamma}g\|_{L^2(E)}^2\,dt
		\end{equation}
		for some $C_{\text{obs}}>0$ and every $g\in L^2(\R^{n+m})$.  Accordingly, we call $E$ an \emph{observable set} if \eqref{equ-obs} holds for some $T>0$.

        \item \textbf{Minimal observable time.}
		We define
		$$
		T_{\min}(E) :=\inf\{T>0:\,\eqref{Equ-Gru-Cont} \text{ is null controllable from $E$ in time } T\},\,\quad \inf\emptyset=\infty.		
        $$
		If $T_{\min}(E)=+\infty$, the observability inequality fails for every $T>0$; we call such an $E$ an \emph{unobservable set}.
		When $0<T_{\min}(E)<+\infty$, the infimum need not be attained, so null controllability at $T=T_{\min}(E)$ may fail.
		
		\item \textbf{Stabilizability.}
		System \eqref{Equ-Gru-Cont} is {\it exponentially stabilizable} from $E$ if there is a bounded feedback $h=-Ku$ such that
		$$
		\|e^{-t(-\Delta_\gamma+\1_EK)}\|_{\mathcal L(L^2(\R^{n+m}))}
		\le Me^{-\mu t}
		\qquad (t\ge0)
		$$
		for some $M,\mu>0$.
		It is {\it completely stabilizable} from $E$  if such a feedback exists for every decay rate $\mu>0$.
		
		By \cite[Theorem 1]{TWX20}, exponential stabilizability from $E$ is equivalent to the following {\it weak observability} inequality: for every $\alpha\in(0,1)$, there exist $T>0$ and $C>0$ such that
		\begin{equation}\label{equ-w-obs}
			\|e^{T\Delta_\gamma}f\|_{L^2(\R^{n+m})}^2
			\le C\int_0^T
			\|e^{t\Delta_\gamma}f\|_{L^2(E)}^2\,dt
			+\alpha\|f\|_{L^2(\R^{n+m})}^2
		\end{equation}
		for every $f\in L^2(\R^{n+m})$. Accordingly, we call such an $E$ an \emph{exponentially stabilizable set}, or simply a \emph{stabilizable set}.
		
	\end{enumerate}

	Our goal is twofold: 
	(i) to determine how the degeneracy parameter $\gamma$ affects the geometry of {\it observable sets} and the minimal observable time $T_{\min}$; and (ii) to compare the geometry of null controllable sets with that of {\it stabilizable sets}.
	
	We briefly review some results on null controllability and stabilizability that motivate the present work.
	
	First, Beauchard, Cannarsa and Guglielmi \cite{BCG14} studied the Grushin heat equation on the bounded rectangle $\Omega=(-1,1)\times(0,1)$:
	\begin{equation}\label{Equ-Gru-Cont-bd}
		\left\{
		\begin{array}{ll}
			\partial_t u-(\partial_x^2+|x|^{2\gamma}\partial_y^2)u=\1_\omega h &\text{in }(0,T)\times\Omega,\\
			u=0 &\text{on }(0,T)\times\partial\Omega,\\
			u(0,\cdot)=u_0\in L^2(\Omega).
		\end{array}
		\right.
	\end{equation}
	Let $\omega$ be a nonempty open subset of $(0,1)\times(0,1)$. When $\gamma\ge1$, assume in addition that $\omega=(a,b)\times(0,1)$ with $0<a<b\le1$. Their result gives the trichotomy 
	\begin{equation}\label{eq:BCG-regimes}
		\begin{array}{c|c}
			0<\gamma<1 & \text{null controllability for every }T>0,\\
			\gamma=1 & \text{a positive minimal time }T_*\ge a^2/2\text{ is required},\\
			\gamma>1 & \text{null controllability fails for every }T>0.
		\end{array}
	\end{equation}
	This identifies how the degeneracy parameter changes the null controllability: arbitrary-time control ($0<\gamma < 1$), minimal-time control ($\gamma = 1$), and non-controllability ($\gamma > 1$).
	The exact minimal time $T_{\min} = a^2/2$ was later established by Beauchard, Miller and Morancey \cite{BMM15} for symmetric strips, and subsequently extended by Allonsius, Boyer and Morancey \cite{ABM21} to a single vertical strip via the moments method.
    More recently, in the multidimensional setting, Dard\'e and Trabut \cite{DT26} determined the exact minimal time $T_*=\frac{1}{2d}\inf_{x\in\omega}|x|^2$ for cylindrical observation sets with cross-section $\omega\subset B_R\subset\mathbb R^d$, $d\ge2$. They also proved failure of observability at the critical time $T=T_*$ when $d=2,3$.
    
    For unbounded domains, Alphonse and Seelmann \cite{AS24,AS25} studied heat equations associated with $(-\Delta_\gamma)^s$ on $\mathbb R^d\times\mathbb T^d$, where $\gamma\ge1$ is an integer and $s>(1+\gamma)/2$. They proved null controllability in every $T>0$ from strips $\omega_x\times\mathbb T^d$ with $\omega_x$ thick in \cite{AS24}, and from general thick sets $\omega\subset\mathbb R^d\times\mathbb T^d$ in \cite{AS25}.
    The integer restriction on $\gamma$ leaves both $0<\gamma<1$ and noninteger $\gamma>1$ outside their scope.
    Recently, Jaming and Wang \cite{JW25} established whole-space null controllability results for generalized fractional Baouendi--Grushin equations from $\omega_x\times\mathbb R^m$ with $\omega_x$ equidistributed.
	For further results on null controllability and minimal observable times for Grushin-type equations, see \cite{AT13,Koe17,LZ17,BDE20,DKR23,Tam26,Lis25} and the references therein.

	For comparison, Wang, Wang, Zhang and Zhang \cite{WWZZ19} and Egidi and Veseli\'c \cite{EV18} established a sharp thickness characterization of null controllability for the classical heat equation on $\R^n$: 
	\begin{equation}\label{eq:heat}
		\partial_t u-\Delta u=\1_E h\quad\text{in }(0,T)\times\R^n,
		\qquad u(0,\cdot)=u_0\in L^2(\R^n).
	\end{equation}
	For a measurable set $E\subset\R^n$,  
	the system \eqref{eq:heat} is null controllable from $E$ in time $T>0$ if and only if $E$ is {\it thick}: there exist $L>0$ and $\delta\in(0,1)$ such that for all $j\in\Z^n$, 
	\begin{equation}\label{equ-thick}
		\frac{|E\cap\Lambda_{L}(j)|}{|\Lambda_L(j)|}\ge \delta,
	\end{equation}
	where $\Lambda_L(j)=j+(-\frac{L}{2},\frac{L}{2})^n$.
	For the free heat equation it completely characterizes the null controllable sets and imposes no positive waiting time.
	
	Thickness also characterizes stabilizability. Huang, Wang and Wang \cite{HWW21} proved this for fractional heat equations on $\R^n$, while Alphonse and Martin \cite{AM22} established exponential and complete stabilizability criteria for broader classes of diffusive equations under suitable assumptions on the Fourier symbol. Thus, null controllability and stabilizability of the classical heat equation share the same geometric criterion.
	
	Motivated by these results, we investigate how the degeneracy parameter affects the geometric conditions for null controllability, stabilizability and the minimal observable time. We characterize stabilizable product sets by thickness of both factors, establish null controllability criteria in the different degeneracy regimes, and determine the exact minimal observable time for strips in the critical case. For strong degeneracy, we obtain further sufficient conditions and necessary conditions for null controllability in terms of local density near the degeneracy.

\subsection{Main Theorems}

	 Thickness expresses a uniform distribution  throughout the whole space: every cube of a fixed size contains a uniformly positive proportion of the set. For the Grushin operator, diffusion in the $y$-directions degenerates on $x=0$, so null controllability may also depend on how the control set is distributed near this degeneracy set. To describe this local geometry, we introduce two further notions: essential distance and local relative density.
	
\begin{enumerate}[label=\textbf{(D$_{\arabic*}$)},start=4,leftmargin=2.4em,itemsep=0.45em]
    \item \textbf{Essential distance.} For a measurable set $E_1\subset\R^n$, define
    \begin{equation}\label{eq:distance}
        a_{\mathrm{ess}}(E_1):=\operatorname*{ess\,inf}_{x\in E_1}|x| =\sup\{a\ge0:\ |E_1\cap B(0,a)|=0\}.
    \end{equation}
    We use the essential distance because changing the control set on a set of measure zero does not affect null controllability.

    \item \textbf{Local relative density.} For a measurable set $E_1\subset\R$ satisfying $a_{\mathrm{ess}}(E_1)=0$, define its local relative density and local density defect by
    \begin{equation}\label{eq:density}
        \delta_{E_1}(r):=\frac{|E_1\cap(-r,r)|}{2r},
        \qquad
        D_{E_1}(r):=\log\frac{1}{\delta_{E_1}(r)},
        \qquad r>0,
    \end{equation}
    and set
    \begin{equation}\label{eq:asym-density}
        L_{E_1}:=\limsup_{r\downarrow0}r^2D_{E_1}(r).
    \end{equation}
    The assumption $a_{\mathrm{ess}}(E_1)=0$ ensures that $\delta_{E_1}(r)>0$ for every $r>0$.
    The quantities $\delta_{E_1}$ and $D_{E_1}$ describe the local density of $E_1$ near the origin, while $L_{E_1}$ measures its asymptotic behavior. They will be used to formulate null controllability conditions in the strongly degenerate case.
\end{enumerate}

	\subsubsection{Stabilizability for $\gamma>0$ and null controllability for $0<\gamma<1$}
	
	
	\begin{theorem}\label{thm:weak}
		Let $\gamma>0$, and $E=E_1\times E_2$, where $E_1\subset\R^n,E_2\subset\R^m$ are measurable sets of positive measure. Then the following are equivalent:
		\begin{enumerate}[label=(\roman*)]
			\item $E_1$ and $E_2$ are thick;
			\item the system \eqref{Equ-Gru-Cont} is completely stabilizable from $E$;
			\item the system \eqref{Equ-Gru-Cont} is exponentially stabilizable from $E$.
		\end{enumerate}
		Furthermore, if $0<\gamma<1$, the above conditions are also equivalent to the null controllability of system \eqref{Equ-Gru-Cont} from $E$ in every time $T>0$.
	\end{theorem}
	
	We make the following remarks on Theorem~\ref{thm:weak}.
	\begin{enumerate}[label={$\boldsymbol{(a_{\arabic*})}$},leftmargin=2.5em,itemsep=0.45em]
		\item For $0<\gamma<1$, null controllable and stabilizable product sets coincide: both factors must be thick. Under this condition, null controllability holds in every positive time, so $T_{\min}(E)=0$. 
		This is consistent with the thickness characterization of null controllability for the classical heat equation on the whole space \cite{EV18,WWZZ19}, here adapted to product control sets.
		Moreover, this extends the results of \cite{JW25} from equidistributed strips to a thickness characterization of product control sets in the regime $0<\gamma<1$.
		
		\item The thickness of the two factors is reflected in two spectral inequalities.
		Thickness of $E_1$ gives the uniform spectral inequality in the $x$-variables from Theorem \ref{thm-spec-H-eta}, while thickness of $E_2$, together with the bound on the Fourier frequencies and the Logvinenko--Sereda theorem, gives the spectral inequality in the $y$-variables. When $0<\gamma<1$, the exponents in both spectral inequalities are strictly less than one, and the Lebeau--Robbiano strategy yields null controllability in every positive time. For complete stabilizability, no such restriction on the exponents is required, so the argument applies to all $\gamma>0$.
		

	\end{enumerate}

\subsubsection{Critical degeneracy: $\gamma=1$}
	
	When $\gamma=1$, thickness of both factors still characterizes stabilizability, but is no longer sufficient for null controllability. In the case $n=m=1$, Lissy \cite{Lis25} proved that, if $\overline{E_2}\neq\R$ (which means that $\R\setminus E_2$ contains at least a nonempty interval), then the system \eqref{Equ-Gru-Cont} is not null controllable from $\R\times E_2$ in any time. This motivates the study of strip sets $E_1\times\R^m$, for which we determine the minimal observable time in terms of the essential distance of $E_1$ from the origin. 
	
	
	\begin{theorem}\label{thm:mini-timel}
		Let $n,m\ge1$ and $\gamma=1$. Let $E_1\subset\R^n$ be a measurable thick set, and write $a=a_{\mathrm{ess}}(E_1)$. Then system \eqref{Equ-Gru-Cont} is null controllable from $E_1\times\R^m$ in time $T>0$ if and only if
		\begin{equation}\label{eq:Tcrit}
			T>\frac{a^2}{2n}.
		\end{equation}
		As a consequence, the minimal observable time of \eqref{Equ-Gru-Cont} with $E=E_1\times\R^m$ is exactly given by
		$$
		T_{\min}(E)=\frac{a^2}{2n}.
		$$
	\end{theorem}

    We comment on the geometry and mechanism behind this threshold.
	\begin{enumerate}[label={$\boldsymbol{(b_{\arabic*})}$},leftmargin=2.5em,itemsep=0.45em]
		
		\item
		For $n=1$, the minimal observable time is $a_{\mathrm{ess}}(E_1)^2/2$, which coincides with the threshold obtained for strips in bounded domains \cite{BMM15,ABM21} and  is consistent  with the recent multidimensional results of Dard\'{e} and Trabut \cite{DT26}.
		Theorem \ref{thm:mini-timel} treats the whole-space setting $\mathbb R^n\times\mathbb R^m$ and applies to arbitrary measurable thick sets $E_1$. 
         It also establishes failure of null controllability at the critical time $T=a_{\mathrm{ess}}(E_1)^2/(2n)$ in every dimension $n\geq1$, provided that $a_{\mathrm{ess}}(E_1)>0$.
		Unlike in the weakly degenerate case ($0<\gamma<1$), thickness no longer guarantees null controllability in arbitrary time: a positive essential distance enforces a strict minimal waiting time, whereas the same thickness condition remains sufficient for stabilizability.
		
		\item
		The time threshold $\frac{a^2}{2n}$ reflects the competition between concentration near the degeneracy set $x=0$ and diffusion in $x$.
		After partial Fourier transform in $y$, the system decouples into a family of harmonic oscillators $H_{1,\eta} = -\Delta_x + |\eta|^2|x|^2$, whose Gaussian ground states become increasingly localized near $x=0$ as $|\eta|\to\infty$. For a solution evolving from a ground state, the total $L^2$-mass at time $T$ is bounded below by $e^{-2nT|\eta|}$. Meanwhile, the semigroup evolution over a time interval of length $T$ reduces the squared $L^2$-norm on $E_1$ of the solution by the factor $e^{-a_{\mathrm{ess}}(E_1)^2|\eta|}$. Comparing these two exponential factors explains the threshold $T > \frac{a_{\ess}(E_1)^2}{2n}$. When $a_{\ess}(E_1)=0$, this concentration cost disappears, and null controllability is restored for all $T>0$.
		
	\end{enumerate}

\subsubsection{Strong degeneracy: $\gamma>1$}


We first point out that Lissy's result \cite[Theorem~1.2]{Lis25} can be extended to every $\gamma>1$:

\begin{proposition}\label{prop:transverse-gap}
		Let $n=m=1$, $\gamma>1$, and let $E_2\subset\R$ be measurable. If $\overline{E_2}\ne\R$, then system \eqref{Equ-Gru-Cont} is not null controllable from $\R\times E_2$ in any time $T>0$. 
	\end{proposition}
	
		The proof for $\gamma=1$ uses prolate spheroidal wave functions \cite{Lis25}, whereas for $\gamma>1$, we  adapt Koenig's semiclassical coherent-state construction \cite{Koe20}.
        For the resulting family of solutions, the observation term decays like  $e^{-c/h}$,  while the terminal energy decays like  $e^{-Ch^{-2/(1+\gamma)}}$ as $h\to0^+$. Since $2/(1+\gamma)<1$, the ratio of the observation term to the terminal energy tends to zero, ruling out the corresponding observability inequality. For completeness, we provide the proof in Appendix \ref{app-1}.

	This negative result leads us to consider full strips of the form $E_1\times\mathbb{R}^m$. The next theorem shows how the distribution of $E_1$ near $x=0$ affects null controllability.
	
	\begin{theorem}\label{thm:strong}
		Let $n=1$, $m\ge1$, and $\gamma>1$. Let $E_1\subset\R$ be thick, and set $E=E_1\times\R^m$.
		\begin{enumerate}[label=(\roman*),leftmargin=2em,itemsep=0.45em,topsep=0.40em]
			\item If $a_{\mathrm{ess}}(E_1)>0$, system \eqref{Equ-Gru-Cont} is not null controllable from $E$ in any time $T>0$.
			\item Let $a_{\mathrm{ess}}(E_1)=0$, and let $L_{E_1}$ be given by \eqref{eq:asym-density}. Then the following statements hold:
			\begin{enumerate}[label=(\alph*),leftmargin=2em,itemsep=0.45em,topsep=0.40em]
				\item If $L_{E_1}=0$, system \eqref{Equ-Gru-Cont} is null controllable from $E$ in every time $T>0$.
				\item If $0<L_{E_1}<\infty$, there exists $T_{\mathrm{suff}}>0$ such that system \eqref{Equ-Gru-Cont} is null controllable from $E$ in every time $T>T_{\mathrm{suff}}$.
				\item If system \eqref{Equ-Gru-Cont} is null controllable from $E$ in some time $T>0$, then
				\begin{equation}\label{eq:necessary-density}
					\limsup_{r\downarrow0}
					r^{\,\frac{2(1+\gamma)}{\gamma-1}}D_{E_1}(r)
					\le C_{T,\gamma}<\infty,
				\end{equation}
				where $D_{E_1}(r)$ is given by \eqref{eq:density}.
			\end{enumerate}
		\end{enumerate}
	\end{theorem}

The following remarks concerning Theorems \ref{thm:weak}--\ref{thm:strong} are in order.
\begin{enumerate}[label={$\boldsymbol{(c_{\arabic*})}$},leftmargin=2.5em,itemsep=0.45em]
\item For product control sets, stabilizability has the same thickness criterion for every $\gamma>0$ by Theorem~\ref{thm:weak}. In contrast, null controllability depends strongly on the degeneracy parameter $\gamma$: 
for $0<\gamma<1$, thickness of both factors remains necessary and sufficient; at $\gamma=1$, Theorem~\ref{thm:mini-timel} shows that the minimal observable time for full strips is completely determined by the essential distance from the degeneracy; for $\gamma>1$, a positive essential distance excludes null controllability in any time, whereas at zero distance, the local density near $x=0$ becomes decisive.
This gives a refined geometry compared to the bounded-domain case \cite{BCG14}: Theorem~\ref{thm:strong} distinguishes sets with $a_{\mathrm{ess}}(E_1)=0$ according to the vanishing rate of their local density near $x=0$. It yields both sufficient conditions and necessary conditions, albeit without providing a complete classification or exact minimal times for general density profiles. 
\item The local density conditions in Theorem~\ref{thm:strong} reflect the concentration of wave packets near $x=0$ at high Fourier frequencies. When the essential distance vanishes, the asymptotic density defect $L_{E_1}$ measures the loss of observability due to the sparsity of $E_1$ near $x=0$. Balancing this loss against the semigroup dissipation yields the sufficient conditions, while testing the observability inequality gives the necessary condition. The two conditions involve different powers and therefore do not give a complete characterization.
\item The following table summarizes how the degeneracy parameter $\gamma$ affects the geometry of null controllable and stabilizable sets and the minimal observable time. For simplicity, we restrict the summary to $n=m=1$ and product control sets $E=E_1\times E_2$.
		
		\begin{center}
			\small
			\setlength{\tabcolsep}{4pt}
			\renewcommand{\arraystretch}{1.5}
			\setlength{\extrarowheight}{3pt}
			\renewcommand{\tabularxcolumn}[1]{m{#1}}
			
			\begin{tabularx}{\linewidth}{
					|>{\centering\arraybackslash}m{0.12\linewidth}
					|>{\centering\arraybackslash}X
					|>{\centering\arraybackslash}m{0.33\linewidth}
					|>{\centering\arraybackslash}m{0.19\linewidth}|
				}
				\hline
				$\gamma$
				& Controllability
				& $T_{\min}$
				& Stabilizability\\
				\hline
				
				$0<\gamma<1$
				& $E_1,E_2$ thick
				& $0$
				& \multirow{4}{=}{\centering
					$E_1,\,E_2$ thick}\\
				\cline{1-3}
				
				$\gamma=1$
				& $E_1$ thick, $E_2=\R$
				& $\frac{a_{\mathrm{ess}}(E_1)^2}{2}$
				&\\
				\cline{1-3}
				
				\multirow{2}{*}{$\gamma>1$}
				& $E_1$ thick, $E_2=\R$, $a_{\mathrm{ess}}(E_1)>0$ 
				& $+\infty$
				&\\
				\cline{2-3}
				
				& $E_1$ thick,  $E_2=\R$, $a_{\mathrm{ess}}(E_1)=0$
				& {\renewcommand{\arraystretch}{1.05}%
					\setlength{\extrarowheight}{0pt}%
					$\begin{array}{@{}l@{}}
						\,\qquad L_{E_1}=0:\ T_{\min}=0,\\
						0<L_{E_1}<\infty:\ T_{\min}<\infty,\\
						\quad\text{\eqref{eq:necessary-density} fails}:
						\ T_{\min}=+\infty
					\end{array}$}
				&\\
				\hline
			\end{tabularx}
			
		\end{center}
\end{enumerate}

\subsection{Outline of the Proof}

The proofs of Theorems~~\ref{thm:weak}, \ref{thm:mini-timel}, and \ref{thm:strong} use a unified framework: Fourier decomposition  reduces the problem to  Schr\"odinger operators, and frequency-uniform spectral inequalities provide the key estimates. A Lebeau--Robbiano iteration then converts these inequalities into observability, which yields null controllability by the Hilbert Uniqueness Method. The same framework also explains how the geometry changes with $\gamma$: low frequencies detect global thickness, while high frequencies detect concentration near the degeneracy. We first present the common Fourier and spectral framework, then discuss stabilizability and null controllability in the different degeneracy regimes.

\par\medskip\noindent\textbf{Fourier decomposition and spectral estimates.}
	Taking the partial Fourier transform in $y$ reduces the operator $-\Delta_\gamma$ to a family of Schrödinger operators on $\mathbb{R}^n$:
	\begin{equation}\label{eq:fiber}
		H_{\gamma,\eta}=-\Delta_x+|\eta|^2|x|^{2\gamma}, \quad \eta\in\mathbb{R}^m.
	\end{equation}
	For $\eta\ne0$, the rescaling $x=|\eta|^{-1/(1+\gamma)}z$ yields the relation
	\begin{equation}\label{eq:scaling}
		H_{\gamma,\eta} \simeq |\eta|^{2/(1+\gamma)}H_{\gamma},
		\qquad H_{\gamma}=-\Delta_z+|z|^{2\gamma}.
	\end{equation}
	The central estimate is a spectral inequality for $E_1$, uniform in $\eta$:
	\begin{equation}\label{eq:spectra}
		\|P_\Lambda(H_{\gamma,\eta})f\|_{L^2_x} \le Ce^{C\Phi_\gamma(\Lambda)} \|\mathbf{1}_{E_1}P_\Lambda(H_{\gamma,\eta})f\|_{L^2_x},
	\end{equation}
	where $P_\Lambda(H_{\gamma,\eta})$ denotes the spectral projection of $H_{\gamma,\eta}$ up to energy $\Lambda$. This estimate stems from the quantitative propagation of smallness for solutions of elliptic equations (see \cite{LM18}).
	For $E_2$, the classical Logvinenko--Sereda theorem provides the corresponding spectral estimate.
	Uniformity in $\eta$ is essential, as fixed-frequency estimates alone are insufficient to reconstruct observability for $-\Delta_\gamma$.

\par\medskip\noindent\textbf{Thickness and stabilization for all $\gamma>0$.}
For the necessity of thickness, we use the weak observability characterization of exponential stabilizability in \cite[Theorem 1]{TWX20}. As $|\eta|\to0$, the operators $H_{\gamma,\eta}$ converge to the free Laplacian in the strong resolvent sense. Passing to the limit in the weak observability inequality gives the corresponding inequality for the heat equation on $E_1$, and hence the thickness of $E_1$. The thickness of $E_2$ follows by choosing suitable Gaussian test functions.
	
Conversely, thickness of both factors gives the spectral estimates needed for complete stabilizability. By \cite[Theorem 1.1]{LWXY22}, complete stabilizability is equivalent to the following weak observability inequality
\begin{equation}\label{eq:complete-weakobs}
	\|e^{T\Delta_\gamma}f\|_{L^2(\R^{n+m})}^2 \le D_\mu\int_0^T \|e^{t\Delta_\gamma}f\|_{L^2(E)}^2\,dt +C_\mu e^{-2\mu T}\|f\|_{L^2(\R^{n+m})}^2.
\end{equation}
The spectral estimates control the low-frequency part, while semigroup decay controls the high-frequency remainder. This yields complete stabilizability, and hence exponential stabilizability, for every $\gamma>0$.

\par\medskip\noindent\textbf{Null controllability for weak and critical degeneracy: $0<\gamma\le1$.}
	If $0<\gamma<1$, the exponent in the uniform spectral estimate \eqref{eq:spectra} grows sublinearly with energy. The Lebeau--Robbiano iteration absorbs this cost in every positive time, yielding null controllability for all $T>0$ and thus $T_{\min}=0$.
	When $\gamma=1$, the operators $H_{1,\eta}$ are harmonic oscillators. We determine the minimal observable time by proving the failure of null controllability for $0<T\le a^2/(2n)$ and positive results for $T>a^2/(2n)$, where $a=a_{\mathrm{ess}}(E_1)$. 
	
	For the negative result, using the explicit kernel of the Hermite semigroup and Gaussian tail estimates, we construct high-frequency test states for which the ratio of the observed energy to the final energy is bounded by
	\begin{equation*}
		C|\eta|^{-2}\exp\bigl((2nT-a^2)|\eta|\bigr).
	\end{equation*}
	This tends to zero as $|\eta|\to\infty$ whenever $0<T\le a^2/(2n)$, including the endpoint $T=a^2/(2n)$ due to the $|\eta|^{-2}$ term. This contradicts the observability inequality and proves the failure of null controllability in these times. 
	
	For the positive result, the key is to establish observability for $H_{1,\eta}$ uniformly in $\vert{}\eta\vert{}\ge1$. Treating $\eta$ as a parameter, we decompose the solution $e^{-tH_{1,\eta}}g$ into low and high spectral modes of $H_{1,\eta}$. For the high modes, semigroup decay compensates for the observability cost $e^{d_0\vert{}\eta\vert{}}$ from \cite{AS24}, yielding a uniform estimate when the spectral cutoff is chosen sufficiently large. For the low modes, Hermite eigenfunction estimates provide a spectral inequality with a cost of $e^{a^2\vert{}\eta\vert{}}$ (up to an arbitrarily small loss in the exponent), which is absorbed by the decay of the semigroup whenever $2nT>a^2$. Combining these estimates with the standard argument for $0<\vert{}\eta\vert{}<1$ yields null controllability for every $T>a^2/(2n)$.

\par\medskip\noindent\textbf{Null controllability for strong degeneracy: $\gamma>1$.}
	If $a_{\mathrm{ess}}(E_1)>0$, high-frequency states concentrate near $x=0$ outside the observed region, so observability fails for every $T>0$. Thus, in the strong degeneracy regime, positive results require $a_{\mathrm{ess}}(E_1)=0$. However, distance alone is insufficient,  because it does not control the relative measure of $E_1$  near $x=0$.
	
	The local density conditions arise from rescaling  $H_{\gamma,\eta}$ at the scale  $r=|\eta|^{-1/(1+\gamma)}$. Under this change of variables and a subsequent time scaling, the observability estimate for $H_{\gamma,\eta}$ on $E_1$  reduces to an estimate for $H_{\gamma}$ on  $r^{-1}E_1$,  The intersection of this scaled set with  $(-1,1)$ has measure $2\delta_{E_1}(r)$. A quantitative spectral inequality on sets of positive measure \cite{MZ25} yields the desired estimate with explicit dependence on $D_{E_1}(r)$. Upon returning to the original variables, this dependence is expressed in terms of  $r^2D_{E_1}(r)$,  which motivates the definition $L_{E_1}=\limsup_{r\downarrow0}r^2D_{E_1}(r).$ If $L_{E_1}=0$, semigroup decay gives an observability constant uniform in $\eta$ for every fixed $T>0$. If $0<L_{E_1}<\infty$, such a uniform constant exists for sufficiently large $T$. Together with the low-frequency estimates provided by thickness, this yields the  sufficient conditions for null controllability.
    
    These conditions allow the local relative density to vanish at different rates. For example, a polynomial lower bound $\delta_{E_1}(r)\ge cr^k$ implies $L_{E_1}=0$, and hence null controllability in every positive time. In contrast,  the behavior $\delta_{E_1}(r)\asymp e^{-c/r^2}$ gives $L_{E_1}=c$, so the result guarantees null controllability for sufficiently large times. In the latter case, this sufficient condition does not determine the exact value of $T_{\min}$.

	Finally, testing  the observability inequality on rescaled ground states gives the necessary condition \eqref{eq:necessary-density}. If this condition fails, null controllability fails in every positive time. The sufficient and necessary conditions involve different powers of $r$ and therefore do not give a complete characterization.
	
	The rest of the paper is organized as follows: Section~\ref{sec2} establishes the frequency-uniform spectral inequalities used throughout. Sections~\ref{sec3}--\ref{sec5} prove Theorems~\ref{thm:weak}, \ref{thm:mini-timel}, and \ref{thm:strong}, respectively. The appendix contains the proof of Proposition~\ref{prop:transverse-gap}.
    
\section{Uniform spectral inequality}\label{sec2}

    In this section, we establish the frequency-uniform spectral inequalities required for the subsequent controllability and stabilization results. The proof has three stages. We first take the partial Fourier transform in $y$ to obtain a family of Schr\"{o}dinger operators and relate them to $H_\gamma=-\Delta_x+|x|^{2\gamma}$. We then establish spatial decay and localization estimates for their spectral subspaces. Finally, we combine these estimates with the propagation of smallness to transfer information from the thick observation set and obtain the uniform spectral inequalities in Theorem~\ref{thm-spec-H-eta}.
	
	The main point of the proof is to track the explicit dependence of the estimates on the spectral cutoff $\lambda$ and the nonzero Fourier frequency $\eta$, thereby establishing a spectral inequality uniform in $\eta$. The separate cases in the last subsection address the different localization scales that arise when $\gamma$ and $|\eta|$ vary.
    
\subsection{Fourier decomposition and scaling}\label{sec2.1}
	We begin by taking a partial Fourier transform of $-\Delta_\gamma$ with respect to $y$, which turns the $y$-diffusion into a frequency-dependent potential. A subsequent scaling then reduces the resulting operator at each nonzero frequency to the model operator $H_\gamma$. This scaling will later allow us to track how the spectral estimates depend on $|\eta|$. We follow some ideas from \cite{AS25}.
	
	Recall that
	\begin{equation}\label{Gru-op-1}
		-\Delta_\gamma=-\Delta_x-|x|^{2\gamma}\Delta_y,
		\qquad (x,y)\in\R^n\times\R^m.
	\end{equation}
	We use the following convention for the Fourier transform:
	\begin{equation*}
		\hat g(\zeta)=(\mathcal Fg)(\zeta)
		=(2\pi)^{-d/2}\int_{\R^d}e^{-iz\cdot\zeta}g(z)\,dz,
		\qquad g\in\mathcal S(\R^d).
	\end{equation*}
	For $g\in\mathcal S(\R^{n+m})$, the partial Fourier transform in $y$ is therefore given by
	\begin{equation*}
		\hat{g}(x,\eta)=(\mathcal F_yg)(x,\eta)
		=(2\pi)^{-m/2}\int_{\R^m}e^{-iy\cdot\eta}g(x,y)\,dy,
		\qquad (x,\eta)\in\R^n\times\R^m.
	\end{equation*}
	Let
	\begin{equation*}
		H_{\gamma,\eta}=-\Delta_x+|\eta|^2|x|^{2\gamma},
		\quad\quad H_\gamma=-\Delta_x+|x|^{2\gamma},
	\end{equation*}
	where $x\in\R^n$.
	By \eqref{Gru-op-1}, we have
	\begin{equation*}
		\mathcal{F}_y(-\Delta_\gamma)(\eta)=-\Delta_x+|\eta|^2|x|^{2\gamma}=H_{\gamma,\eta}.
	\end{equation*}
	Recall that when $|\eta|\neq 0$, each of these operators, equipped with the respective domain
	\begin{equation*}
		D(H_{\gamma,\eta})=\{\varphi\in L^2(\R^n):\;	H_{\gamma,\eta}\varphi\in L^2(\R^n)\},
	\end{equation*}
	has purely discrete spectrum consisting of a sequence of positive eigenvalues diverging to $+\infty$.
	Let $\{\phi_{\gamma,j}^\eta\}_{j=1}^\infty$ be an orthonormal basis of $L^2(\R^n)$ consisting of eigenfunctions of $H_{\gamma,\eta}$, with corresponding eigenvalues $\lambda_{\gamma,j}^\eta$ (counted with multiplicity), satisfying
	\begin{equation}\label{equ-diag-H-eta}
		H_{\gamma,\eta}\phi_{\gamma,j}^\eta = \lambda_{\gamma,j}^\eta \phi_{\gamma,j}^\eta \quad \text{in } \R^n.
	\end{equation}
	For any $\phi \in \text{Ran}(P_\lambda(H_{\gamma,\eta}))$, we can write
	\begin{equation*}
		\phi = \sum_{\lambda_{\gamma,j}^\eta \le \lambda} \alpha_j \phi_{\gamma,j}^\eta, \quad \text{with } \alpha_j = \langle \phi, \phi_{\gamma,j}^\eta \rangle.
	\end{equation*}
	Similarly, let $\{\phi_{\gamma,j}\}_{j=1}^\infty$ be an orthonormal basis of eigenfunctions for $H_\gamma$ with eigenvalues $\lambda_{\gamma,j}$ (counted with multiplicity), satisfying
	\begin{equation}\label{equ-diag-H}
		H_{\gamma}\phi_{\gamma,j} = \lambda_{\gamma,j} \phi_{\gamma,j} \quad \text{in } \R^n.
	\end{equation}
	
	For every non-zero frequency $\eta\in\R^m\setminus\{0\}$, define the unitary operator $\mathcal M_{\gamma,\eta}$ on $L^2(\R^n)$ by
	\begin{equation*}
		\mathcal{M}_{\gamma,\eta} g=|\eta|^{\frac{n}{2(1+\gamma)}}g(|\eta|^{\frac{1}{1+\gamma}}\cdot),\quad g\in L^2(\R^n).
	\end{equation*}
	A straightforward computation shows that
	\begin{equation}\label{equ-scaling-H}
		(\mathcal{M}_{\gamma,\eta})^*H_{\gamma,\eta} \mathcal{M}_{\gamma,\eta}=|\eta|^{\frac{2}{1+\gamma}}H_\gamma.
	\end{equation}
	Thus, we may choose the eigenfunctions of $H_{\gamma,\eta}$ so that
	\begin{equation}\label{equ-eigen-connect}
		\phi_{\gamma,j}^\eta=\mathcal{M}_{\gamma,\eta}\phi_{\gamma,j},\quad\mbox{ and }\quad\lambda_{\gamma,j}^\eta=|\eta|^{\frac{2}{1+\gamma}}\lambda_{\gamma,j}.
	\end{equation}

	\subsection{Spatial decay and localization estimates}\label{sec:doubling}
	We first establish localization estimates and then derive a propagation of smallness estimate. These estimates play a key role in the proof of the uniform spectral inequality in Theorem~\ref{thm-spec-H-eta}.
	
	\subsubsection{Localization estimates}
	We derive the decay property of $\phi\in \text{Ran}(P_\lambda(H_{\gamma,\eta}))$ at infinity with explicit dependence on $\lambda$ and $|\eta|$ in both the constants and the spatial regions involved. When $0<\gamma<1/2$, the potential is not Lipschitz at the origin, so we use an annular localization estimate away from the origin instead. These regions will serve as the domains for the propagation of smallness argument below.

	\begin{lemma}\label{lem:doub-phi}
		Let $\gamma>0$, $|\eta|>0$. There exists $\hat{C}>1$ depending on $n,\gamma$, such that for any $\lambda\ge\lambda_{\gamma,1}^\eta$ and every $\phi\in \text{Ran}(P_\lambda(H_{\gamma,\eta}))$, the following estimates hold. Define
		\begin{equation}\label{equ-M}
			M_{\lambda,\eta}=({\lambda}{|\eta|^{-2/(1+\gamma)}}+1)^{\frac{1}{2\gamma}}+\ln (|\eta|+1),
		\end{equation}
		and set 
		\begin{equation}\label{equ-R-lambda-eta}
			R_{\lambda,\eta}:=\hat C|\eta|^{-\frac{1}{1+\gamma}}M_{\lambda,\eta}.
		\end{equation}
		Then
		\begin{equation}\label{equ-L2bound-L}
			\|\phi\|_{L^2(\R^n)}^2\le 2 \|\phi\|_{L^2( B(0, R_{\lambda,\eta}))}^2,
		\end{equation}
		and
		\begin{equation}\label{equ-H1bound-L}
			\|\phi\|_{H^1(\R^n)}^2\le 2\|\phi\|_{H^1(B(0,R_{\lambda,\eta}))}^2.
		\end{equation}
	\end{lemma}
	\begin{proof}
		Following the proof of \cite[Theorem 1.4]{DSV24}, we show that there exists $C_2>0$, depending only on $n$ and $\gamma$, such that for every $r'>0$, the following decay estimate
		\begin{equation}\label{equ-local-estim}
			\|\phi\|_{H^1\bigl(\R^n\setminus B(0,2|\eta|^{-\frac{1}{1+\gamma}}(r'+C_2M_{\lambda,\eta}))\bigr)}^2
			\le e^{-r'}\|\phi\|_{L^2(\R^n)}^2
		\end{equation}
		holds for all $\lambda\ge\lambda_{\gamma,1}^\eta$ and $\phi\in\operatorname{Ran}(P_\lambda(H_{\gamma,\eta}))$.
		Indeed, for any $r>0$,
		\begin{align}\label{eq:doub5}
			\|\phi\|_{H^1(\R^n\setminus  B_r)}^2&=\|\phi\|_{L^2(\R^n\setminus B_r)}^2+\||\nabla\phi|\|_{L^2(\R^n\setminus B_r)}^2\notag\\
			&\le e^{-\rho_\eta  r/2}\Big(\|e^{\rho_\eta |x|/2}\phi\|_{L^2(\R^n)}^2+\|e^{\rho_\eta |x|/4}|\nabla\phi|\|_{L^2(\R^n)}^2\Big),
		\end{align}
		where we write
		\begin{equation*}
			{\rho_\eta}:=|\eta|^{1/(1+\gamma)}.
		\end{equation*}
		We rewrite every $\phi\in \text{Ran}(P_\lambda(H_{\gamma,\eta}))$ as
		\begin{equation}\label{equ-lem1-2}
			\phi=\sum_{j=1}^{N_{\lambda,\eta}}\phi_j,\quad\;\mbox{where }N_{\lambda,\eta}=\max\{j\in\N^+:\;\lambda_{\gamma,j}^\eta\le \lambda\}\quad\mbox{and }\;\phi_j=\alpha_j\phi_{\gamma,j}^\eta.
		\end{equation}
		Let
		\begin{equation*}
			\Lambda_{\lambda,\eta}:=\lambda{|\eta|^{-2/(1+\gamma)}},
		\end{equation*}
		it follows from (3.9) in \cite{DSV24} that
		\begin{align*}
			N_{\lambda,\eta}=\#\{j:\,|\eta|^\frac{2}{1+\gamma}\lambda_{\gamma,j}\le\lambda\}=\#\{j:\,\lambda_{\gamma,j}\le\Lambda_{\lambda,\eta}\}
			\lesssim_{n,\gamma}\big(1+\Lambda_{\lambda,\eta}\big)^{1+\frac{n}{2}(1+\frac{1}{\gamma})}.
		\end{align*}
		On the other hand, using the expansion \eqref{equ-lem1-2} and H\"{o}lder's inequality, we may estimate
		\begin{align}\label{eq:doub1}
			\|e^{\rho_\eta|x|/2}\phi\|_{L^2(\R^n)}^2
			&\le  \Big(\sum_{j=1}^{N_{\lambda,\eta}}\|e^{\rho_\eta|x|/2}\phi_j\|_{L^2(\R^n)}\Big)^2\notag\\
			&\le N_{\lambda,\eta}\cdot \sum_{j=1}^{N_{\lambda,\eta}}\|e^{\rho_\eta|x|/2}\phi_j\|_{L^2(\R^n)}^2\notag\\
			&=N_{\lambda,\eta}\cdot \sum_{j=1}^{N_{\lambda,\eta}}|\alpha_j|^2\ \|e^{|x|/2}\phi_{\gamma,j}\|_2^2,
		\end{align}
		where the last equality follows from \eqref{equ-lem1-2} and a change of variables.
		Similarly, taking into account $|\nabla\phi|\le\sum_{j=1}^{N_{\lambda,\eta}}|\nabla\phi_j|$, by scaling and gradient computation, we have
		\begin{align}\label{eq:doub2}
			\big\|e^{\rho_\eta|x|/4}|\nabla\phi|\big\|_{L^2(\R^n)}^2\le |\eta|^{\frac{2}{1+\gamma}}N_{\lambda,\eta}\cdot \sum_{j=1}^{N_{\lambda,\eta}}|\alpha_j|^2\ \big\|e^{|x|/4}|\nabla\phi_{\gamma,j}|\big\|_2^2.
		\end{align}
		
		We now estimate $\big\|e^{|x|/2}\phi_{\gamma,j}\big\|_2$ and $\big\|e^{|x|/4}|\nabla\phi_{\gamma,j}|\big\|_2$.
		It follows from \cite[Remark 3.3]{ZZG25} (see also \cite{DSV24}) that
		if $H_\gamma f=\lambda f$, then for $R$ satisfying $R^{2\gamma}\ge\max\{ \lambda+2,1\}$,
		\begin{equation}\label{eq-weight}
			\|e^{|x|/2}f\|_{L^2(\mathbb R^n)}^2
			\le 7e^{R+1}\|f\|_{L^2(\mathbb R^n)}^2,
		\end{equation}
		and
		\begin{equation}\label{eq-weight-grad}
			\|e^{|x|/4}\nabla f\|_{L^2(\mathbb R^n)}^2
			\le C(1+R^{2\gamma})e^{2(R+1)}\|f\|_{L^2(\mathbb R^n)}^2.
		\end{equation}
		Noting that
		$$
		\phi_{\gamma,j}\in\text{Ran}(P_{\Lambda_{\lambda,\eta}}(H_\gamma)),\qquad j\le N_{\lambda,\eta},
		$$
		we choose $R=\big(\Lambda_{\lambda,\eta}+2\big)^{\frac{1}{2\gamma}}$.
		Then, combining \eqref{eq:doub1}, \eqref{eq:doub2}, \eqref{eq-weight}, and \eqref{eq-weight-grad} yields
		\begin{align}\label{eq:doub3}
			\|e^{\rho_\eta|x|/2}\phi\|_2^2
			&\le C_\gamma N_{\lambda,\eta}\sum_{j=1}^{N_{\lambda,\eta}}\exp\big({C_\gamma(\Lambda_{\lambda,\eta}+1)^{\frac{1}{2\gamma}}}\big)|\alpha_j|^2\notag\\
			&\le C_\gamma \exp\big({C_\gamma(\Lambda_{\lambda,\eta}+1)^{\frac{1}{2\gamma}}}\big)  \|\sum_{j=1}^{N_{\lambda,\eta}}\alpha_j\phi_{\gamma,j}\|_{L^2(\R^n)}^2
		\end{align}
		and
		\begin{align}\label{eq:doub4}
			\|e^{\rho_\eta|x|/4}|\nabla\phi|\|_2^2 &\le C_\gamma|\eta|^{\frac{2}{1+\gamma}}N_{\lambda,\eta}\sum_{j=1}^{N_{\lambda,\eta}}(\Lambda_{\lambda,\eta}+1)\exp\big({C_\gamma(\Lambda_{\lambda,\eta}+1)^{\frac{1}{2\gamma}}}\big)|\alpha_j|^2\notag\\
			&\le
			C_\gamma \exp\Big({C_\gamma\big((\Lambda_{\lambda,\eta}+1)^{\frac{1}{2\gamma}}}+\ln (|\eta|+1)\big)\Big)
			\|\sum_{j=1}^{N_{\lambda,\eta}}\alpha_j\phi_{\gamma,j}\|_{L^2(\R^n)}^2.
		\end{align}
		Since $$\big\|\sum_{j=1}^{N_{\lambda,\eta}}\alpha_j\phi_{\gamma,j}\big\|_2^2=\big\|\sum_{j=1}^{N_{\lambda,\eta}}\alpha_j\phi_{\gamma,j}^\eta\big\|_2^2=\|\phi\|_2^2,$$
		it follows from \eqref{eq:doub5}, \eqref{eq:doub3} and \eqref{eq:doub4} that
		\begin{equation*}
			\|\phi\|_{H^1(\R^n\setminus B_r)}^2\le e^{-|\eta|^{\frac{1}{1+\gamma}}r/2}\cdot e^{C_2M_{\lambda,\eta}}\|\phi\|_{L^2(\R^n)}^2,
		\end{equation*}
		with $C_2$ depending on $n,\gamma$.
		Let $r=2|\eta|^{-\frac{1}{1+\gamma}}(r'+C_2M_{\lambda,\eta})$, we deduce \eqref{equ-local-estim}.
		
		Taking $r'=\ln 2$ in \eqref{equ-local-estim} and let $r=2|\eta|^{-\frac{1}{1+\gamma}}(\ln 2+C_2M_{\lambda,\eta})$, then
		\begin{align*}
			\|\phi\|_{L^2(\R^n)}^2&=\|\phi\|_{L^2(\R^n\setminus B_r)}^2+\|\phi\|_{L^2( B_r)}^2\notag\\
			&\le \frac{1}{2}\|\phi\|_{L^2(\R^n)}^2+\|\phi\|_{L^2( B_r)}^2.
		\end{align*}
		Thus,
		\begin{equation*}
			\|\phi\|_{L^2(\R^n)}^2\le 2\|\phi\|_{L^2( B_r)}^2\le 2\|\phi\|_{L^2( B(0,\hat{C}|\eta|^{-\frac{1}{1+\gamma}}M_{\lambda,\eta}))}^2
		\end{equation*}
		with $\hat{C}$ depending only on $n,\gamma$.
		On the other hand, by using \eqref{equ-local-estim} again,
		\begin{align*}
			\|\phi\|_{H^1(\R^n)}^2&=\|\phi\|_{H^1(\R^n\setminus B_r)}^2+\|\phi\|_{H^1( B_r)}^2\\
			&\le \frac{1}{2}\|\phi\|_{H^1(\R^n)}^2+\|\phi\|_{H^1( B_r)}^2,
		\end{align*}
		which implies \eqref{equ-H1bound-L}.
	\end{proof}
	
	For $0<\gamma<1/2$, the potential $|x|^{2\gamma}$ is not Lipschitz near the origin. The next estimate localizes the spectral mass to an annulus, where the propagation of smallness argument can be applied.
	\begin{lemma}\label{lem:doub-phi-ball}
		Assume that $0<\gamma<\frac12$, $\eta\in\R^m\setminus\{0\}$, and fix
		$\varepsilon>0$ and $a>0$. There exist constants
		$C_1=C_1(n,\gamma,a,\varepsilon)>1$ and
		$C_2=C_2(n,\gamma)>1$ with the following property.
		For every $\lambda\ge\lambda_{\gamma,1}^\eta$, define
		\begin{equation}\label{equ-J-smallGam}
			\mathcal{J}_\eta(\lambda)
			:=
			\begin{cases}
				\bigl(|\eta|^{1/3}+\lambda^{1/2}+1\bigr)^{\frac{1}{2\gamma}},
				& 0<|\eta|<a,\\
				\bigl(
				|\eta|^{\frac43+2\varepsilon}
				+|\eta|^{-\frac{2-\gamma}{3(1+\gamma)}+2\varepsilon}
				\lambda^{1/2}+1
				\bigr)^{\frac{1}{2\gamma}},
				& |\eta|\ge a.
			\end{cases}
		\end{equation}
		Using $M_{\lambda,\eta}$ from \eqref{equ-M}, set
		\begin{equation*}
			\begin{aligned}
				\bar{R}_{\lambda,\eta}
				&:=2|\eta|^{-\frac{1}{1+\gamma}}
				\bigl(2C_1\mathcal{J}_\eta(\lambda)+C_2M_{\lambda,\eta}\bigr),\\
				\widetilde{R}_{\lambda,\eta}
				&:=2|\eta|^{-\frac{1}{1+\gamma}}
				\bigl(2C_1\bigl(\mathcal{J}_\eta(\lambda)+\ln(1+\lambda)\bigr)
				+C_2M_{\lambda,\eta}\bigr).
			\end{aligned}
		\end{equation*}
		Then, for every $\phi\in\operatorname{Ran}(P_\lambda(H_{\gamma,\eta}))$,
		the following estimates hold \footnote{The annulus
			$B(0,\bar{R}_{\lambda,\eta})\setminus B(0,70\sqrt n)$ is nonempty:
			$\bar{R}_{\lambda,\eta}\gtrsim (1+\lambda)^{1/(4\gamma)}$ when $0<|\eta|<a$, while $\bar{R}_{\lambda,\eta}\gtrsim (1+|\eta|)^{(2-\gamma)/(3\gamma(1+\gamma))}
			>(1+|\eta|)^{1/(2\gamma(1+\gamma))}$ when $|\eta|>a$. Choosing the implicit constant in $\bar{R}_{\lambda,\eta}$ sufficiently large ensures that the annulus is nonempty in both cases.}:
		\begin{equation}\label{equ-L2bound}
			\|\phi\|_{L^2(\R^n)}^2
			\le 2\exp\!\left\{C_1\mathcal{J}_\eta(\lambda)^{2\gamma}\right\} \|\phi\|_{L^2\bigl(B(0,\bar{R}_{\lambda,\eta}) \setminus B(0,70\sqrt{n})\bigr)}^2,
		\end{equation}
		and
		\begin{equation}\label{equ-H1bound}
			\|\phi\|_{H^1(\R^n)}^2
			\le 2\exp\!\left\{
			C_1\bigl(\mathcal{J}_\eta(\lambda)^{2\gamma} +\ln(1+\lambda)\bigr) \right\}
			\|\phi\|_{H^1\bigl(B(0,\widetilde{R}_{\lambda,\eta}) \setminus B(0,70\sqrt{n})\bigr)}^2.
		\end{equation}
	\end{lemma}

	\begin{proof}
		For any $\eps>0,a>0$, applying \cite[Proposition 3.8]{JW25} with $V(x)=|x|^{2\gamma}$, $r=|\eta|^2$ and $\sigma=0$, there exists $C>1$ depending on $n,\gamma,\eps,a$, such that
		for every$\lambda\ge \lambda_{\gamma,1}^\eta$, $\phi\in \text{Ran}(P_\lambda(H_{\gamma,\eta}))$,
		\begin{equation*}
			\|\phi\|_{L^2(\R^n)}\le e^{C\mathcal{J}_\eta(\lambda)^{2\gamma}}\|\phi\|_{L^2(B(0,1/3)^c)},
		\end{equation*}
		where $\mathcal{J}_\eta(\lambda)$ is given by  \eqref{equ-J-smallGam}.
		To replace $B(0,{1/3})^c$ by $B(0,{70\sqrt{n}})^c$, set
		$$
		b=210\sqrt{n},\quad
		(U_b\phi)(z)=b^{n/2}\phi(bz),\quad
		\eta'=b^{1+\gamma}\eta,\quad
		\lambda'=b^2\lambda.
		$$
		Since
		$
		H_{\gamma,\eta'}U_b=b^2U_bH_{\gamma,\eta},
		$
		we have
		$$
		U_b\phi\in\operatorname{Ran} \bigl(P_{\lambda'}(H_{\gamma,\eta'})\bigr),
		\quad \|U_b\phi\|_2=\|\phi\|_2, \quad
		\|U_b\phi\|_{L^2(B_{1/3}^c)} =\|\phi\|_{L^2(B(0,{70\sqrt{n}})^c)}.
		$$
		We apply \cite[Proposition 3.8 and Remark 3.9]{JW25} to $U_b\phi$, choosing the splitting point at $r=b^{2(1+\gamma)}a^2$, equivalently $|\eta'|=b^{1+\gamma}a$.
		Since $b$ is fixed, substituting $\eta'$ and $\lambda'$ changes only the constants in the two branches of $\mathcal J_\eta(\lambda)^{2\gamma}$.
		Consequently,
		\begin{equation}\label{equ-ball-compl}
			\|\phi\|_{L^2(\mathbb R^n)}^2 \le e^{C_1\mathcal J_\eta(\lambda)^{2\gamma}} \|\phi\|_{L^2(B(0,{70\sqrt{n}})^c)}^2,
		\end{equation}
		where $C_1>1$ depends only on $n,\gamma,\varepsilon,a$.

		We first prove \eqref{equ-L2bound}. For simplicity, set
		$$
		\mathscr{J}_{\lambda,\eta}=C_1\mathcal J_\eta(\lambda)^{2\gamma}.
		$$
		Applying \eqref{equ-local-estim} with
		$r'=2C_1\mathcal J_\eta(\lambda)$ and using
		$0<2\gamma<1$ and $\mathcal J_\eta(\lambda)>1$, we obtain
		\begin{equation*}
			\|\phi\|_{L^2(B(0,{\bar{R}_{\lambda,\eta}})^c)}^2 \le e^{-2C_1\mathcal J_\eta(\lambda)}\|\phi\|_2^2 \le e^{-2\mathscr{J}_{\lambda,\eta}}\|\phi\|_2^2.
		\end{equation*}
		This, combining with \eqref{equ-ball-compl}, implies
		\begin{align*}
			\|\phi\|_{L^2(B(0,{\bar{R}_{\lambda,\eta}})\setminus B(0,{70\sqrt{n}}))}^2
			&= \|\phi\|_{L^2(B(0,{70\sqrt{n}})^c)}^2 -
			\|\phi\|_{L^2(B(0,{\bar{R}_{\lambda,\eta}})^c)}^2\\
			&\ge e^{-\mathscr{J}_{\lambda,\eta}}(1-e^{-\mathscr{J}_{\lambda,\eta}})\|\phi\|_2^2.
		\end{align*}
		Since $\mathscr{J}_{\lambda,\eta}\ge\ln 2$, it follows that
		$$
		\|\phi\|_2^2 \le 2e^{C_1\mathcal J_\eta(\lambda)^{2\gamma}} \|\phi\|_{L^2(B(0,{\bar{R}_{\lambda,\eta}})\setminus B(0,{70\sqrt{n}}))}^2.
		$$
		
		We next prove \eqref{equ-H1bound}. Since $\phi\in\operatorname{Ran}(P_\lambda(H_{\gamma,\eta}))$ and the potential is nonnegative,
		$$
		\|\phi\|_{H^1(\mathbb R^n)}^2 \le (1+\lambda)\|\phi\|_2^2.
		$$
		Combining this with \eqref{equ-ball-compl}, we obtain
		\begin{equation}\label{eq:H1-out70}
			\|\phi\|_{H^1(\mathbb R^n)}^2 \le e^{C_1\big(\mathcal J_\eta(\lambda)^{2\gamma}+\ln(1+\lambda)\big)} \|\phi\|_{H^1(B(0,{70\sqrt{n}})^c)}^2.
		\end{equation}
		Set
		$$
		\tilde{\mathscr{J}}_{\lambda,\eta}= C_1 \bigl(\mathcal J_\eta(\lambda)^{2\gamma} +\ln(1+\lambda)\bigr).
		$$
		Applying \eqref{equ-local-estim} with
		$r'=2 C_1 \bigl(\mathcal J_\eta(\lambda)+\ln(1+\lambda)\bigr),$
		we have
		\begin{align*}
			\|\phi\|_{H^1(B(0,{\widetilde R_{\lambda,\eta}})^c)}^2 \le e^{-2 C_1 (\mathcal J_\eta(\lambda)+\ln(1+\lambda))} \|\phi\|_2^2 \le e^{-2\tilde{\mathscr{J}}_{\lambda,\eta}}
			\|\phi\|_{H^1(\mathbb R^n)}^2,
		\end{align*}
		where the last inequality follows from the facts that $0<2\gamma<1$ and $\mathcal{J}_\eta(\lambda)>1$.
		This, combining with \eqref{eq:H1-out70}, implies
		\begin{align*}
			\|\phi\|_{H^1(B(0,{\widetilde R_{\lambda,\eta}}) \setminus B(0,{70\sqrt{n}}))}^2 \ge e^{-\tilde{\mathscr{J}}_{\lambda,\eta}}(1-e^{-\tilde{\mathscr{J}}_{\lambda,\eta}}) \|\phi\|_{H^1(\mathbb R^n)}^2 \ge \frac12 e^{-\tilde{\mathscr{J}}_{\lambda,\eta}} \|\phi\|_{H^1(\mathbb R^n)}^2.
		\end{align*}
		Thus,
		$$
		\|\phi\|_{H^1(\mathbb R^n)}^2 \le 2e^{ C_1 (\mathcal J_\eta(\lambda)^{2\gamma}+\ln(1+\lambda))} \|\phi\|_{H^1(B(0,{\widetilde R_{\lambda,\eta}}) \setminus B(0,{70\sqrt{n}}))}^2.
		$$
		Therefore the proof of Lemma \ref{lem:doub-phi-ball} is complete.
	\end{proof}

	\subsubsection{Propagation of smallness}
	The localization estimates above concentrate the spectral mass to the regions where propagation is needed. We now derive a propagation of smallness estimate that transfers smallness from the thick observation set to these regions, using an elliptic extension of each function in the spectral subspace.
	
	Let $\phi=\sum_{\lambda_{\gamma,j}^\eta\le\lambda}\alpha_j\phi_{\gamma,j}^\eta\in \text{Ran}(P_{\lambda}(H_{\gamma,\eta}))$, where the eigenpairs $(\phi_{\gamma,j}^\eta,\lambda_{\gamma,j}^\eta)$ satisfy \eqref{equ-diag-H-eta}.
	For convenience, we introduce
	$$\mathcal{S}_a(s)={\sinh\left(\sqrt{a}s\right)}/{\sqrt{a}},\qquad a>0,\, s\in\R.$$
	Set
	\begin{equation}\label{eq:def-hat-Phi}
		\hat{\Phi}=\hat{\Phi}_\phi(x,s):=\sum_{\lambda_{\gamma,j}^\eta\le\lambda}\alpha_j\phi_{\gamma,j}^\eta\mathcal{S}_{\lambda_{\gamma,j}^\eta}(s).
	\end{equation}
	%
	By construction, $\hat{\Phi}$ is the elliptic extension of $\phi$ satisfying
	\begin{equation}\label{eq:hat-Phi}
		\begin{cases}
			-\Delta\hat{\Phi}+|\eta|^2|x|^{2\gamma}\hat{\Phi}=0\quad\mbox{in }\R^{n+1},\\
			\partial_{s}\hat{\Phi}(x,0)=\phi(x),\,\hat{\Phi}(x,0)=0\quad\mbox{in }\R^n,
		\end{cases}
	\end{equation}
	where  $\Delta\hat{\Phi}=\sum_{i=1}^{n}\partial_{x_i}^2\hat{\Phi}+\partial_{s}^2\hat{\Phi}$.
	
	Let $V_\eta(x)=|\eta|^2|x|^{2\gamma}$. For a bounded domain $\omega\subset\R^n$ such that $V_\eta\in W^{1,\infty}(\omega)$, define
	\begin{equation*}
		\|V_\eta\|_{C^{0,1}(\omega)} :=\|V_\eta\|_{L^\infty(\omega)} +\||\nabla V_\eta|\|_{L^\infty(\omega)}.
	\end{equation*}
	Set
	$$
	\tau_0=\bigl(5+4\|V_\eta\|_{C^{0,1}(\omega)}\bigr)^{1/2}, \qquad I_\rho=(-\rho,\rho),\quad \rho>0.
	$$
	Suppose that $\Omega\subset\R^n$ is a $(\delta,1)$ thick set. For $j\in\mathbb Z^n$, set $Q_j=\Lambda_2(j)$ and $\Omega_j:=\Omega\cap\Lambda_1(j)$, and take $\omega=Q_j$. Thus, $\Omega_j$ is the portion of the observation set in the unit cube $\Lambda_1(j)\subset Q_j$; the sets $\Omega_j$ are pairwise disjoint up to their boundaries.
	
	\begin{lemma}\label{lem:propagation-smallness}
        Let $\gamma>0$, $\eta\ne0$, $\lambda>0$, $\phi\in\operatorname{Ran}\bigl(P_\lambda(H_{\gamma,\eta})\bigr)$, and $j\in\mathbb{Z}^n$.
        If $0<\gamma<\tfrac12$, assume additionally that $4Q_j\subset\mathbb{R}^n\setminus B(0,1)$. Then the elliptic extension $\widehat{\Phi}$ of $\phi$ satisfies
		\begin{equation}\label{equ-three-sphere}
			\|\hat\Phi\|_{H^1(Q_j\times I_\rho)} \le C\bigl(5+8\|V_\eta\|_{C^{0,1}(\omega)}\bigr)e^{10\tau_0\rho}\delta^{-\alpha/2}\|\phi\|_{L^2(\Omega_j)}^\alpha\|\hat\Phi\|_{H^1(4Q_j\times I_{4\rho})}^{1-\alpha},
		\end{equation}
    where $C>0$ and $\alpha\in(0,1)$ depend only on $n,\delta,\gamma,\rho$.    
	\end{lemma}
	
	\begin{remark}\label{rmk:three-ball}
    If $\gamma\ge1/2$, then $V_\eta(x)=|\eta|^2|x|^{2\gamma}$ is locally Lipschitz on $\mathbb R^n$, so no additional restriction on $j$ is needed. If $0<\gamma<1/2$, its gradient is unbounded near the origin. The condition $4Q_j\subset\mathbb R^n\setminus B(0,1)$ is a sufficient condition ensuring the Lipschitz regularity on $4Q_j$ and the uniform coefficient bounds used in the proof.
	\end{remark}
	
	\begin{proof}
		We follow the argument in \cite{Zhu24}. For completeness, we outline the proof and specify the dependence of the constants. Fix a cube $Q_j$ satisfying the assumptions of the lemma. Clearly,
		\begin{equation}\label{eq:potential-enlarged-cube}
			\|V_\eta\|_{C^{0,1}(4Q_j)}\le C_{n,\gamma}\|V_\eta\|_{C^{0,1}(Q_j)}.
		\end{equation}
		Let
		\begin{equation*}
				\Phi(x,s,t,y)= e^{\tau_0 y}e^{i2(\|V_\eta\|_{C^{0,1}(\omega)}+1)^{1/2}t} \hat{\Phi}(x,s),\quad \mathcal{Q}_j=Q_j\times I_\rho\times I_\rho\times I_\rho.
		\end{equation*}
		By \eqref{eq:hat-Phi}, we have
		\begin{equation*}
			-\Delta_{x,s}\Phi - \frac{V_\eta(x)+5+4\|V_\eta\|_{C^{0,1}(\omega)}}{4(\|V_\eta\|_{C^{0,1}(\omega)}+1)}\partial_{t}^2\Phi-\partial_{y}^2\Phi=0\quad \text{in } \mathbb{R}^{n+3}.
		\end{equation*}
		Equivalently, $-\operatorname{div}(A\nabla\Phi)=0$, where
		\begin{equation}\label{eq:matrix-A}
			A=\operatorname{diag}\left(I_n,1,1+\frac{V_\eta(x)+1}{4(\|V_\eta\|_{C^{0,1}(\omega)}+1)},1\right).
		\end{equation}
		By \eqref{eq:potential-enlarged-cube}, the matrix $A$ satisfies
		\begin{equation*}
			|\xi|^2\le\langle A(X)\xi,\xi\rangle\le C_{n,\gamma}|\xi|^2,\qquad|A(X)-A(Y)|\le C_{n,\gamma}|X-Y|,\quad X,Y\in4\mathcal Q_j,
		\end{equation*}
		where $C_{n,\gamma}$ depends only on $n$ and $\gamma$.
		As in \cite[(2.31)--(2.32)]{Zhu24}, direct calculations show that
		\begin{equation}\label{eq:Zhu-2.31}
				\|\nabla\Phi\|_{L^2(4\mathcal Q_j)}^2\le C\rho^2(5+8\|V_\eta\|_{C^{0,1}(\omega)})e^{8\tau_0\rho}\|\hat\Phi\|_{H^1(4Q_j\times I_{4\rho})}^2,
		\end{equation}
		and
		\begin{equation}\label{eq:Zhu-2.32}
			\|\nabla\Phi\|_{L^2(\mathcal Q_j)}^2\ge C\rho^2 e^{-2\tau_0\rho}\|\hat{\Phi}\|_{H^1(Q_j\times(-\rho,\rho))}^2.
		\end{equation}
		On the other hand, let $\mathcal{E}_j=\Omega_j\times\{0\}\times(-\frac{\rho}{2},\frac{\rho}{2})\times(-\frac{\rho}{2},\frac{\rho}{2})$. By rescaling estimates, we can identify the cubes $\mathcal Q_j$ with the standard cube $Q^{n+3}$ of side length $2$, and $\mathcal E_j$ with the set $E$ in \cite[Corollary 2]{Zhu24}. 
		For fixed $\rho$, these rescalings preserve the ellipticity and Lipschitz bounds for $A$ in \eqref{eq:matrix-A}, with constants depending only on $n,\gamma,\rho$.
		Let
		\begin{equation*}
			\hat\phi(x,t,y)=e^{\tau_0 y}e^{i2(\|V_\eta\|_{C^{0,1}(\omega)}+1)^{1/2}t}\phi(x).
		\end{equation*}
		Since $\nabla\Phi(x,0,t,y)=(0,\hat\phi(x,t,y),0,0)$, we have
		\begin{equation}\label{eq:Zhu-2.54}
			\|\hat{\phi}\|_{L^2(\mathcal{E}_j)}^2\le C\rho^2e^{\tau_0\rho}\|\phi\|_{L^2(\Omega_j)}^2.
		\end{equation}
		By \cite[Corollary 2]{Zhu24},
		\begin{equation}\label{eq:Zhu-2.53}
			\|\nabla\Phi\|_{L^2(\mathcal Q_j)}\le C\left(\frac{|\mathcal Q_j|}{|\mathcal{E}_j|}\right)^{\frac{\alpha}{2}}\|\hat{\phi}\|_{L^2(\mathcal{E}_j)}^\alpha\|\nabla\Phi\|_{L^2(4\mathcal Q_j)}^{1-\alpha},
		\end{equation}
		where the exponent $\alpha\in(0,1)$ can be chosen to depend only on $n,\delta,\gamma,\rho$ according to \cite[Remark 3]{Zhu24}.
		Since $\Omega$ is a $(\delta,1)$-thick set, we have
		\begin{equation*}
			\frac{|\mathcal Q_j|}{|\mathcal{E}_j|}\le C_{n,\rho} \delta^{-1}.
		\end{equation*}
		This, combining with \eqref{eq:Zhu-2.54}, \eqref{eq:Zhu-2.53}, implies
		\begin{equation}\label{eq:Zhu-2.53-54}
			\|\nabla\Phi\|_{L^2(\mathcal Q_j)}\le C\rho^\alpha e^{{\tau_0 \rho\alpha}/{2}}\delta^{-{\alpha}/{2}}\|\phi\|_{L^2(\Omega_j)}^\alpha\|\nabla\Phi\|_{L^2(4\mathcal Q_j)}^{1-\alpha}.
		\end{equation}
		Therefore, by \eqref{eq:Zhu-2.31}, \eqref{eq:Zhu-2.32} and \eqref{eq:Zhu-2.53-54}, we derive
		\begin{equation*}
				\|\hat{\Phi}\|_{H^1(Q_j\times I_\rho)}\le C\bigl(5+8\|V_\eta\|_{C^{0,1}(\omega)}\bigr)e^{10\tau_0\rho}\delta^{-\frac{\alpha}{2}}\|\phi\|_{L^2(\Omega_j)}^\alpha\|\hat{\Phi}\|_{H^1(4Q_j\times I_{4\rho})}^{1-\alpha},
		\end{equation*}
		where the constants $C>0$ and $\alpha\in(0,1)$ depend only on $n,\delta,\gamma,\rho$.
	\end{proof}

	\subsection{Uniform spectral inequality for $H_{\gamma,\eta}$}
	
	We combine the localization estimates with the propagation of smallness to obtain a spectral inequality that is uniform in the nonzero parameter $\eta\in\R^m$.
	\begin{theorem}\label{thm-spec-H-eta}
		Let $\gamma>0$, $L>0$, $\eta\in\R^m\setminus\{0\}$, and let $\Omega\subset\R^n$ be a $(\delta,L)$ thick set. For every $\eps>0$, there exists a constant $C=C(n,\gamma,\delta,L,\eps)>0$ such that, for every $\lambda>0$ and every $\phi\in\operatorname{Ran}(P_\lambda(H_{\gamma,\eta}))$,
		\begin{equation}\label{equ-spec-H-eta}
			\|\phi\|_{L^2(\R^n)}\le
			\begin{cases}
				C\exp\!\left\{C\bigl(1+\lambda^{\frac{1+\gamma}{2}+\eps}\bigr)\right\} \|\phi\|_{L^2(\Omega)},
				& \gamma\ge\frac12,\\[4pt]
				C\exp\!\left\{C\bigl(1+\lambda^{\frac{5+2\gamma}{6}+\eps}\bigr)\right\}\|\phi\|_{L^2(\Omega)},
				& 0<\gamma<\frac12.
			\end{cases}
		\end{equation}
	\end{theorem}
	
	\begin{remark}
		Spectral inequalities for Schr\"odinger operators with growing potentials have been studied extensively; see, for example, \cite{AS24,BJP21,DSV24,MZ25,Zhu24,ZZG25}.
		Zhu and Zhuge \cite{ZZG25} considered general potentials with power growth and sets $\Omega$ satisfying
		\[
		\forall j\in\mathbb Z^n,\quad\exists z_j\in\Lambda_1(j)\quad\text{such that}\quad B_{\delta^{1+|j|^\sigma}}(z_j)\subset\Omega\cap\Lambda_1(j),
		\]
		where $\delta\in(0,1/2)$ and $\sigma\ge0$. When $\sigma=0$, this condition requires a ball of fixed radius in each unit cube and corresponds to the usual equidistribution condition. They established the spectral inequality
		\[
		\|\phi\|_{L^2(\mathbb R^n)}\le\big(\frac{1}{\delta}\big)^{C\lambda^{\frac{\sigma}{\beta_1}+\frac{\beta_2}{2\beta_1}}}\|\phi\|_{L^2(\Omega)},\qquad\phi\in\operatorname{Ran}(P_\lambda(-\Delta+V)),
		\]
		where $\beta_1,\beta_2$ are the growth exponents in their
		assumptions on $V$.
		For $V(x)=|x|^{2\gamma}$ with $\gamma>0$,
		one can take $\beta_1=\beta_2=2\gamma$.
		Under a different set of assumptions on the potential, which also includes $V(x)=|x|^{2\gamma}$ for every $\gamma>0$, Zhu \cite{Zhu24} obtained spectral inequalities with the same exponent for measurable sets satisfying
		\[
		|\Omega\cap\Lambda_1(j)|\ge \delta^{1+|j|^\sigma},\qquad j\in\mathbb Z^n,\quad 0\le\sigma<1.
		\]
		For $\sigma=0$, this condition implies that $\Omega$ is thick.
		
		Based on \cite{ZZG25}, Jaming and Wang \cite{JW25} derived spectral inequalities for $-\Delta+rV$ with explicit dependence on $r>0$.
		Taking $r=|\eta|^2$ and $V(x)=|x|^{2\gamma}$ yields spectral inequalities for $H_{\gamma,\eta}$ with explicit dependence on $\eta$, in particular on equidistributed sets.
		Inspired by \cite{Zhu24}, Theorem~2.3 establishes a spectral inequality for $H_{\gamma,\eta}$ on measurable thick sets, with constants independent of $\eta\in\mathbb{R}^m\setminus\{0\}$.
	\end{remark}

	
	\begin{proof}
        For simplicity of notation, we present the proof only for $L=1$. The same argument applies to any fixed $L>0$, with the constants depending on $L$.
		We adapt the localization and propagation of smallness argument in the proof of \cite[Theorem~1, p.~5699]{Zhu24} to $V_\eta(x)=|\eta|^2|x|^{2\gamma}$. The four cases below are distinguished according to whether $\gamma\ge\frac12$ or $0<\gamma<\frac12$, and whether $0<|\eta|<1$ or $|\eta|\ge1$. In each case, Step 1 establishes the localization and doubling estimates, while Step 2 applies \eqref{equ-three-sphere} via a covering argument to obtain the global estimate.
		Throughout this proof, $\hat{\Phi}$ denotes the elliptic extension of $\phi$ defined in \eqref{eq:def-hat-Phi}. In particular, for every $s$, $\hat{\Phi}(\cdot,s)$ and $\partial_s\hat{\Phi}(\cdot,s)$ belong to $\operatorname{Ran}(P_\lambda(H_{\gamma,\eta}))$, and $\partial_s\hat{\Phi}(x,0)=\phi(x)$. We may assume that $\lambda\ge\lambda_{\gamma,1}^{\eta}$, since otherwise $\operatorname{Ran}(P_\lambda(H_{\gamma,\eta}))=\{0\}$ and the desired inequality is trivial. Under this assumption, the identity $\lambda_{\gamma,1}^\eta=|\eta|^{{2}/({1+\gamma})}\lambda_{\gamma,1}$ gives
		\begin{equation}\label{eq:eta-bound}
			|\eta|\le \left(\frac{\lambda}{\lambda_{\gamma,1}}\right)^{\frac{1+\gamma}{2}} \le C\lambda^{\frac{1+\gamma}{2}},
		\end{equation}
		where $C$ depends only on $n$ and $\gamma$.
		
		First, we compare the $L^2$-norm of $\phi$ and the $H^1$-norm of $\hat{\Phi}$ (see \cite{JL99, ZZG25} for details). More precisely, for any $\lambda \ge \lambda_{\gamma,1}^\eta$,
		\begin{equation}\label{equ-comp-phi-Phi}
			2\rho\|\phi\|_{L^2(\R^n)}^2\le \|\hat{\Phi}\|_{H^1(\R^n\times I_\rho)}^2\le 2\rho(1+\frac{\rho^2}{3}(1+\lambda))e^{2\rho\sqrt{\lambda}}\|\phi\|_{L^2(\R^n)}^2,
		\end{equation}
		where $\rho>0$ is the same as in \eqref{equ-three-sphere}. In what follows, we fix $\eps>0$ and the choice of $\rho$. All constants below are independent of $\lambda$ and $\eta$, but may depend on $n,\gamma,\delta,\eps,\rho$.

		\noindent\emph{Case 1. $\gamma\ge 1/2$ and $0<|\eta|<1$.}
		
		\emph{Step 1. Localization and doubling estimates.}
		We claim
		\begin{equation}\label{equ-doub-Phi-HL}
			\|\hat{\Phi}\|_{H^1(\R^n\times I_{\rho/2})}^2\le C \big(1+\lambda\big)\|\hat{\Phi}\|_{L^2\big(B(0,2R_{\lambda,\eta})\times I_\rho\big)}^2,
		\end{equation}
		and
		\begin{equation}\label{equ-doub-Phi-HH}
			\|\hat{\Phi}\|_{H^1(\R^n\times I_{4\rho})}^2\le C  e^{9\rho\sqrt{\lambda}}\|\hat{\Phi}\|_{H^1\big(B(0,R_{\lambda,\eta})\times I_{\rho/2}\big)}^2,
		\end{equation}
		where $R_{\lambda,\eta}$ denotes the radius defined in \eqref{equ-R-lambda-eta}.
		
		Indeed, since $|\eta|<1$, we have
		\begin{equation}\label{eq:V-L-infty1}
			\|V_\eta\|_{L^\infty\big(B(0,2R_{\lambda,\eta})\big)}\le C \big(1+\lambda+|\eta|^{\frac{2}{1+\gamma}}(\ln(|\eta|+1))^{2\gamma}\big)\le C (1+\lambda).
		\end{equation}
		The estimates follow by adapting the proof of Lemma~3 in \cite{Zhu24}; we sketch the main steps below.
		
		For every $s$, both $\hat{\Phi}(\cdot,s)$ and
		$\partial_s\hat{\Phi}(\cdot,s)$ belong to
		$\operatorname{Ran}(P_\lambda(H_{\gamma,\eta}))$. Apply the $H^1$ localization estimate \eqref{equ-H1bound-L} to $\hat{\Phi}(\cdot,s)$
		and the $L^2$ localization estimate \eqref{equ-L2bound-L} to $\partial_s\hat{\Phi}(\cdot,s)$,
		respectively, and integrate in $s$. This gives
		\begin{equation}\label{eq:Phi-HL1}
			\|\hat{\Phi}\|_{H^1(\R^n\times I_{\rho/2})}^2 \le C\|\hat{\Phi}\|_{H^1(B(0,R_{\lambda,\eta})\times I_{\rho/2})}^2.
		\end{equation}
		On the other hand, the Caccioppoli inequality for $-\Delta\hat{\Phi}+V_\eta\hat{\Phi}=0$, together with \eqref{eq:V-L-infty1}, yields
		\begin{equation*}
			\|\hat{\Phi}\|_{H^1(B(0,R_{\lambda,\eta})\times I_{\rho/2})}^2 \le C(1+\lambda) \|\hat{\Phi}\|_{L^2(B(0,2R_{\lambda,\eta})\times I_\rho)}^2.
		\end{equation*}
		These two inequalities prove the first estimate. For the second, \eqref{equ-comp-phi-Phi} gives
		\begin{equation}\label{eq:Phi-H1-small-rho}
			\|\hat{\Phi}\|_{H^1(\R^n\times I_{4\rho})}^2\le 8(1+\frac{16\rho^2}{3}(1+\lambda))e^{8\rho\sqrt{\lambda}}\|\hat{\Phi}\|_{H^1(\R^n\times I_{\rho/2})}^2.
		\end{equation}
		Then, \eqref{equ-doub-Phi-HH} follows from \eqref{eq:Phi-HL1} and \eqref{eq:Phi-H1-small-rho}.

		\emph{Step 2. Uniform spectral inequality for $H_{\gamma,\eta}$.}
		We take the finite family of cubes $Q_j=\Lambda_2(j)$ that meet $B(0,2R_{\lambda,\eta})$. This covering has finite overlap depending only on $n$, and
		\begin{equation}\label{equ-cover1}
			B(0,2R_{\lambda,\eta})\times(-\rho,\rho)\subset\bigcup_j Q_j\times(-\rho,\rho),
		\end{equation}
		Moreover, the dilated cubes satisfy
		\begin{equation}\label{equ-cover2}
			\bigcup_j 4Q_j\times(-4\rho,4\rho)\subset
			B(0,8R_{\lambda,\eta})\times(-4\rho,4\rho).
		\end{equation}
		Indeed, for $0<|\eta|<1$, we have $|\eta|^{-\frac{1}{1+\gamma}}M_{\lambda,\eta}\gtrsim(1+\lambda)^{\frac{1}{2\gamma}}\ge1$. Thus, after increasing the fixed constant $\hat C$ in \eqref{equ-R-lambda-eta} if necessary, \eqref{equ-cover2} holds.
		
		Since $\omega=Q_j\subset B(0,8R_{\lambda,\eta})$, we have
		$$\|V_\eta\|_{C^{0,1}(\omega)}\le C\big(1+\lambda+|\eta|^{\frac{2}{1+\gamma}}(\ln(|\eta|+1))^{2\gamma}\big)\le C(1+\lambda).$$
		It follows from \eqref{equ-three-sphere} that
		\begin{equation*}
			\|\hat{\Phi}\|_{H^1(Q_j\times I_\rho)}\le Ce^{C (1+\lambda^{1/2})}\|\phi\|_{L^2(\Omega_j)}^\alpha\|\hat{\Phi}\|_{H^1(4Q_j\times I_{4\rho})}^{1-\alpha}.
		\end{equation*}
		Summing the local estimate over the covering cubes gives the following chain; the justification for each inequality is stated immediately below it.
		\begin{align*}
			\|\hat{\Phi}\|_{H^1(B(0,2R_{\lambda,\eta})\times I_\rho)}^2&\le C\sum_j \|\hat{\Phi}\|_{H^1(Q_j\times I_\rho)}^2\\
			&\le Ce^{C (1+\lambda^{1/2})}\sum_j\big(\|\phi\|_{L^2(\Omega_j)}^{2\alpha} \|\hat{\Phi}\|_{H^1(4Q_j\times I_{4\rho})}^{2(1-\alpha)}\big)\\
			&\le Ce^{C (1+\lambda^{1/2})}\big(\sum_j \|\phi\|_{L^2(\Omega_j)}^2\big)^\alpha\big(\sum_j \|\hat{\Phi}\|_{H^1(4Q_j\times I_{4\rho})}^2\big)^{1-\alpha}\\
			&\le Ce^{C (1+\lambda^{1/2})}\|\phi\|_{L^2(\Omega)}^{2\alpha} \|\hat{\Phi}\|_{H^1(B(0,8R_{\lambda,\eta})\times I_{4\rho})}^{2(1-\alpha)}\\
			&\le Ce^{C (1+\lambda^{1/2})}\|\phi\|_{L^2(\Omega)}^{2\alpha} \|\hat{\Phi}\|_{H^1(B(0,R_{\lambda,\eta})\times I_{\rho/2})}^{2(1-\alpha)}.
		\end{align*}
		Here, the first inequality follows from \eqref{equ-cover1}, the second from \eqref{equ-three-sphere}, and the third from Hölder's inequality. The fourth relies on \eqref{equ-cover2}, combined with the disjointness of $\Omega_j$ and the finite overlap of $4Q_j$. The last follows from \eqref{equ-doub-Phi-HH}. Therefore, we conclude that
		\begin{equation}\label{equ-spec-Pphi}
			\|\hat{\Phi}\|_{H^1(B(0,2R_{\lambda,\eta})\times I_\rho)}\le Ce^{C (1+\lambda^{1/2})} \|\phi\|_{L^2(\Omega)},
		\end{equation}
		where $C$ depends on $n,\gamma,\delta,\rho$. Further, by \eqref{equ-doub-Phi-HL} and \eqref{equ-spec-Pphi}, we obtain
		\begin{align*}
			\|\hat{\Phi}\|_{H^1(\R^n\times I_{\rho/2})}\le Ce^{C (1+\lambda^{1/2})} \|\phi\|_{L^2(\Omega)}.
		\end{align*}
		At last, applying the first inequality in \eqref{equ-comp-phi-Phi}, we arrive at
		\begin{equation}\label{eq:spec-case1}
			\|\phi\|_{L^2(\R^n)}\le Ce^{C (1+\lambda^{1/2})} \|\phi\|_{L^2(\Omega)}.
		\end{equation}
		
		\noindent\emph{Case 2. $\gamma\ge 1/2$ and $|\eta|\ge 1$.}
		
		For $|\eta|\ge1$, set $A_{\lambda,\eta}:=\hat{C}M_{\lambda,\eta}$\footnote{We enlarge $R_{\lambda,\eta}$ to $A_{\lambda,\eta}$ to ensure a uniform positive lower bound on the radius, so that the covering argument in Case~1 remains valid.}, where $M_{\lambda,\eta}$ is defined in \eqref{equ-M}. 
		By Lemma~\ref{lem:doub-phi} and the fact that  $R_{\lambda,\eta}\le A_{\lambda,\eta}$ for $|\eta|\ge1$, we obtain, for any $\lambda\ge\lambda_{\gamma,1}^\eta$ and every $\phi\in \text{Ran}(P_\lambda(H_{\gamma,\eta}))$,
		\begin{equation}\label{eq:phi-doub-L2}
			\|\phi\|_{L^2(\R^n)}^2\le 2 \|\phi\|_{L^2(B(0,A_{\lambda,\eta}))}^2,
		\end{equation}
		and
		\begin{equation}\label{eq:phi-doub-H2}
			\|\phi\|_{H^1(\R^n)}^2\le 2\|\phi\|_{H^1(B(0,A_{\lambda,\eta}))}^2.
		\end{equation}
		Thanks to \eqref{eq:eta-bound}, the potential estimate on the enlarged ball becomes
		\begin{equation}\label{eq:V-bound-2}
			\|V_\eta\|_{C^{0,1}(B(0,8A_{\lambda,\eta}))} \le C\left(1+\lambda|\eta|^{\frac{2\gamma}{1+\gamma}}+|\eta|^2(\ln(1+|\eta|))^{2\gamma}\right) \le C(1+\lambda^{1+\gamma+2\eps}).
		\end{equation}
		This also bounds $\|V_\eta\|_{L^\infty(B(0,2A_{\lambda,\eta}))}$. 
		
		Using \eqref{eq:phi-doub-L2}, \eqref{eq:phi-doub-H2}, and \eqref{eq:V-bound-2}, we argue as in Step 1 of Case 1 to obtain \eqref{equ-doub-Phi-HL} and \eqref{equ-doub-Phi-HH}, with $R_{\lambda,\eta}$ replaced by $A_{\lambda,\eta}$ and the factor $1+\lambda$ in \eqref{equ-doub-Phi-HL} replaced by $1+\lambda^{1+\gamma+2\eps}$.
		Then, using \eqref{eq:V-bound-2} and following Step 2 of Case 1, we obtain
		\begin{equation*}
			\|\phi\|_{L^2(\R^n)} \le C\exp\left(C\left(1+\lambda^{\frac{1+\gamma}{2}+\eps}\right)\right)\|\phi\|_{L^2(\Omega)}.
		\end{equation*}
	
		This completes the proof of \eqref{equ-spec-H-eta} for $\gamma\ge 1/2$.
		
		For the remaining cases, $0<\gamma<1/2$. These cases differ essentially from Cases 1 and 2: the potential $|x|^{2\gamma}$ is not Lipschitz at the origin, so so the propagation of smallness argument used there cannot be applied on regions containing the origin. We instead use the annular localization in Lemma~\ref{lem:doub-phi-ball}; accordingly, the propagation of smallness argument and the cube covering below are carried out on annuli at a fixed positive distance from the origin.
		Choose
		$$
		0<\eps_0< \min\left\{ \frac{2\eps}{1+\gamma}, \frac{1}{1+\gamma}-\frac23 \right\}.
		$$
		In \eqref{equ-J-smallGam}, we take $a=1$ and denote the auxiliary parameter $\varepsilon$ by $\varepsilon_0$.

		\noindent\emph{Case 3. $0<\gamma<1/2$ and $0<|\eta|<1$.}

		\emph{Step 1. Doubling type estimate for $\hat{\Phi}$.}
		Unlike the ball estimates in Cases 1 and 2, the localization estimate used here is the annular estimate of Lemma~\ref{lem:doub-phi-ball}: it controls the global norms by the norms on $B(0,R)\setminus B(0,70\sqrt n)$, where the potential is Lipschitz. Applying it to $\hat{\Phi}(\cdot,s)$ and $\partial_s\hat{\Phi}(\cdot,s)$ for each fixed $s$ gives the annular bound below; the subsequent covering is chosen in the same annular region.
		We claim there exists $C>0$ (independent of $\lambda$ and $\eta$) such that
		\begin{equation}\label{equ-doub-Phi-HL-ball}
			\|\hat{\Phi}\|_{H^1(\R^n\times I_{\rho/2})}^2\le Ce^{C(1+\lambda^{1/2})}\|\hat{\Phi}\|_{L^2\big((B(0,2\hat{C}\hat{R}_{\lambda,\eta})\setminus B_{50\sqrt{n}})\times I_\rho\big)}^2,
		\end{equation}
		and
		\begin{equation}\label{equ-doub-Phi-HH-ball}
			\|\hat{\Phi}\|_{H^1(\R^n\times I_{4\rho})}^2\le Ce^{C(1+\lambda^{1/2})}\|\hat{\Phi}\|_{H^1\big((B(0,\hat{C}\hat{R}_{\lambda,\eta})\setminus B_{70\sqrt{n}})\times I_{\rho/2}\big)}^2,
		\end{equation}
		where $\hat{R}_{\lambda,\eta}=|\eta|^{-\frac{1}{1+\gamma}}(\mathcal{J}_\eta(\lambda)+\ln(1+\lambda)+M_{\lambda,\eta})$ with $\mathcal{J}_\eta(\lambda),\,M_{\lambda,\eta}$ given by \eqref{equ-J-smallGam} and \eqref{equ-M}, respectively.
		
		Indeed, since $\hat{\Phi}(\cdot,s)$ and $\partial_s\hat{\Phi}(\cdot,s)$ belong to $\operatorname{Ran}(P_\lambda(H_{\gamma,\eta}))$ for every $s$, applying Lemma~\ref{lem:doub-phi-ball} to these functions and integrating in $s$, together with $0<|\eta|<1$, gives
		\begin{equation}\label{eq4}
			\|\hat{\Phi}\|_{H^1(\R^n\times I_{\rho/2})}^2
			\le Ce^{C(1+\lambda^{1/2})}
			\|\hat{\Phi}\|_{H^1\big((B(0,\hat C\hat R_{\lambda,\eta})\setminus B_{70\sqrt n})\times I_{\rho/2}\big)}^2.
		\end{equation}
		Using $0<|\eta|<1$ again, by calculation
		\begin{align*}
			\|V_\eta\|_{L^\infty(B(0,2\hat{C}\hat{R}_{\lambda,\eta})\setminus B_{50\sqrt{n}})}&\le C|\eta|^{\frac{2}{1+\gamma}}\big(1+|\eta|^{1/3}+\lambda^{1/2}+\lambda|\eta|^{-\frac{2}{1+\gamma}}+(\ln(1+\lambda+|\eta|)^{2\gamma})\big)\\
			&\le  C(1+\lambda).
		\end{align*}
		Combining this with \eqref{eq4}, repeating the strategy in \cite[Lemma 3]{Zhu24}, we obtain \eqref{equ-doub-Phi-HL-ball}.
		On the other hand, by \eqref{eq:Phi-H1-small-rho} and \eqref{eq4}, we obtain \eqref{equ-doub-Phi-HH-ball}.
		
		\emph{Step 2. Uniform spectral inequality for $H_{\gamma,\eta}$.}
		Increasing $\hat C$ if necessary, we may assume that $\hat C\hat R_{\lambda,\eta}>70\sqrt n$. We choose a finite family of cubes $Q_j$ satisfying
		\begin{equation}\label{equ-cover1-ball}
			(B(0,2\hat{C}\hat{R}_{\lambda,\eta})\setminus B_{50\sqrt{n}})\times(-\rho,\rho)\subset\bigcup_j Q_j\times(-\rho,\rho),
		\end{equation}
		and
		\begin{equation}\label{equ-cover2-ball}
			\bigcup_j 4Q_j\times(-4\rho,4\rho)\subset
			(B(0,8\hat{C}\hat{R}_{\lambda,\eta})\setminus B_{30\sqrt{n}})\times(-4\rho,4\rho),
		\end{equation}
		with overlap bounded by a constant depending only on $n$.

		On this annulus, $\|V_\eta\|_{C^{0,1}(Q_j)}\le C(1+\lambda)$.
		Moreover, $4Q_j\subset\R^n\setminus B_1$, so the local estimate \eqref{equ-three-sphere} gives
		\begin{equation*}
			\|\hat{\Phi}\|_{H^1(Q_j\times I_{\rho})}\le Ce^{C (1+\lambda^{1/2})}\|\phi\|_{L^2(\Omega_j)}^\alpha\|\hat{\Phi}\|_{H^1(4Q_j\times I_{4\rho})}^{1-\alpha}.
		\end{equation*}
		Now repeat the bounded-overlap and H\"older steps from Case 1, using the annular coverings \eqref{equ-cover1-ball}--\eqref{equ-cover2-ball} and \eqref{equ-doub-Phi-HH-ball} for the final localization. This gives
		\begin{equation*}
			\begin{aligned}
				\|\hat{\Phi}\|_{H^1((B(0,2\hat{C}\hat{R}_{\lambda,\eta})\setminus B_{50\sqrt{n}})\times I_\rho)}^2
				&\le Ce^{C(1+\lambda^{1/2})}\|\phi\|_{L^2(\Omega)}^{2\alpha}\\
				&\quad\times\|\hat{\Phi}\|_{H^1((B(0,\hat{C}\hat{R}_{\lambda,\eta})\setminus B_{70\sqrt{n}})\times I_{\rho/2})}^{2(1-\alpha)},
			\end{aligned}
		\end{equation*}
		which implies
		\begin{equation}\label{equ-spec-Pphi-ball}
			\|\hat{\Phi}\|_{H^1((B(0,2\hat{C}\hat{R}_{\lambda,\eta})\setminus B_{50\sqrt{n}})\times I_{\rho})}\le Ce^{C(1+\lambda^{1/2})} \|\phi\|_{L^2(\Omega)}.
		\end{equation}
		Further, by \eqref{equ-doub-Phi-HL-ball} and \eqref{equ-spec-Pphi-ball}, we obtain
		\begin{align*}
			\|\hat{\Phi}\|_{H^1(\R^n\times I_{\rho/2})}\le Ce^{C(1+\lambda^{1/2})} \|\phi\|_{L^2(\Omega)}.
		\end{align*}
		Finally, applying the first inequality in \eqref{equ-comp-phi-Phi} again, we derive
		\begin{equation*}
			\|\phi\|_{L^2(\R^n)}\le Ce^{C (1+\lambda^{1/2})} \|\phi\|_{L^2(\Omega)}.
		\end{equation*}

		\emph{Case 4. $0<\gamma<1/2$ and $|\eta|\ge1$.}
		This case uses the same annular localization from Lemma~\ref{lem:doub-phi-ball} as Case 3; in particular, the cubes in the propagation step remain away from the origin. The difference lies in the large-$\vert{}\eta\vert{}$ bound on $\mathcal J_\eta(\lambda)$, which in turn affects the powers of $\lambda$ in the annular localization and propagation of smallness estimates.
		We set
		$
		\kappa=\big(\frac23+\eps_0\big)(1+\gamma).
		$
		By the choice of $\eps_0$, we deduce that
		\begin{equation}\label{eq:case4-specbd}
			\frac{1+\kappa}{2}<\frac{5+2\gamma}{6}+\eps.
		\end{equation}
		Using \eqref{equ-J-smallGam} and \eqref{eq:eta-bound}, we obtain
		\begin{equation*}
			\mathcal J_\eta(\lambda)^{2\gamma}+\ln(1+\lambda)
			\le C(1+\lambda^\kappa).
		\end{equation*}
		This implies that \eqref{eq4} holds with $\lambda^{1/2}$ replaced by $\lambda^\kappa$.
		Consequently, the analogues of the doubling estimates \eqref{equ-doub-Phi-HL-ball} and \eqref{equ-doub-Phi-HH-ball} hold with $\lambda^{1/2}$ replaced by $\lambda^\kappa$.
		On the other hand, the definition of $\hat R_{\lambda,\eta}$ and \eqref{eq:eta-bound} give
		\begin{align*}
			&\|V_\eta\|_{C^{0,1}(B(0,8\hat C\hat R_{\lambda,\eta})\setminus B_{30\sqrt n})}
			\le C(1+\lambda^{1+\kappa}).
		\end{align*}
		Thus, \eqref{equ-three-sphere} yields
		\begin{equation*}
			\|\hat\Phi\|_{H^1(Q_j\times I_\rho)} \le C\exp\big({C(1+\lambda^{(1+\kappa)/2})}\big) \|\phi\|_{L^2(\Omega_j)}^\alpha\|\hat\Phi\|_{H^1(4Q_j\times I_{4\rho})}^{1-\alpha}.
		\end{equation*}
		Then, proceeding as in \emph{Case 3}, it follows from \eqref{eq:case4-specbd} that
		\begin{equation*}
			\|\phi\|_{L^2(\R^n)}\le C\exp\big({C(1+\lambda^{\frac{5+2\gamma}{6}+\eps})}\big) \|\phi\|_{L^2(\Omega)}.
		\end{equation*}
		
		This completes the proof of Theorem \ref{thm-spec-H-eta}.
	\end{proof}

	\section{Proof of Theorem \ref{thm:weak}}\label{sec3}
	
	The proof has two parts. We first establish the equivalence of the three stabilizability conditions. For $0<\gamma<1$, we then prove the additional equivalence with null controllability.
	
	We begin by introducing the spectral projections and establishing the associated dissipation estimate.
	For $\lambda\ge0$ and $g\in L^2(\mathbb R^n\times\mathbb R^m)$, define
	\begin{equation}\label{equ-def-pi}
		\Pi_\lambda g=\mathcal F_\eta^{-1}\Big(\sum_{\lambda_{\gamma,j}^\eta\le\lambda}\big\langle (\mathcal F_y g)(\cdot,\eta),\phi_{\gamma,j}^\eta\big\rangle_{L^2(\mathbb R^n)}\phi_{\gamma,j}^\eta\Big),
	\end{equation}
	where $(\lambda_{\gamma,j}^\eta,\phi_{\gamma,j}^\eta)$ are the eigenpairs in \eqref{equ-eigen-connect}. The sum is understood for almost every $\eta\ne0$ and is zero when $\lambda<\lambda_{\gamma,1}^\eta$. Orthonormality and Plancherel's theorem show that $\Pi_\lambda$ is an orthogonal projection.
	
	For every $\lambda\ge0$, $t>0$, and $g\in L^2(\R^n\times\R^m)$, we have the following dissipation estimate:
	\begin{equation}\label{equ-dissi-pi}
		\|(1-\Pi_\lambda)(e^{-t(-\Delta_\gamma)}g)\|_{L^2(\R^n\times\R^m)}\le e^{-t\lambda}\|g\|_{L^2(\R^n\times\R^m)}.
	\end{equation}
	Indeed, after taking the partial Fourier transform in $y$, $(1-\Pi_\lambda)$ retains only the eigenmodes of $H_{\gamma,\eta}$ with eigenvalues greater than $\lambda$. This implies $e^{-tH_{\gamma,\eta}}\le e^{-t\lambda}$ on this subspace, and the estimate then follows from orthonormality and Plancherel's theorem.

	\subsection{Stabilizability}
	We prove the implications (i) $\implies$ (ii), (ii) $\implies$ (iii), and (iii) $\implies$ (i), which establish the equivalence in Theorem~\ref{thm:weak}.
	
	\noindent\emph{Proof of (i) $\implies$ (ii).}
	For $0<\gamma<1$, null controllability, proved in Subsection~\ref{subsec:null-control}, implies complete stabilizability.
	For $\gamma\ge1$, we use the following result:
	\begin{theorem}[{\cite[Theorem 4.1]{LWXY22}}]\label{thm-rapid-stab}
		Let $X$ and $U$ be Hilbert spaces, and let $A:D(A)\subset X\to X$ be such that $-A$ generates a $C_0$-semigroup $S(t)=e^{-tA}$ on $X$. Let $B\in\mathcal L(U;X)$, and let $(P_k)_{k\in\N^*}$ be a family of orthogonal projections on $X$. Suppose that for every $k\in\N^*$ there are constants $M_k>0$ and $\alpha_k>0$, with $\alpha_k\to+\infty$ as $k\to\infty$, such that
		\begin{equation}\label{equ-dissip-stab}
			\|(I-P_k)S(t)^*\varphi\|_X\le M_ke^{-\alpha_kt}\|\varphi\|_X,
			\qquad t\ge0,\quad \varphi\in X.
		\end{equation}
		Assume in addition that for every $k\in\N^*$ there is a constant $C_k>0$ such that
		\begin{equation}\label{equ-spectral-stab}
			\|P_k\varphi\|_X\le C_k\|B^*P_k\varphi\|_U,
			\qquad \varphi\in X.
		\end{equation}
		Then the system
		\begin{equation}\label{equ-abstract-stab}
			\partial_tu(t)+Au(t)=Bh(t)
		\end{equation}
		is completely stabilizable.
	\end{theorem}
	For $\gamma\ge1$, apply Theorem~\ref{thm-rapid-stab} with $X=U=L^2(\R^{n+m})$, $A=-\Delta_\gamma$, $B g=\mathds{1}_E g$, and $P_k=\Pi_k$. Since $A$ is self-adjoint, $S(t)^*=e^{-tA}$, so \eqref{equ-dissi-pi} gives \eqref{equ-dissip-stab}. The spectral inequality \eqref{equ-spec-pi}, established in Subsection~\ref{subsec:null-control} under the thickness assumption, gives condition \eqref{equ-spectral-stab}, because $\|B^*\Pi_k g\|_{L^2(\R^{n+m})}=\|\Pi_k g\|_{L^2(E)}$. Hence system \eqref{Equ-Gru-Cont} is completely stabilizable from $E$, proving (i) $\implies$ (ii).

	\noindent\emph{Proof of (ii) $\implies$ (iii).} Complete stabilizability is stronger than exponential stabilizability, so this implication follows directly from the definitions.
	
	\noindent\emph{Proof of (iii) $\implies$ (i).}
	Assume that the system \eqref{Equ-Gru-Cont} is exponentially stabilizable from $E$. We show that both $E_1$ and $E_2$ are thick.

	We will use the characterizations of exponential stabilizability in terms of observability estimates, taken from the work \cite{TWX20}. According to Theorem 1 in \cite{TWX20}, the system \eqref{Equ-Gru-Cont} is exponentially stabilizable from $E$ if and only if for any $\alpha\in(0,1)$, there exist $T>0$ and $C>0$ such that for all $g\in L^2(\R^n\times\R^m)$,
	\begin{equation}\label{equ-weak-obs}
		\|e^{-T(-\Delta_\gamma)}g\|_{L^2(\R^n\times\R^m)}^2\le C\int_0^T\|e^{-t(-\Delta_\gamma)}g\|_{L^2(E)}^2dt+\alpha\|g\|_{L^2(\R^n\times\R^m)}^2.
	\end{equation}
	The proof proceeds in two steps.
	
	\noindent\emph{Step 1: Thickness of $E_1$.}

	Choose $\psi\in\mathcal{S}(\R^m)$ with $\|\psi\|_{2}^2=1$, and $\text{supp }\hat{\psi}\subset B(0,1)$. Let $R\ge1$ and define 
	$$
	\psi_R(y):=R^{-m/2}\psi(y/R), 
	$$
	then
	$\|\psi_R\|_{2}^2=1, \text{supp }\widehat{\psi_R}\subset B(0,1/R)$.
	For every $f(x)\in L^2(\R^n)$, let
	\begin{equation*}
		G_{f,R}(x,y)=f(x)\psi_R(y).
	\end{equation*}
	Then, 
	\begin{equation}\label{eq:L2-gR}
		\|G_{f,R}\|_2^2=\|f\|_{L^2(\R^n)}^2.
	\end{equation}
	
	First, we estimate $\|e^{T\Delta_\gamma}G_{f,R}\|_2^2$. By Plancherel's theorem,
	\begin{align}\label{equ-LHS-A}
		\|e^{T\Delta_\gamma}G_{f,R}\|_2^2&=\int_{\R^n}\int_{\R^m}|e^{-TH_{\gamma,\eta}}f(x)R^{m/2}\hat\psi(R\eta)|^2dxd\eta\notag\\
		&=\int_{\R^m}\int_{\R^n}|e^{-TH_{\gamma,{\eta}/{R}}}f(x)|^2|\hat{\psi}(\eta)|^2dxd\eta,
	\end{align}
	where $H_{\gamma,{\eta}/{R}}=-\Delta_x+R^{-2}|\eta|^2|x|^{2\gamma}$. Notice that $C_0^\infty(\mathbb R^n)$ is a common operator core for $H_{\gamma,\eta/R}$ and $-\Delta_x$. For every $\varphi\in C_0^\infty(\R^n)$, we have 
	$$
	\|H_{\gamma,{\eta}/{R}}\,\varphi-(-\Delta_x)\varphi\|_{L^2(\R^n)}\to 0,\qquad\mbox{as }R\to\infty.
	$$
	Then, $H_{\gamma,{\eta}/{R}}\to -\Delta_x$ in the strong resolvent sense (see \cite[Theorem VIII.25]{RS72}). Since $s\mapsto e^{-ts}$ is bounded and continuous on $[0,\infty)$, functional calculus (\cite[Theorem VIII.20]{RS72}) yields that for every $f\in L^2(\R^n)$
	\begin{equation}\label{equ-appro-A}
		\|e^{-tH_{\gamma,{\eta}/{R}}}\, f-e^{t\Delta_x}f\|_{L^2(\R^n)}\to 0,\quad\mbox{as }R\to\infty.
	\end{equation}
	Therefore, it follows from \eqref{equ-LHS-A} and \eqref{equ-appro-A} that
	\begin{align}\label{equ-lowbd-A}
		\lim_{R\to\infty}\|e^{T\Delta_\gamma}G_{f,R}\|_2^2=\int_{\R^m}\int_{\R^n}|e^{T\Delta_x}f(x)|^2|\hat{\psi}(\eta)|^2dxd\eta
		=\|e^{T\Delta_x}f\|_{L^2(\R^n)}^2.
	\end{align}
	In \eqref{equ-lowbd-A}, the first equality follows from the dominated convergence theorem.
	
	Second, we estimate $\int_0^T\|e^{t\Delta_\gamma}G_{f,R}\|_{L^2(E_1\times\R^m)}^2dt$.
	Let $v(t,x):=e^{t\Delta_x}f(x)$, and $\tilde{u}(t,x,y):=v(t,x)\psi_R(y)$.
	Then,
	\begin{align*}
		\int_0^T\|e^{-tH_{\gamma,\eta}}\widehat{G_{f,R}}(\eta)-\hat{\tilde{u}}(\eta)\|_2^2dt
		=\int_0^T\int_{\R^m}\|e^{-tH_{\gamma,\eta/R}}f-e^{t\Delta_x}f\|_{L^2(\R^n_x)}^2\ |\hat{\psi}(\eta)|^2d\eta dt.
	\end{align*}
	By \eqref{equ-appro-A}, the integrand converges pointwise to zero, while the contraction properties of $e^{-tH_{\gamma,\eta/R}}$ and $e^{t\Delta_x}$ provide a uniform integrable bound. Hence, by the dominated convergence theorem,
	\begin{equation*}
		\int_0^T\|e^{-tH_{\gamma,\eta}}\widehat{G_{f,R}}(\eta)-\hat{\tilde{u}}(\eta)\|_2^2dt
		\to0.
	\end{equation*}
	Therefore,
	\begin{align*}
		\int_0^T\|e^{t\Delta_\gamma}G_{f,R}\|_{L^2(E_1\times\R^m)}^2dt
		&\le 2\int_0^T\|\tilde{u}\|_{L^2(E_1\times\R^m)}^2dt+2\int_0^T\|e^{t\Delta_\gamma}G_{f,R}-\tilde{u}\|_{L^2(E_1\times\R^m)}^2dt\\
		&\le 2\int_0^T\|e^{t\Delta_x}f\|_{L^2(E_1)}^2dt+{o}_R(1),
	\end{align*}
	which implies
	\begin{equation}\label{equ-uppbd-A}
		\lim_{R\to\infty}\int_0^T\|e^{t\Delta_\gamma}G_{f,R}\|_{L^2(E_1\times\R^m)}^2dt\le 2\int_0^T\|e^{t\Delta_x}f\|_{L^2(E_1)}^2dt.
	\end{equation}
	
	Since the system \eqref{Equ-Gru-Cont} is exponentially stabilizable from $E\subset E_1\times\R^m$, we have
	\begin{equation}\label{equ-weak-obs1}
		\|e^{-T(-\Delta_\gamma)}G_{f,R}\|_2^2\le C\int_0^T\|e^{-t(-\Delta_\gamma)}G_{f,R}\|_{L^2(E_1\times\R^m)}^2dt+\alpha\|G_{f,R}\|_2^2.
	\end{equation}
	Letting $R\to\infty$, it follows from \eqref{equ-lowbd-A}, \eqref{equ-uppbd-A} and \eqref{eq:L2-gR} that
	\begin{equation}\label{equ-obs-Lapla}
		\|e^{T\Delta_x}f\|_{L^2(\R^n)}^2\le 2C \int_0^T\|e^{t\Delta_x}f\|_{L^2(E_1)}^2dt+\alpha\|f\|_{L^2(\R^n)}^2,
	\end{equation}
	for every $f\in L^2(\R^n)$.
	That means the system
	\begin{equation}\label{Equ-Laplac}
		\begin{cases}
			\partial_{t}f(t,x)-\Delta_x f(t,x)=\1_{E_1} h(t,x),\;\;(t,x)\in(0,T)\times\R^n,\\
			f(0,x)=f_0(x)\in L^{2}(\R^n),
		\end{cases}
	\end{equation}
	is exponentially stabilizable from $E_1\subset\R^n$. Thus, it follows from \cite[Theorem 1.1]{HWW21} or \cite[Theorem 2.1]{AM22} that $E_1\subset\R^n$ is a thick set.

	\noindent\emph{Step 2: Thickness of $E_2$.}
	
	Now, we fix $\alpha\in(0,1/2)$ and $T>0$ in \eqref{equ-weak-obs}, and define
	\begin{equation}\label{eq:defsigma}
		\sigma:=\frac 1 2\Big(\frac{\ln(\frac{1}{2\alpha})}{2T\lambda_{\gamma,1}}\Big)^{\frac{1+\gamma}{2}}>0,
	\end{equation}
	where $\lambda_{\gamma,1}$ is the smallest eigenvalue of $H_{\gamma}$.
	Fix $ y_0\in\R^m$ and consider the following function:
	\begin{equation}\label{equ-examp-contr}
		g(x,y)=\mathcal{F}_\eta^{-1}\big(\phi_{\gamma,1}^\eta(x) (\mathcal{F}_yh)(\eta)\chi_{(\sigma,2\sigma)}(|\eta|)\big)(x,y),\quad\forall\,(x,y)\in\R^n\times\R^m,
	\end{equation}
	where $h(y)=\exp(-\frac{|y-y_0|^2}{2})$ and $\phi_{\gamma,1}^\eta$ is the eigenfunction of $H_{\gamma,\eta}$ associated with the eigenvalue $\lambda_{\gamma,1}^\eta$.
	By calculation,
	\begin{equation}\label{equ-h-Four}
		(\mathcal{F}_yh)(\eta)=\exp(-iy_0\cdot \eta-\frac{|\eta|^2}{2}),
	\end{equation}
	and 
	\begin{equation}\label{equ-thm1-8}
		\|g(x,y)\|_{L^2(\R^n\times\R^m)}^2=\int_{\sigma<|\eta|<2\sigma}e^{-|\eta|^2}d\eta:=c_\sigma.
	\end{equation}

	First, we establish a lower bound for $\|e^{-T(-\Delta_\gamma)}g\|_{L^2(\R^n\times\R^m)}^2$. In fact, by Plancherel's theorem, together with \eqref{equ-examp-contr} and \eqref{equ-h-Four}, we have
	\begin{align*}\label{equ-thm1-1}
		\|e^{-T(-\Delta_\gamma)}g\|_{L^2(\R^n\times\R^m)}^2&=\|e^{-TH_{\gamma,\eta}}(\mathcal{F}_yg)(x,\eta)\|_{L^2(\R^n\times\R^m_\eta)}^2\notag\\
		&=\int_{\sigma<|\eta|<2\sigma}e^{-|\eta|^2}d\eta\int_{\R^n}\big|e^{-T\lambda_{\gamma,1}^\eta}\phi_{\gamma,1}^\eta(x)\big|^2dx\notag\\
		&\ge \exp\big({-2T\lambda_{\gamma,1}(2\sigma)^{\frac{2}{1+\gamma}}}\big)\cdot\int_{\sigma<|\eta|<2\sigma}e^{-|\eta|^2}d\eta \notag\\
		&= c_\sigma\cdot \exp\big({-2T\lambda_{\gamma,1}(2\sigma)^{\frac{2}{1+\gamma}}}\big).
	\end{align*}
	This, combining with \eqref{equ-weak-obs}, \eqref{equ-thm1-8} and \eqref{eq:defsigma}, implies  
	\begin{equation}\label{equ-thm1-10}
		\frac{c_\sigma}{2}\cdot \exp\big({-2T\lambda_{\gamma,1}(2\sigma)^{\frac{2}{1+\gamma}}}\big)\le C\int_0^{T}\|e^{-t(-\Delta_\gamma)}g\|_{L^2(E)}^2dt.
	\end{equation}
	
	Next, we estimate an upper bound for $\int_{0}^{T}\|e^{-t(-\Delta_\gamma)}g\|_{L^2(E)}^2dt$. By \eqref{equ-examp-contr}, \eqref{equ-h-Four},
	\begin{align}\label{equ-thm1-5}
		\int_{0}^{T}\|e^{-t(-\Delta_\gamma)}g\|_{L^2(E)}^2dt&=\int_{0}^{T}\|\mathcal{F}_\eta^{-1}e^{-tH_{\gamma,\eta}}(\mathcal{F}_yg)(x,\eta)\|_{L^2(E)}^2dt\notag\\
		&\le\int_{0}^{T}\int_{\R^n}\int_{E_2}\Big|\mathcal{F}_\eta^{-1}\big(e^{-tH_{\gamma,\eta}}\phi_{\gamma,1}^\eta(x)e^{-\frac{|\eta|^2}{2}}\chi_{(\sigma,2\sigma)}(|\eta|)\big)(x,y-y_0) \Big|^2dxdydt\notag\\
		&=\int_{0}^{T}\int_{\R^n}\int_{E_2-y_0}\Big|\mathcal{F}_\eta^{-1}\big(e^{-t\lambda_{\gamma,1}^\eta}\phi_{\gamma,1}^\eta(x)e^{-\frac{|\eta|^2}{2}}\chi_{(\sigma,2\sigma)}(|\eta|)\big)(x,y) \Big|^2dxdydt\notag\\
		&\le \int_{0}^{T}\int_{(E_2-y_0)\cap[-L_y,L_y]^m}\|{\psi}(x,y,t)\|_{L^2(\R^n_x)}^2dydt\notag\\
		&\quad+\int_{0}^{T}\int_{\R^n}\int_{\R^m\setminus B(0,L_y)}| {\psi}(x,y,t)|^2dxdydt\notag\\
		&:=I+II,
	\end{align}
	where ${\psi}(x,y,t)$ is given by 
	\begin{equation*}
		{\psi}(x,y,t)=\mathcal{F}_\eta^{-1}\big(e^{-t\lambda_{\gamma,1}^\eta}\phi_{\gamma,1}^\eta(x)e^{-\frac{|\eta|^2}{2}}\chi_{(\sigma,2\sigma)}(|\eta|)\big)(x,y).
	\end{equation*}
	We claim there exists a positive constant $C_1<\infty$, depending on $m,\gamma,\alpha,T$, such that
	\begin{equation}\label{equ-psi-sup}
		\|{\psi}(x,y,t)\|_{L^2(\R^n_x)}^2\le {C_1}<\infty, \quad\forall y\in\R^m,\,t>0.
	\end{equation}
	In fact, using the $L^1-L^\infty$ boundedness of the Fourier transform and H\"{o}lder inequality, we have
	\begin{align*}
		\quad\|{\psi}(x,y,t)\|_{L^2(\R^n_x)}^2
		&\le (2\pi)^{-m} \int_{\R_x^n}\Big|\int_{\R^m_\eta}|e^{-t\lambda_{\gamma,1}^\eta}\phi_{\gamma,1}^\eta(x)e^{-\frac{|\eta|^2}{2}}\chi_{(\sigma,2\sigma)}(|\eta|) |d\eta\Big|^2dx\\
		&\le C_{m,\sigma}\int_{\R_x^n} \int_{\sigma<|\eta|<2\sigma}|e^{-t\lambda_{\gamma,1}^\eta}\phi_{\gamma,1}^\eta(x)e^{-\frac{|\eta|^2}{2}}|^2dxd\eta\\
		&\le C_{m,\sigma}\int_{\sigma<|\eta|<2\sigma}|e^{-\frac{|\eta|^2}{2}}|^2d\eta\\
		&\le C_1.
	\end{align*}
	It follows from \eqref{equ-psi-sup} that
	\begin{align}\label{equ-bound-I}
		I\le C_1T |(E_2-y_0)\cap[-L_y,L_y]^m|.
	\end{align}
	On the other hand, by Plancherel's theorem,
	\begin{align*}
		\|\psi(x,y,t)\|_{L^2(\R^n\times\R^m\times(0,T))}^2&=\int_{0}^T\int_{\R^n}\int_{\sigma<|\eta|<2\sigma} \big|e^{-t\lambda_{\gamma,1}^\eta}\phi_{\gamma,1}^\eta(x)e^{-\frac{|\eta|^2}{2}}\big|^2dxd\eta dt\\
		&\le T\int_{\R^m_\eta} e^{-|\eta|^2}d\eta=\pi^{m/2}T<\infty,
	\end{align*}
	which means $\psi\in L^2(\R^n\times\R^m\times(0,T))$.
	As a consequence, we can take $L_y>1$ (independent of $y_0$) sufficiently large such that
	\begin{equation}\label{equ-bound-II-stab}
		II\le \frac{c_\sigma}{4C}\cdot \exp\big({-2T\lambda_{\gamma,1}(2\sigma)^{\frac{2}{1+\gamma}}}\big).
	\end{equation}
	
	By \eqref{equ-thm1-10}, \eqref{equ-thm1-5}, \eqref{equ-bound-I}, and \eqref{equ-bound-II-stab}, we obtain that 
	\begin{equation*}
		\frac{c_\sigma}{4C_1CT}\cdot \exp\big({-2T\lambda_{\gamma,1}(2\sigma)^{\frac{2}{1+\gamma}}}\big)\le |(E_2-y_0)\cap[-L_y,L_y]^m|,\quad\forall y_0\in\R^m.
	\end{equation*}
	This shows that $E_2$ is thick.

\subsection{Null controllability for $0<\gamma<1$}\label{subsec:null-control}
	
	We first prove that thickness implies null controllability, using the classical Lebeau--Robbiano strategy.
	\begin{theorem}(\cite{BPS18})\label{thm-L-R-method}
		Let $\mathcal{O}$ be an open subset of $\R^n$, $\Omega$ be a measurable subset of $\mathcal{O}$, $(\pi_k)_{k\in\N^*}$ be a family of orthogonal projections defined on $L^2(\mathcal{O})$, $(e^{-tA})_{t\ge 0}$ be a strongly continuous contraction semigroup on $L^2(\mathcal{O})$; $c_1,c_2,a,b,t_0,l>0$ be positive constants with $a<b$. If the following spectral inequality
		\begin{equation}\label{equ-spectral}
			\forall g\in L^2(\mathcal{O}),\,\forall k\ge 1,\quad\quad\|\pi_kg\|_{L^2(\mathcal{O})}\le e^{c_1k^a}\|\pi_kg\|_{L^2(\Omega)},
		\end{equation}
		and the following dissipation estimate
		\begin{equation}\label{equ-dissipation}
			\forall g\in L^2(\mathcal{O}),\,\forall k\ge 1,\;\forall 0<t<t_0,\;\; \|(1-\pi_k)(e^{-tA}g)\|_{L^2(\mathcal{O})}\le \frac{1}{c_2}e^{-c_2t^l k^b}\|g\|_{L^2(\mathcal{O})},
		\end{equation}
		hold, then there exists a positive constant $C>1$ such that the following observability estimate holds
		\begin{equation}\label{equ-observability}
			\forall T>0,\,\forall g\in L^2(\mathcal{O}),\quad \|e^{-TA}g\|_{L^2(\mathcal{O})}^2\le C \exp{(\frac{C}{T^{\frac{al}{b-a}}})}\int_{0}^{T}\|e^{-tA}g\|_{L^2(\Omega)}^2dt.
		\end{equation}
	\end{theorem}
	
	Noting that for all $g\in L^2(\R^n\times\R^m)$,
	\begin{equation}\label{eq:Ran-H-gam-eta}
		\mathcal{F}_y(\Pi_\lambda g)=\sum_{\lambda_{\gamma,j}^\eta\le\lambda}\langle(\mathcal{F}_yg)(x,\eta),\phi_{\gamma,j}^\eta\rangle\phi_{\gamma,j}^\eta\in \text{Ran}(P_\lambda(H_{\gamma,\eta})).
	\end{equation}
	For fixed $x\in\R^n$, it follows from \eqref{equ-eigen-connect} that
	\begin{equation}\label{equ-eta-range}
		\mbox{supp}\;\mathcal{F}_y(\Pi_\lambda g)(\eta)\subset \Big\{\eta\in\R^m:\;|\eta|\le (\frac{\lambda}{\lambda_{\gamma,1}})^{\frac{1+\gamma}{2}}\Big\}.
	\end{equation}
	We proceed in two cases:
	
	\noindent\emph{Case 1. $\gamma\ge1/2$.}
	By \eqref{eq:Ran-H-gam-eta}, Plancherel's theorem and Theorem \ref{thm-spec-H-eta}, for all $\lambda\ge 0$, $\eps>0$, and $\gamma\ge 1/2$, we have
	\begin{align}\label{equ-spec-pi-part}
		\|\Pi_\lambda g\|_{L^2(\R^n\times\R^m)}^2&=\int_{\R^m}\big\|\mathcal{F}_y(\Pi_\lambda g) \big\|^2_{L^2(\R^n)}d\eta\notag\\
		&\le Ce^{C(1+
		\lambda^{\frac{1+\gamma}{2}+\eps})}
		\int_{\R^m}\|\mathcal{F}_y(\Pi_\lambda g)\|_{L^2(E_1)}^2d\eta\notag\\
		&=Ce^{C(1+
			\lambda^{\frac{1+\gamma}{2}+\eps})}
		\int_{\R^m}\|\Pi_\lambda g\|_{L^2(E_1)}^2dy.
	\end{align}
	On the other hand, by \eqref{equ-eta-range} and quantitative Logvinenko–Sereda Theorem (see \cite{Kov01}), we derive that 
	\begin{align}\label{eq-LS}
		\|\Pi_\lambda g(x,\cdot)\|_{L^2(\R^m)}^2\le e^{C(1+\lambda^{\frac{1+\gamma}{2}})}
		\|\Pi_\lambda g(x,\cdot)\|_{L^2(E_2)}^2.
	\end{align}
	It follows from \eqref{equ-spec-pi-part} and \eqref{eq-LS} that for any $\gamma\ge 1/2$,
	\begin{equation}\label{equ-spec-pi}
		\|\Pi_\lambda g\|_{L^2(\R^n\times\R^m)}^2\le Ce^{C(1+\lambda^{\frac{1+\gamma}{2}+\eps})} \|\Pi_\lambda g\|_{L^2(E_1\times E_2)}^2.
	\end{equation}
	
	Due to $\gamma<1$, we can take $\eps=\frac{1-\gamma}{4}>0$. Then
	\begin{equation}\label{eq5}
		\frac{1+\gamma}{2}+\eps=\frac{3+\gamma}{4}<1.
	\end{equation}
	Applying Theorem \ref{thm-L-R-method} with $d=n+m,\mathcal{O}=\R^n\times\R^m,\Omega=E_1\times E_2$, and $A=-\Delta_\gamma$. By \eqref{equ-spec-pi}, \eqref{eq5}, and \eqref{equ-dissi-pi},  we conclude that 
	the following observability estimate holds: for any $T>0$, every $g\in L^2(\R^n\times\R^m)$,
	\begin{equation}\label{equ-obs-Gru-part}
		\|e^{-T(-\Delta_\gamma)}g\|_{L^2(\R^n\times\R^m)}^2\le C\exp(CT^{-\frac{3+\gamma}{1-\gamma}})\int_{0}^{T}\|e^{-t(-\Delta_\gamma)}g\|_{L^2(E_1\times E_2)}^2dt.
	\end{equation}
	
	\noindent\emph{Case 2. $\gamma \in (0,1/2)$.}
	By \eqref{eq:Ran-H-gam-eta}, Plancherel's theorem and Theorem \ref{thm-spec-H-eta}, for all $\lambda\ge 0$, and $\eps>0$, we have
	\begin{align}\label{equ-spec-pi-part1}
		\|\Pi_\lambda g\|_{L^2(\R^n\times\R^m)}^2\le Ce^{C(1+\lambda^{\frac{5+2\gamma}{6}+\eps})}\int_{\R^m}\|\Pi_\lambda g\|_{L^2(E_1)}^2dy
	\end{align}
	Since $\gamma<1/2$, we can choose $\eps=\frac{1-2\gamma}{12}>0$. Then, 
	\begin{equation}\label{eq1}
		\frac{1+\gamma}{2}<\frac{3}{4}<\frac{5+2\gamma}{6}+\eps=\frac{11+2\gamma}{12}<1.
	\end{equation}
	It follows from \eqref{equ-spec-pi-part1}, \eqref{eq1} and \eqref{eq-LS} that
	\begin{equation*}
		\|\Pi_\lambda g\|_{L^2(\R^n\times\R^m)}^2\le Ce^{C(1+\lambda^{\frac{11+2\gamma}{12}})} \|\Pi_\lambda g\|_{L^2(E_1\times E_2)}^2.
	\end{equation*}
	Using Theorem \ref{thm-L-R-method} again, we obtain: for any $T>0$, every $g\in L^2(\R^n\times\R^m)$,
	\begin{equation}\label{equ-obs-Gru-part1}
		\|e^{-T(-\Delta_\gamma)}g\|_{L^2(\R^n\times\R^m)}^2\le C\exp(CT^{-\frac{11+2\gamma}{1-2\gamma}})\int_{0}^{T}\|e^{-t(-\Delta_\gamma)}g\|_{L^2(E_1\times E_2)}^2dt.
	\end{equation}
	
	By the Hilbert Uniqueness Method, the observability estimates \eqref{equ-obs-Gru-part} and \eqref{equ-obs-Gru-part1} imply null controllability from $E=E_1\times E_2$ for every $T>0$. Thus, (i) implies null controllability.

Conversely, null controllability in a fixed positive time yields complete stabilizability by \cite[Proposition~21]{TWX20}, so it implies (ii). The stabilizability subsection gives (ii) $\implies$ (iii) $\implies$ (i). Therefore, for $0<\gamma<1$, null controllability is equivalent to each of (i)--(iii) in Theorem~\ref{thm:weak}.
	
	\section{Proof of Theorem \ref{thm:mini-timel}}\label{sec4}
	
	Let $E=E_1\times\mathbb R^m$, where $E_1\subset\mathbb R^n$ is thick, and set $a=a_{\mathrm{ess}}(E_1)$. We prove the time characterization in Theorem~\ref{thm:mini-timel}: if $a>0$, null controllability fails for $0<T\le \frac{a^2}{2n}$ and holds for $T>\frac{a^2}{2n}$; if $a=0$, it holds for every $T>0$.
	
	\subsection{Hermite preliminaries}
	We collect the Hermite basis and a finite-dimensional estimate used in the uniform observability argument below.
	For $k\in\N$, the one-dimensional Hermite function is given by
	$$
	h_k(t)=\frac{1}{(2^k k!\sqrt{\pi})^{1/2}}
	H_k(t)e^{-t^2/2},\qquad t\in\R,
	$$
	where $H_k$ denotes the $k$-th Hermite polynomial. For a multi-index
	$\alpha=(\alpha_1,\ldots,\alpha_n)\in\N^n$, set
	$$
	\Phi_\alpha(x)=\prod_{j=1}^n h_{\alpha_j}(x_j),
	\qquad x=(x_1,\ldots,x_n)\in\R^n.
	$$
	Then $\{\Phi_\alpha\}_{\alpha\in\N^n}$ forms an orthonormal basis of $L^2(\R^n)$ and
	$$
	(-\Delta+|x|^2)\Phi_\alpha=(2|\alpha|+n)\Phi_\alpha,\qquad|\alpha|=\alpha_1+\cdots+\alpha_n.
	$$
	
	Now, we consider $H_s:=-\Delta+s^2|x|^2$, where $s>0$.
	Defining
	$$
	\Phi_{\alpha,s}(x)=s^{n/4}\Phi_\alpha(\sqrt{s}\,x),
	$$
	then, by scaling, we obtain
	\begin{equation}\label{eq:eigen-Hs}
		(-\Delta+s^2|x|^2)\Phi_{\alpha,s}=s(2|\alpha|+n)\Phi_{\alpha,s},
	\end{equation}
	and $\{\Phi_{\alpha,s}\}_{\alpha\in\N^n}$ is again an orthonormal
	basis of $L^2(\R^n)$. 
	We have
	$$
	\Phi_{\alpha,s}(x)=s^{n/4}\kappa_\alpha H_\alpha(\sqrt{s}\,x)e^{-s|x|^2/2},
	$$
	where
	$$
	H_\alpha(x)=\prod_{j=1}^n H_{\alpha_j}(x_j)
	$$
	is a multivariate Hermite polynomial, and
	$$
	\kappa_\alpha=\frac{1}{\pi^{n/4}2^{|\alpha|/2}(\alpha!)^{1/2}},\qquad\alpha!=\alpha_1!\cdots\alpha_n!.
	$$
	Moreover, the operator $-H_s$ generates the Hermite semigroup $e^{-tH_s}$ whose kernel is explicitly given by
	\begin{equation*}
		K_{t,s}(x,y)=\sum_{\alpha\in\mathbb{N}^n}e^{-(2|\alpha|+n)st}\Phi_{\alpha,s}(x)\Phi_{\alpha,s}(y).
	\end{equation*}
	By Mehler's formula, for $t>0$,
	\begin{equation}\label{eq:ker-Hs}
		\begin{aligned}
			K_{t,s}(x,y)=\Big(\frac{s}{2\pi\sinh(2s t)}\Big)^{n/2}\cdot\exp\Big(-\frac{s\tanh(s t)}{4}|x+y|^2-\frac{s}{4\tanh(s t)}|x-y|^2\Big).
		\end{aligned}
	\end{equation}

	For $N\in\N$ and $s>0$, let
	\begin{equation}\label{eq:XNs}
		X_{N,s}:=\operatorname{span}\{\Phi_{\alpha,s}:|\alpha|\le N\},
	\end{equation}
	and let $P_{N,s}$ be the orthogonal projection onto $X_{N,s}$. The next lemma gives the Hermite spectral inequality needed below. Unlike Theorem~\ref{thm-spec-H-eta}, which gives a spectral inequality for $H_{\gamma,\eta}$ uniform in $\eta$ under a thickness assumption, the lemma applies to $X_{N,s}$ with $N$ fixed and only assumes that $E_1$ is measurable. It also differs from general Hermite spectral estimates for thick sets such as \cite[Theorem 2.1]{BJP21}: its exponential rate explicitly records the essential distance $a_{\rm ess}(E_1)$, through the factor $e^{b^2s}$ for every $b>a_{\rm ess}(E_1)$.
	\begin{lemma}\label{lem:sharp-low-spec}
		Assume $E_1\subset\R^n$ is measurable and let $a_{\rm ess}=a_{\rm ess}(E_1)$. Fix $N\in\N$ and $b>a_{\rm ess}$. There is $C=C(N,b,E_1)>0$ such that
		\begin{equation}\label{eq:sharp-low-mode}
			\|f\|_{L^2(\R^n)}^2\le C s^{-n/2}e^{b^2s}\|f\|_{L^2(E_1)}^2
		\end{equation}
		for every $s\ge1$ and every $f\in X_{N,s}$.
	\end{lemma}
	\begin{proof}
		Note that $\{H_\alpha:|\alpha|\le N\}$ is a basis of the space $\mathcal P_N$ of polynomials on $\R^n$ of degree at most $N$. Then, every $f\in X_{N,s}$ can be written as
		\begin{equation}\label{eq:hermi-expan}
			f(x)=\sum_{|\alpha|\le N}c_\alpha\Phi_{\alpha,s}=s^{n/4}e^{-s|x|^2/2}q(\sqrt{s}x),
		\end{equation}
		with
		\begin{equation}\label{eq:poly-q}
			q(x)=\sum_{|\alpha|\le N}c_\alpha\kappa_\alpha
			H_\alpha(x)=\sum_{|\nu|\le N}a_\nu x^\nu.
		\end{equation}
		
		Let $F_b:=E_1\cap B(0,b)$. Since $b>a_{\rm ess}$, we have $|F_b|>0$. As a nonzero polynomial cannot vanish on a set of positive measure, $\|\cdot\|_{L^2(F_b)}$ defines a norm on $\mathcal P_N$.
		Hence, by equivalence of norms, there exists
		$C_{N,F_b}>0$ (independent of $s$), such that
		\begin{equation}\label{eq:poly-norm-equivalence}
			\int_{F_b}|p(y)|^2\,dy\ge  C_{N,F_b} \sum_{|\nu|\le N}|a_\nu|^2,\qquad \forall p(y)=\sum_{|\nu|\le N} a_\nu y^\nu\in\mathcal P_N.
		\end{equation}
		Due to
		$$
		q(\sqrt{s}x)=\sum_{|\nu|\le N}a_\nu s^{|\nu|/2} x^\nu\in \mathcal P_N,
		$$
		it follows from \eqref{eq:poly-norm-equivalence} and the fact $s\ge 1$ that
		\begin{equation}\label{eq:B1}
			\int_{F_b}|q(\sqrt{s}x)|^2\,dx\ge  C_{N,F_b} \sum_{|\nu|\le N}s^{|\nu|}|a_\nu|^2\ge C_{N,F_b} \sum_{|\nu|\le N}|a_\nu|^2.
		\end{equation}
		Notice that both $\{H_\alpha:|\alpha|\le N\}$ and $\{y^\nu:|\nu|\le N\}$ are bases of the space $\mathcal P_N$. Thus, there exists $\tilde{C}_{N}>0$, such that
		\begin{equation}\label{eq:B2}
			\sum_{|\nu|\le N}|a_\nu|^2 \ge\tilde{C}_{N}
			\sum_{|\alpha|\le N}|c_\alpha\kappa_\alpha|^2\ge {C}_{N} \sum_{|\alpha|\le N}|c_\alpha|^2,
		\end{equation}
		where the last inequality follows from the fact that $\kappa_\alpha \ge C_N'$ whenever $|\alpha| \le N$.
		By \eqref{eq:B1}, \eqref{eq:B2} and the orthonormality of $\Phi_{\alpha,s}$, we derive  
		\begin{equation*}
			\int_{F_b}|q(\sqrt{s}x)|^2\,dx\ge C_1 \sum_{|\alpha|\le N}|c_\alpha|^2= C_1 \|f(x)\|_{L^2(\R^n)}^2,
		\end{equation*} 
		where $C_1$ is independent of $s$.
		Thus, for every $f\in X_{N,s}$
		\begin{align*}
			\|f\|_{L^2(\R^n)}^2\le\frac{1}{C_1}\int_{F_b}|q(\sqrt{s}x)|^2\,dx\le \frac{s^{-n/2}}{C_1}e^{b^2s}\|f\|_{L^2(F_b)}^2\le C s^{-n/2}e^{b^2s}\|f\|_{L^2(E_1)}^2.
		\end{align*}
		Therefore, we complete the proof.
	\end{proof}
	
	\subsection{Necessity: failure for $0<T\le a^2/(2n)$}
	Assume $a>0$ and fix $0<T\le a^2/(2n)$. By the definition of essential distance, $E_1\subset B(0,a)^c$ up to a null set. For $s>0$, consider the operator
	$$
	H_s=-\Delta_x+s^2|x|^2.
	$$
	We consider the normalized Gaussian
	\begin{equation}\label{eq:Psi-s}
		\Psi_s(x):=\big(\frac{2s}{\pi}\big)^{n/4}e^{-s|x|^2},\qquad\|\Psi_s\|_2=1.
	\end{equation}
	Recall that the kernel $K_{t,s}$ of $e^{-tH_s}$ is given by \eqref{eq:ker-Hs}. Then we have
	\begin{align*}
		\big(e^{-tH_s}\Psi_s\big)(x)&=\int_{\R^n}\ K_{t,s}(x,y)\Psi_s(y)dy\\
		&=\big(\frac{2s}{\pi}\big)^{n/4}\Big(\frac{s}{2\pi\sinh(2s t)}\Big)^{n/2}e^{-\frac{s}{4} (r+\frac{1}{r})|x|^2}\int_{\R^n}e^{-b|y|^2+h\cdot y}dy\\
		&=\big(\frac{2s}{\pi}\big)^{n/4}\Big(\frac{s}{2\pi\sinh(2s t)}\Big)^{n/2}e^{-\frac{s}{4} (r+\frac{1}{r})|x|^2}(\frac{\pi}{b})^{n/2}e^{\frac{|h|^2}{4b}},
	\end{align*}
	where
	$$
	r:=\tanh(st),\;b:=\frac{s}{4}(r+\frac{1}{r}+4),\;h:=\frac{s}{2}(\frac{1}{r}-r)x.
	$$
	A direct computation gives
	\begin{equation}\label{eq:solu-Psi}
		\big(e^{-tH_s}\Psi_s\big)(x)=\Big(\sqrt{\frac{2}{\pi}}\frac{2}{3-q}\Big)^{n/2}s^{n/4} e^{-nst}e^{-\frac{s}{2} \frac{3+q}{3-q}|x|^2},\quad\mbox{with }q=e^{-4st}.
	\end{equation}
	
	Then, by \eqref{eq:solu-Psi} and the fact $q\in(0,1)$, there is $C_n>0$ depending only on $n$, such that
	\begin{equation}\label{eq:lowbd-Psi}
		\|e^{-TH_s}\Psi_s\|_2^2=\Big(\sqrt{\frac{2}{\pi}}\frac{2}{3-q}\Big)^{n}\big(\frac{\pi(3-q)}{3+q}\big)^{n/2}  e^{-2nsT}\ge C_n e^{-2nsT}.
	\end{equation}
	On the other hand, observe that 
	\begin{equation*}
		\int_{|x|\ge a}e^{-A|x|^2}dx=c_n A^{-n/2}\Gamma(\frac{n}{2},Aa^2),
	\end{equation*}
	where $\Gamma(\alpha,z)=\int_z^\infty t^{\alpha-1}e^{-t}dt$ is the upper incomplete gamma function. By the standard asymptotic expansion of $\Gamma(\alpha,z)$ (see, e.g. \cite[section 8.11]{OLBC10}), we have
	\begin{equation}\label{eq:Gauss-tail}
		\int_{|x|\ge a}e^{-A|x|^2}dx\sim C_{n,a}A^{-1}e^{-Aa^2},\qquad A\to\infty.
	\end{equation}
	Noting that $E_1\subset B(0,a)^c$ up to a null set, by 
	\eqref{eq:solu-Psi}, \eqref{eq:Gauss-tail} and the fact $q\in(0,1)$, there is $s_0>0$ such that for any $s>s_0$ 
	\begin{align}\label{eq:uppbd-Psi}
		&\quad\int_0^T\|e^{-tH_s}\Psi_s\|_{L^2(E_1)}^2dt\le C_{n,a}s^{n/2-1}e^{-{sa^2}} \int_0^T e^{-2nst}e^{- \frac{2sq}{3-q}a^2}dt\notag\\
		&\le C_{n,a}s^{n/2-1} e^{-{sa^2}}\int_0^T e^{-2nst}\exp\big(-\frac{2a^2s}{3}e^{-4st}\big)dt\notag\\
		&\le C_{n,a}s^{n/2-2} e^{-{sa^2}}\int_0^\infty e^{-2n\tau}\exp\big(-\frac{2a^2s}{3}e^{-4\tau}\big)d\tau\qquad(\mbox{setting }\tau=st)\notag\\
		&= C_{n,a}s^{n/2-2} e^{-{sa^2}}\frac{1}{4}\int_0^1 w^{\frac{n}{2}-1}e^{-\frac{2a^2}{3}sw}dw\qquad(\mbox{setting }w=e^{-4\tau})\notag\\
		&\le C_{n,a}s^{-2}e^{-{sa^2}}.
	\end{align}
	
	Suppose, for contradiction, that null controllability holds for this $T$. By the Hilbert Uniqueness Method, there would then exist $C_{E_1,T}>0$ such that
	\begin{equation}\label{eq:obs-critical}
		\|e^{T\Delta_1}f\|_2^2\le C_{E_1,T}\int_0^T \|e^{t\Delta_1}f\|_{L^2(E_1\times\R^m)}^2dt,
	\end{equation}
	for every $f\in L^2(\R^n\times\R^m)$. 
	Let $k>s_0$, and set
	$$
	f_k(x,y)=\mathcal{F}_\eta^{-1}\big(\Psi_{|\eta|}(x)\chi_{(k,k+1)}(|\eta|)\big)(x,y)\in L^2(\R^n\times\R^m).
	$$
	Applying \eqref{eq:obs-critical} with $f=f_k$, by Plancherel's theorem, \eqref{eq:lowbd-Psi} and \eqref{eq:uppbd-Psi}, we obtain
	\begin{equation*}
		C_n\int_{k<|\eta|<k+1}e^{-2nT|\eta|}d\eta\le C_{E_1,T}\ C_{n,a} \int_{k<|\eta|<k+1}|\eta|^{-2}e^{-a^2|\eta|}d\eta.
	\end{equation*}
	Since $0<T\le \frac{a^2}{2n}$, it follows that
	$$
	\frac{C_ne^{-2nT}}{C_{E_1,T}C_{n,a}}\le k^{-2} e^{-(a^2-2nT)k}\to 0, \qquad\mbox{as }k\to\infty.
	$$
	This contradiction shows that system \eqref{Equ-Gru-Cont} is not null controllable from $E_1\times\R^m$ in time $0<T\le\frac{a^2}{2n}$.

	\subsection{Sufficiency: null controllability for $T>a^2/(2n)$}
	Fix $T>a^2/(2n)$. By the Hilbert Uniqueness Method, it suffices to prove the observability estimate
	\begin{equation}\label{equ-obs-Tlarge}
		\|e^{T\Delta_1}f\|_{L^2(\mathbb{R}^{n+m})}^2 \le C_{\mathrm{obs}} \int_0^T \|e^{t\Delta_1}f\|_{L^2(E_1\times\mathbb{R}^m)}^2 \, dt,
	\end{equation}
	where $\Delta_1 = \Delta_x + |x|^2\Delta_y$.
	
	The key estimate is the following observability inequality for $H_s$, uniform in $s>0$:
	\begin{lemma}\label{lem:uniform-obs-Hs}
		Let $s>0$, $H_s:=-\Delta_x+s^2|x|^2$, and let $E_1\subset\R^n$ be a thick set. For any $T>\frac{a^2}{2n}$, there is $C_T>0$, independent of $s$, such that
		\begin{equation}\label{eq:uniform-fiber-critical}
			\|e^{-TH_s}g\|_{L^2(\R^n)}^2\le C_T\int_0^T \|e^{-tH_s}g\|_{L^2(E_1)}^2\,dt
		\end{equation}
		for every $g\in L^2(\R^n)$.
	\end{lemma}
	
	\begin{proof}
		We first consider $s\ge1$. 
		For every $T>0$, it follows from \cite[Proposition 33]{AS24} that there exist $C,\,d_0>0$ depending only on $n,\,E_1$, such that for all $s\ge 1,\,T>0$, and $g\in L^2(\R^n)$,
		\begin{equation}\label{equ-obs-hermite}
			\|e^{-TH_{s}}g\|_{L^2(\R^n)}^2\le \frac{C }{T}e^{d_0(1+s+T^{-1})}\int_0^T\|e^{-tH_{s}}g\|_{L^2(E_1)}^2dt.
		\end{equation}
		Taking $b>0$ such that
		\begin{equation}\label{eq:b-critical-choice}
			a<b<\sqrt{2nT},
		\end{equation}
		and choosing
		\begin{equation}\label{eq:time-critical-choice}
			0<t_0<t_1<t_2<T-\frac{b^2}{2n}.
		\end{equation}
		Set
		\begin{equation}\label{eq:critical-decomposition}
			u_s(t):=e^{-tH_s}g,\qquad v_s(t):=P_{N,s}u_s(t),\qquad
			w_s(t):=(I-P_{N,s})u_s(t),
		\end{equation}
		where $P_{N,s}$ is the orthogonal projection onto the spectral subspace $X_{N,s}$ of $H_s$, defined in \eqref{eq:XNs} above. Here, $N$ is a fixed, sufficiently large integer to be specified later.
		
		By dissipativity and \eqref{equ-obs-hermite}, for $t\ge t_1>t_0$,
		\begin{align*}
			\|w_s(t)\|_2^2&=\|(I-P_{N,s})e^{-(t-t_0)H_s}e^{-t_0H_s}g\|_2^2\\
			&\le  e^{-2(2N+n)s(t-t_0)} \|u_s(t_0)\|_2^2\\
			&\le C_1 e^{-\big(2(2N+n)(t-t_0)-d_0\big)s}\int_0^{t_0}\|u_s(\tau)\|_{L^2(E_1)}^2d\tau.
		\end{align*}
		Fixing $N$ sufficiently large such that 
		\begin{equation*}\label{eq:N-critical-choice}
			2(2N+n)(t_1-t_0)>d_0,
		\end{equation*}
		we obtain
		\begin{equation}\label{eq:high-mode-uniform}
			\|w_s(t)\|_2^2\le C_1\int_0^{t_0} \|u_s(\tau)\|_{L^2(E_1)}^2\,d\tau, \qquad t\ge t_1,\quad s\ge1,
		\end{equation}
		where $C_1>0$ is independent of $s$ and $t$.
		
		On the other hand, by \eqref{eq:eigen-Hs}, Lemma \ref{lem:sharp-low-spec} and the fact $s\ge 1$, we obtain, for every $0<t<T$
		\begin{align*}
			\|v_s(T)\|_2^2&=\|P_{N,s}e^{-(T-t)H_s}e^{-tH_s}g\|_2^2\le e^{-2(T-t)ns}\|v_s(t)\|_2^2\\
			&\le C_2e^{-\big(2(T-t)n-b^2\big)s}\|v_s(t)\|_{L^2(E_1)}^2,
		\end{align*}
		where $C_2>0$ is independent of $s$ and $t$. By \eqref{eq:time-critical-choice}, we have 
		\begin{equation*}
			2(T-t)n-b^2>0
		\end{equation*}
		for all $t\in[t_1,t_2]$.
		Thus,
		\begin{equation*}
			\|v_s(T)\|_2^2\le C_2\|v_s(t)\|_{L^2(E_1)}^2,\qquad t\in[t_1,t_2].
		\end{equation*} 
		Integrating the above inequality over $t\in[t_1,t_2]$ gives
		\begin{equation}\label{eq:low-mode}
			\|v_s(T)\|_2^2 \le C_2' \int_{t_1}^{t_2} \|v_s(t)\|_{L^2(E_1)}^2\,dt.
		\end{equation}
		
		Then, by triangle inequality, \eqref{eq:low-mode} and \eqref{eq:high-mode-uniform}, we have
		\begin{align*}
			&\quad\|u_s(T)\|_2^2\le 2\|v_s(T)\|_2^2+2\|w_s(T)\|_2^2\\
			&\le 4C_2'\int_{t_1}^{t_2} \|u_s(t)\|_{L^2(E_1)}^2\,dt+4C_2'\int_{t_1}^{t_2} \|w_s(t)\|_2^2\,dt+2\|w_s(T)\|_2^2\\
			&\le 4C_2'\int_{t_1}^{t_2} \|u_s(t)\|_{L^2(E_1)}^2\,dt+C_3\int_0^{t_0}\|u_s(\tau)\|_{L^2(E_1)}^2\,d\tau\\
			&\le C\int_0^T\|u_s(t)\|_{L^2(E_1)}^2\,dt,
		\end{align*}
		for all $T>\frac{b^2}{2n}$. For any given $T > \frac{a^2}{2n}$, by choosing $b$ such that $a < b < \sqrt{2nT}$, we conclude that \eqref{eq:uniform-fiber-critical} holds for all $s \ge 1$.

		Next, we consider $0<s<1$. It follows from \eqref{eq:spec-case1} that for every $g\in \text{Ran}(P_\lambda(H_s))$,
		\begin{equation*}\label{equ-spec-H1}
			\|g\|_{L^2(\R^n)}\le Ce^{C(1+\sqrt{\lambda})}\|g\|_{L^2(E_1)}.
		\end{equation*}
		On the other hand, we have the dissipation estimate
		\begin{equation*}
			\|(1-P_\lambda)(e^{-tH_s}g)\|_{L^2(\R^n)}\le e^{-t\lambda}\|g\|_{L^2(\R^n)}.
		\end{equation*}
		Combining these with the Lebeau--Robbiano strategy (see, e.g Theorem \ref{thm-L-R-method}), we obtain, for all $0<s<1,\,T>0$ and $g\in L^2(\R^n)$,
		\begin{equation}\label{equ-obs-hermite2}
			\|e^{-TH_{s}}g\|_{L^2(\R^n)}^2\le Ce^{CT^{-1}}\int_0^T\|e^{-tH_{s}}g\|_{L^2(E_1)}^2dt,
		\end{equation}
		with $C>0$ independent of $s$. Thus, \eqref{eq:uniform-fiber-critical} holds for all $T>0$ provided $0<s<1$.
		
		Therefore, the proof of Lemma \ref{lem:uniform-obs-Hs} is complete.
	\end{proof}

	We now prove that \eqref{equ-obs-Tlarge} holds for every $T>\frac{a^2}{2n}$. Since $\mathcal{F}_y(-\Delta_1)=-\Delta+|\eta|^2|x|^2$, applying Plancherel's theorem together with Lemma~\ref{lem:uniform-obs-Hs} yields, for any $T > \frac{a^2}{2n}$,
	\begin{align*}
		\|e^{T\Delta_1}f\|_{L^2(\R
			^n\times\R^m)}^2&=\int_{\R^m_\eta}\|e^{-TH_{1,\eta}}\hat{f}(x,\eta)\|_{L^2(\R^n_x)}^2d\eta\\
		&\le C_T \int_0^T \int_{\R^m_\eta}\|e^{-tH_{1,\eta}}\hat{f}(x,\eta)\|_{L^2(E_1)}^2d\eta dt\\
		&=C_T \int_0^T\|e^{t\Delta_1}f\|_{L^2(E_1\times\R^m)}^2dt.
	\end{align*}
	This completes the proof of \eqref{equ-obs-Tlarge}.

	The Hilbert Uniqueness Method now gives null controllability from $E_1\times\mathbb R^m$ for every $T>a^2/(2n)$. Combined with the necessity argument, this proves the equivalence in Theorem~\ref{thm:mini-timel}. In particular,
	$$
	T_{\min}(E_1\times\mathbb R^m)=\frac{a^2}{2n}.
	$$

	\section{Proof of Theorem \ref{thm:strong}}\label{sec5}
	Throughout this section, let $n=1$, $m\ge1$, and $\gamma>1$, and write $E=E_1\times\mathbb R^m$, where $E_1\subset\mathbb R$ is thick.
	
	\subsection{Proof of Theorem~\ref{thm:strong}(i): positive essential distance}
	
	Assume that $a_{\mathrm{ess}}(E_1)>0$. Choose $a_0$ with $0<a_0<a_{\mathrm{ess}}(E_1)$; then $E_1\subset\{|x|\ge a_0\}$ up to a null set. Let $\phi_{\gamma,1}$ be the normalized ground state of $H_\gamma=-\partial_x^2+|x|^{2\gamma}$, with eigenvalue $\lambda_{\gamma,1}>0$. For $n=1$, the unitary dilation $\mathcal M_{\gamma,\eta}$ is given by
	$$
	(\mathcal M_{\gamma,\eta}g)(x) =|\eta|^{\frac{1}{2(1+\gamma)}} g\bigl(|\eta|^{\frac{1}{1+\gamma}}x\bigr), \qquad g\in L^2(\mathbb R).
	$$
	Taking $j=1$ in the scaling relation \eqref{equ-eigen-connect}, we see that $\mathcal M_{\gamma,\eta}\phi_{\gamma,1}$ is the normalized ground state of $H_{\gamma,\eta}$, with eigenvalue $\lambda_{\gamma,1}|\eta|^{2/(1+\gamma)}$.
	
	If null controllability held in some time $T>0$, then there exists some constant $C_{\mathrm{obs}}>0$ such that
	\begin{equation}\label{eq:strong-obs-test}
		\|e^{T\Delta_\gamma}f\|_{L^2(\mathbb R^{1+m})}^2 \le C_{\mathrm{obs}}\int_0^T\|e^{t\Delta_\gamma}f\|_{L^2(E_1\times\mathbb R^m)}^2\,dt
	\end{equation}
	for all $f\in L^2(\mathbb R^{1+m})$. For $k\ge1$, set
	$$
	S_k=\{\eta\in\mathbb R^m:k<|\eta|<k+1\},\qquad \widehat f_k(x,\eta)=\mathcal M_{\gamma,\eta}\phi_{\gamma,1}(x)\,\mathbf 1_{S_k}(\eta).
	$$
	For each $\eta\in S_k$, note that
	$$
	e^{-TH_{\gamma,\eta}}\mathcal M_{\gamma,\eta}\phi_{\gamma,1} =e^{-T\lambda_{\gamma,1}|\eta|^{2/(1+\gamma)}} \mathcal M_{\gamma,\eta}\phi_{\gamma,1}.
	$$
	Thus, by Plancherel's theorem, the normalization of the ground state, and $|\eta|<k+1$ on $S_k$, we derive
	\begin{equation}\label{eq:strong-obs-test-final}
		\|e^{T\Delta_\gamma}f_k\|_2^2 \ge |S_k|\,\exp{\big(-2T\lambda_{\gamma,1}(k+1)^{2/(1+\gamma)}\big)}.
	\end{equation}
	For $V_\gamma(z)=|z|^{2\gamma}$, the Agmon distance from a fixed compact set to $z$ grows like
	$$
	\int_{\lambda_{\gamma,1}^{1/(2\gamma)}}^{|z|}\sqrt{t^{2\gamma}-\lambda_{\gamma,1}}\,dt \sim \frac{|z|^{1+\gamma}}{1+\gamma} \qquad(|z|\to\infty).
	$$
	Thus the classical Agmon estimate (see, e.g., \cite{Agmon82}) gives some $c_0>0$ such that
	$$
	\int_{\mathbb R}e^{2c_0|z|^{1+\gamma}}|\phi_{\gamma,1}(z)|^2\,dz<\infty.
	$$
	In particular, for every $R>0$,
	$$
	\int_{|z|\ge R}|\phi_{\gamma,1}(z)|^2\,dz \le C e^{-2c_0R^{1+\gamma}}.
	$$
	Now set $z=|\eta|^{1/(1+\gamma)}x$ in the scaling formula for $\mathcal M_{\gamma,\eta}$. The normalization factor cancels the Jacobian, so the preceding tail estimate yields
	$$
	\int_{|x|\ge a_0}|\mathcal M_{\gamma,\eta}\phi_{\gamma,1}(x)|^2\,dx=\int_{|z|\ge a_0|\eta|^{1/(1+\gamma)}}|\phi_{\gamma,1}(z)|^2\,dz\le C e^{-c a_0^{1+\gamma}|\eta|}.
	$$
	Since for each $t\ge0$,  Plancherel's theorem in $y$ gives
	$$
	\begin{aligned}
		\|e^{t\Delta_\gamma}f_k\|_{L^2(E_1\times\mathbb R^m)}^2= \int_{S_k}e^{-2t\lambda_{\gamma,1}|\eta|^{2/(1+\gamma)}}\|\mathcal M_{\gamma,\eta}\phi_{\gamma,1}\|_{L^2(E_1)}^2\,d\eta.
	\end{aligned}
	$$
	Applying the Agmon tail estimate to the $E_1$-norm in the preceding identity, and using that $E_1\subset\{|x|\ge a_0\}$ up to a null set, $|\eta|>k$ on $S_k$, and $e^{-2t\lambda_{\gamma,1}|\eta|^{2/(1+\gamma)}}\le1$, we obtain, uniformly for $t\ge0$,
	\begin{equation}\label{eq:strong-obs-test-control}
		\int_0^T\|e^{t\Delta_\gamma}f_k\|_{L^2(E_1\times\mathbb R^m)}^2\,dt\le CT|S_k|e^{-c a_0^{1+\gamma}k}.
	\end{equation}
	Combining \eqref{eq:strong-obs-test}--\eqref{eq:strong-obs-test-control} and cancelling $|S_k|$ yields
	$$
	1\le C_{\mathrm{obs}}CT \exp\!\left(-c a_0^{1+\gamma}k+2T\lambda_{\gamma,1}(k+1)^{2/(1+\gamma)}\right).
	$$
	Because $\gamma>1$, we have $2/(1+\gamma)<1$, so the right-hand side tends to zero as $k\to\infty$. This contradiction proves part (i).
	
	\subsection{Proof of Theorem~\ref{thm:strong}(ii)}
	
	Assume henceforth that $a_{\mathrm{ess}}(E_1)=0$. For $r>0$, abbreviate $D(r)=D_{E_1}(r)$, as defined in \eqref{eq:density}. We first prove the following scaled observability inequality.
	
	\begin{lemma}\label{lem:scaling-obs}
		Let $\gamma>1$, $T>0$ and let $E\subset\mathbb R$ satisfy $a_{\mathrm{ess}}(E)=0$. For $r>0$, set $\delta_E(r)=|E\cap(-r,r)|/(2r)$ and $D_E(r)=\log(1/\delta_E(r))$.
		Fix $0<\varepsilon<(\gamma-1)/(2\gamma)$ and set
		$$
		\beta=\frac{1}{2\gamma}+\frac12+\varepsilon<1.
		$$
		Then, for every $r>0$ and $f\in L^2(\mathbb R)$,
		\begin{equation}\label{eq:obs-E-over-r}
			\|e^{-T H_\gamma}f\|_2^2 \le K_T(r)\int_0^T \|e^{-tH_\gamma}f\|_{L^2(r^{-1}E)}^2\,dt,
		\end{equation}
		where
		\begin{equation}\label{eq:K-Tr}
			K_T(r)=\frac{C}{T} \exp\!\left(\frac{C(1+D_E(r))^{1/(1-\beta)}}{T^{\beta/(1-\beta)}}\right).
		\end{equation}
		The constant $C$ may depend on $\gamma$ and $\varepsilon$, but is independent of $r$ and $T$.
	\end{lemma}
	
	\begin{proof}
		Set $F_r=(r^{-1}E)\cap(-1,1)$. By the definition of $\delta_E$,
		$$
		|F_r|=2\delta_E(r)>0.
		$$
		For a spectral parameter $\Lambda\ge1$, let $\widehat I_\Lambda=[-\widehat C\Lambda^{1/(2\gamma)},\widehat C\Lambda^{1/(2\gamma)}]$, with $\widehat C$ as in \cite[Corollary 1]{MZ25}, and rescale $F_r$ to $\widetilde F_{r,\Lambda}=|\widehat I_\Lambda|^{-1}F_r$. Thus $|\widetilde F_{r,\Lambda}|\asymp\Lambda^{-1/(2\gamma)}\delta_E(r)$. The spectral inequality in the proof of \cite[Theorem 3]{MZ25}, applied to $\widetilde H=-\partial_x^2+1+|x|^{2\gamma}$, gives
		$$
		\|u\|_2\le e^{C\Lambda^{1/(2\gamma)+1/2}} \|u\|_{L^2(F_r)}^\alpha\|u\|_2^{1-\alpha}, \qquad u\in\operatorname{Ran}(P_\Lambda(\widetilde H)),
		$$
		where
		$$
		\alpha^{-1}\le C\bigl(1+\log\Lambda+D_E(r)\bigr).
		$$
		Here we used $|F_r|=2\delta_E(r)$ and $D_E(r)=\log(1/\delta_E(r))$. Dividing by $\|u\|_2^{1-\alpha}$ (the case $u=0$ is immediate) and raising to the power $1/\alpha$, then absorbing the logarithmic factor into $\Lambda^\varepsilon$, yields
		$$
		\|u\|_2\le \exp\!\big(C(1+D_E(r))\Lambda^\beta\big) \|u\|_{L^2(r^{-1}E)},\qquad u\in\operatorname{Ran}(P_\Lambda(H_\gamma)).
		$$
		We used $\widetilde H=H_\gamma+1$ to pass between the two spectral projections. The Lebeau--Robbiano strategy (see \cite[Theorem 2.12]{NTTV20}) applied to this spectral inequality gives \eqref{eq:obs-E-over-r}--\eqref{eq:K-Tr}.
	\end{proof}
	
	\subsubsection{Reduction to uniform observability inequalities.}
	
	For $|\eta|\ge1$, put
	$$
	r=|\eta|^{-1/(1+\gamma)}\in(0,1].
	$$
	The unitary scaling of $H_{\gamma,\eta}$ gives
	$$
	\mathcal M_{\gamma,\eta}^*H_{\gamma,\eta}\mathcal M_{\gamma,\eta} =r^{-2}H_\gamma,\qquad \|\mathcal M_{\gamma,\eta}u\|_{L^2(E_1)}=\|u\|_{L^2(r^{-1}E_1)}.
	$$
	Consequently, \eqref{eq:obs-E-over-r}, after the change of time variable, implies
	\begin{equation}\label{eq:obs-large-eta}
		\|e^{-TH_{\gamma,\eta}}f\|_2^2 \le r^{-2}K_{T/r^2}(r)\int_0^T \|e^{-tH_{\gamma,\eta}}f\|_{L^2(E_1)}^2\,dt.
	\end{equation}
	For $0<|\eta|<1$, the spectral inequality \eqref{eq:spec-case1} and the Lebeau--Robbiano strategy yield, uniformly in $\eta$,
	\begin{equation}\label{eq:obs-small-eta}
		\|e^{-TH_{\gamma,\eta}}f\|_2^2 \le C e^{C/T}\int_0^T\|e^{-tH_{\gamma,\eta}}f\|_{L^2(E_1)}^2\,dt.
	\end{equation}
	
	The bottom of the spectrum of $H_{\gamma,\eta}$ is $\lambda_{\gamma,1}/r^2$. Splitting the semigroup at time $T/2$ and using \eqref{eq:obs-large-eta}, we obtain, for $|\eta|\ge1$,
	$$
	\|e^{-TH_{\gamma,\eta}}f\|_2^2\le \mathcal A_T(r)\int_0^T\|e^{-tH_{\gamma,\eta}}f\|_{L^2(E_1)}^2\,dt,
	$$
	where
	\begin{align}
		\mathcal A_T(r)
		&:=r^{-2}K_{T/(2r^2)}(r)e^{-T\lambda_{\gamma,1}/r^2}\notag\\
		&=\frac{C}{T}\exp\!\left[\frac{1}{r^2}\left\{\frac{C\bigl(r^2(1+D(r))\bigr)^{1/(1-\beta)}}{T^{\beta/(1-\beta)}}-T\lambda_{\gamma,1}\right\}\right].\label{eq:strong-fiber-cost}
	\end{align}
	Thus, if $\mathcal A_T$ is bounded on $(0,1]$, then Plancherel's theorem and \eqref{eq:obs-small-eta} give the global observability estimate
	\begin{equation}\label{eq:strong-global-obs}
		\|e^{T\Delta_\gamma}g\|_{L^2(\mathbb R^{1+m})}^2 \le C_T\int_0^T \|e^{t\Delta_\gamma}g\|_{L^2(E_1\times\mathbb R^m)}^2\,dt.
	\end{equation}
	
	\subsubsection{Proof of Theorem~\ref{thm:strong}(ii)(a): $L_{E_1}=0$}
	
	If $L_{E_1}=0$, then $r^2(1+D(r))\to0$ as $r\downarrow0$. For each fixed $T>0$, the expression in braces in \eqref{eq:strong-fiber-cost} is therefore at most $-T\lambda_{\gamma,1}/2$ for all sufficiently small $r$. Hence $\mathcal A_T(r)$ is bounded near zero. On every compact subinterval of $(0,1]$, the function $D$ is continuous and finite, so $\mathcal A_T$ is bounded there as well. Thus \eqref{eq:strong-global-obs} holds for every $T>0$. The Hilbert Uniqueness Method gives null controllability in every positive time, proving (a).
	
	\subsubsection{Proof of Theorem~\ref{thm:strong}(ii)(b): $0<L_{E_1}<\infty$}
	
	If $0<L_{E_1}<\infty$, continuity of $D$ away from zero and the definition of $L_{E_1}$ imply
	$$
	M_E:=\sup_{0<r\le1}r^2(1+D(r))<\infty.
	$$
	Choose $T_{\mathrm{suff}}>0$ sufficiently large that
	$$
	\frac{C M_E^{1/(1-\beta)}}{T^{\beta/(1-\beta)}} \le \frac{T\lambda_{\gamma,1}}{2} \qquad\text{for every }T>T_{\mathrm{suff}}.
	$$
	Then \eqref{eq:strong-fiber-cost} gives $\sup_{0<r\le1}\mathcal A_T(r)<\infty$ for every $T>T_{\mathrm{suff}}$. By \eqref{eq:strong-global-obs} and the Hilbert Uniqueness Method, null controllability holds for all such $T$. This proves (b).
	
	\subsubsection{Proof of Theorem~\ref{thm:strong}(ii)(c): necessary density bound}
	
	Assume that null controllability holds in some time $T>0$. The corresponding observability estimate has a constant $C_{\mathrm{obs}}$. Let $\phi_{\gamma,1}$ again denote the normalized ground state of $H_\gamma$. For $k\ge1$, choose $f_k$ by
	$$
	\widehat f_k(x,\eta) =\mathcal M_{\gamma,\eta}\phi_{\gamma,1}(x)\,\mathbf 1_{\{k<|\eta|<k+1\}}(\eta).
	$$
	Writing $p=2/(1+\gamma)$ and
	$V_{k,m}=|\{\eta\in\mathbb R^m:k<|\eta|<k+1\}|$, Plancherel's theorem gives
	\begin{equation}\label{eq:strong-density-final}
		\|e^{T\Delta_\gamma}f_k\|_2^2 \ge V_{k,m}e^{-2T\lambda_{\gamma,1}(k+1)^p}.
	\end{equation}
	For $s>0$, set $r=|\eta|^{-1/(1+\gamma)}$ and write $z=x/r$ for the scaled spatial variable. The scaling formula gives
	\[
	\|\mathcal M_{\gamma,\eta}\phi_{\gamma,1}\|_{L^2(E_1)}^2 =\int_{r^{-1}E_1}|\phi_{\gamma,1}(z)|^2\,dz.
	\]
	We split this integral into the regions $|z|\le s/r$ and $|z|>s/r$; this gives
	\begin{align}
		\|\phi_{\gamma,1}\|_{L^2(r^{-1}E_1)}^2
		&\le \frac{2s\delta_{E_1}(s)}{r} \|\phi_{\gamma,1}\|_{L^\infty(\mathbb R)}^2 +C_\gamma e^{-c_\gamma(s/r)^{1+\gamma}}. \label{eq:strong-density-mass}
	\end{align}
	Indeed, the inner region has measure $|(r^{-1}E_1)\cap(-s/r,s/r)|=r^{-1}|E_1\cap(-s,s)|=2s\delta_{E_1}(s)/r$, so the integral over this region is bounded by the first term. The outer region is bounded by the Agmon tail at radius $s/r$, giving the second term.
	
	On the shell $k<|\eta|<k+1$, we have $r^{-1}\le(k+1)^{1/(1+\gamma)}$ and $(s/r)^{1+\gamma}\ge s^{1+\gamma}k$. Applying observability to $f_k$, using \eqref{eq:strong-density-final} and \eqref{eq:strong-density-mass}, and cancelling $V_{k,m}$, yields
	\begin{equation}\label{eq:strong-density-inequality}
		e^{-2T\lambda_{\gamma,1}(k+1)^p} \le C_T\left[s\delta_{E_1}(s)(k+1)^{1/(1+\gamma)}+e^{-c_\gamma s^{1+\gamma}k}\right].
	\end{equation}
	Let
	$$
	q=\frac{(1+\gamma)^2}{\gamma-1},\qquad k_s=A_{T,\gamma}s^{-q},
	$$
	where $A_{T,\gamma}$ is large enough that, for all sufficiently small $s>0$,
	$$
	c_\gamma s^{1+\gamma}k_s \ge 4T\lambda_{\gamma,1}(k_s+1)^p.
	$$
	The exponential remainder in \eqref{eq:strong-density-inequality} is then negligible compared with its left-hand side. It follows that
	$$
	\delta_{E_1}(s) \ge c_Ts^{-1}(k_s+1)^{-1/(1+\gamma)} e^{-2T\lambda_{\gamma,1}(k_s+1)^p}.
	$$
	Taking logarithms and using $q p=2(1+\gamma)/(\gamma-1)$, we obtain
	$$
	D_{E_1}(s)\le C_{T,\gamma}s^{-2(1+\gamma)/(\gamma-1)}
	$$
	for all sufficiently small $s>0$. Therefore
	$$
	\limsup_{s\downarrow0}s^{2(1+\gamma)/(\gamma-1)}D_{E_1}(s)<\infty,
	$$
	which is precisely the necessary condition in Theorem~\ref{thm:strong}(ii)(c). The proof of Theorem~\ref{thm:strong} is complete.

	\appendix
	\section{Failure of null controllability for $\gamma>1$}\label{app-1}
	\begin{proof}[Proof of Proposition~\ref{prop:transverse-gap}]
		Let $n=m=1$ and $E=\mathbb R\times E_2$. Since $\overline{E_2}\ne\mathbb R$, it follows that $\mathbb R\setminus{E_2}$ contains an open interval. 
        By the translation invariance of the equation in the $y$-variable, we may assume, without loss of generality, that
 $$
 E_2\cap(-a,a)=\varnothing\qquad \mbox{for some}\,\,a>0.
 $$ 
 We show that null controllability from $E$ fails for every $T>0$. By the Hilbert Uniqueness Method, it suffices to construct initial data whose solutions violate the observability estimate. We therefore seek a family $f_{0,h}(x,y)$ such that the quotient
		\begin{equation}\label{equ-quotient}
			\frac{\int_{0}^{T}\int_\R\int_{E_2}|f_h(t,x,y)|^2dxdydt}{\int_\R\int_{\R}|f_h(T,x,y)|^2dxdy} \to 0,\quad\mbox{as }h\to 0.
		\end{equation}
		
		We consider $f_{0,h}(x,y)\in L^2(\R\times \R)$ defined by
		\begin{equation}\label{equ-exam-data}
			f_{0,h}(x,y)=(2\pi)^{-1/2}h\int_{\R}\chi(h\eta-\eta_0)e^{-(h\eta-\eta_0)^2/{2h}+iy\eta} v_x(\eta)d\eta,
		\end{equation}
		where $v_x(\eta):=\phi_{\gamma,1}\big(|\eta|^{\frac{1}{1+\gamma}}x\big)$, $\eta_0>0$ is chosen large enough, and $\chi\in C_0^\infty(-\delta,\delta)$ with $\delta>0$ sufficiently small, satisfying $0\le \chi\le 1$ and $\chi\equiv1$ in a neighborhood of $0$. 
		By \eqref{equ-diag-H-eta}, \eqref{equ-eigen-connect}, and \eqref{equ-exam-data}
		\begin{align}\label{equ-solution}
			f_h(t,x,y)&= e^{-t(-\Delta_\gamma)}f_{0,h}(x,y)= \mathcal{F}_\eta^{-1}\big(e^{-tH_{\gamma,\eta}}\mathcal{F}_y(f_{0,h})(x,\eta)\big)(t,x,y)\notag\\
			&=(2\pi)^{-1/2}h\int_\R\chi(h\eta-\eta_0)e^{-(h\eta-\eta_0)^2/{2h}+iy\eta-t\lambda_{\gamma,1}|\eta|^{\frac{2}{1+\gamma}}} v_x(\eta)d\eta.
		\end{align}
		Then, by Plancherel's theorem and the $L^2$-normalization of $\phi_{\gamma,1}$,
		\begin{align}\label{equ-fh-norm}
			&\quad\|f_h(T,x,y)\|_{L^2(\R^2)}^2
			=\int_{\R_\eta}\big|h\chi(h\eta-\eta_0)e^{-(h\eta-\eta_0)^2/{2h}-T\lambda_{\gamma,1}|\eta|^{\frac{2}{1+\gamma}}} |\eta|^{-\frac{1}{2(1+\gamma)}}\big|^2d\eta\notag\\
			&\ge \int_{|h\eta-\eta_0|\le \sqrt{h}}  \big|h\chi(h\eta-\eta_0)e^{-(h\eta-\eta_0)^2/{2h}-T\lambda_{\gamma,1}|\eta|^{\frac{2}{1+\gamma}}} |\eta|^{-\frac{1}{2(1+\gamma)}}\big|^2d\eta\notag\\
			&\ge Ch^{\frac{3}{2}+\frac{1}{1+\gamma}}\exp\big(-Ch^{-\frac{2}{1+\gamma}}\big),
		\end{align}
		where $C>0$ depends on $T,\gamma$.
		
		We now estimate the numerator of \eqref{equ-quotient} by using Proposition 3.7 in \cite{Koe20}.
		It follows from \eqref{equ-solution} that
		\begin{align}\label{equ-solution-upp}
			f_h(t,x,y)
			=(2\pi)^{-1/2}\int_\R\chi(\eta-\eta_0)e^{-(\eta-\eta_0)^2/{2h}+iy\eta/h- t\lambda_{\gamma,1}|\eta/h|^{\frac{2}{1+\gamma}}} v_x(\eta/h)d\eta.
		\end{align}
		Since $\eta_0\gg1$, the cutoff $\chi$ restricts the frequency to the positive half-line.
		All powers of a complex variable below are taken with the principal branches.
		Since the ground state $\phi_{\gamma,1}$ is even, the corresponding continuation of $v_x$ is
		$$
		v_x(z)= \phi_{\gamma,1}\bigl(z^{1/(1+\gamma)}|x|\bigr),\quad v_0(z)=\phi_{\gamma,1}(0).
		$$
		By Sibuya's theorem in \cite[Section 2]{EG18}, $v_x$ is holomorphic on the right half-plane, and
		$$
		\phi_{\gamma,1}(w) = c_\gamma w^{-\gamma/2} e^{-\frac{w^{\gamma+1}}{\gamma+1}}(1+o(1))
		$$
		uniformly as $|w|\to\infty$ with $|\arg w|\le\pi/(1+\gamma)$.

		For $|x|>1$, it follows that
		\begin{equation}\label{eq:v_x-large}
			v_x(z/h)=e^{-\frac{|x|^{\gamma+1}}{\gamma+1}(z/h)}\tilde v_x(z/h),\qquad\tilde v_x(z/h)=c_\gamma(z/h)^{-\frac{\gamma}{2(1+\gamma)}}|x|^{-\gamma/2}(1+o(1)),
		\end{equation}
		uniformly in $|x|>1$ as $|z/h|\to\infty$ with
		$|\arg z|\le\pi/(1+\gamma)$.
		Choose $K$ sufficiently large and set $\mathcal C=\{z\in\mathbb C:\Re z>K,\ |\Im z|<K^{-1}\Re z\}.$
		The uniform asymptotics give
		$$
		|\tilde v_x(z/h)|\le C\big(\frac{h}{|z|}\big)^{\frac{\gamma}{2(1+\gamma)}}|x|^{-\gamma/2}\le C,\qquad z\in\mathcal C,\quad |x|>1,\quad 0<h\le h_0,
		$$
		where $C$ is independent of $h$, $x$, and $z$.
		For $|x|\le1$, boundedness near the origin and the above sectorial asymptotics yield
		\begin{equation}\label{eq:v_x-small}
			\sup_{|x|\le1}\sup_{z\in\mathcal C}
			|v_x(z/h)|= \sup_{|x|\le1}\sup_{z\in\mathcal C}\left|\phi_{\gamma,1}\bigl((z/h)^{1/(1+\gamma)}|x|\bigr)\right|\le C.
		\end{equation}
		Thus, by \eqref{eq:v_x-large} and \eqref{eq:v_x-small}, we have
		\begin{align}\label{eq:v_x-infty}
			\|v_x(z/h)\|_{L^\infty(\mathcal C)}\le
			\begin{cases}
				C,&|x|\le1,\\
				Ce^{-\frac{|x|^{\gamma+1}}{\gamma+1}Kh^{-1}},&|x|>1.
			\end{cases}
		\end{align}
		Since $z\mapsto v_x(z/h)$ is holomorphic on $\mathcal C$ and satisfies the uniform bounds above, Proposition 3.7 in \cite{Koe20} yields
		\begin{equation}\label{eq:f_h-bound}
			|f_h(t,x,y)|\le \frac{C}{|y|^2}e^{-c/h+\varepsilon_{t}(h)/h}\,\|v_x(z/h)\|_{L^\infty(\mathcal C)},
		\end{equation}
		uniformly for $t\in[0,T]$, $x\in\R$, $|y|\ge a$, and $0<h\le h_0$.
		Here,
		\begin{align}\label{eq:eps-bound}
			\frac{\varepsilon_{t}(h)}{h}&=\sup_{|\eta|<\delta,|y|<\delta}\Re\big(-t\lambda_{\gamma,1}(\frac{\eta+\eta_0+iy}{h})^{2/(1+\gamma)}\big)\notag\\
			&\le\sup_{|\eta|<\delta,|y|<\delta}\Big(-t\lambda_{\gamma,1}\frac{(\eta+\eta_0)^{2/(\gamma+1)}}{h^{2/(1+\gamma)}}\Big)\notag\\
			&\le -Ct\lambda_{\gamma,1}h^{-\frac{2}{1+\gamma}}.
		\end{align}
		By \eqref{eq:v_x-infty}, \eqref{eq:f_h-bound}, and \eqref{eq:eps-bound},
		\begin{align*}
			&\quad\int_{0}^{T}\int_\R\int_{E_2}|f_h(t,x,y)|^2dydxdt\\
			&\le \int_{0}^{T}\int_{|x|\le 1}\int_{|y|\ge a}\frac{C}{|y|^4}e^{-c/h-Ct\lambda_{\gamma,1}h^{-\frac{2}{1+\gamma}}}dydxdt\\
			&\quad\qquad+\int_{0}^{T}\int_{|x|> 1}\int_{|y|\ge a}\frac{C}{|y|^4}e^{-c/h-Ct\lambda_{\gamma,1}h^{-\frac{2}{1+\gamma}}-\frac{|x|^{\gamma+1}}{\gamma+1}Kh^{-1}}dydxdt\\
			&\le C_{T,a,\gamma}\ e^{-c/h}\Big(1+\int_{\R} e^{-\frac{|x|^{\gamma+1}}{\gamma+1}Kh^{-1}}dx\Big)\\
			&< C_{T,a,\gamma}\ e^{-c/h}(1+Ch^{1/(\gamma+1)})<Ce^{-c/h}.
		\end{align*}
		This, combining with \eqref{equ-fh-norm} and the fact $\gamma>1$, implies \eqref{equ-quotient}. Thus Proposition~\ref{prop:transverse-gap} is proved.
	\end{proof}

	\section*{Acknowledgments}
	S. H. was supported by the National Natural Science Foundation of China under Grant 12622107 and by the Guangdong Basic and Applied Basic Research Foundation under Grant 2026B1515020075.
G. W. was partially supported by the New Cornerstone Science Foundation and the National Natural Science Foundation of China under Grant 12371450.

    \section*{Disclosure of AI assistance}
	The authors began investigating the problem studied in this paper in July 2025. The proof framework was developed by the authors, who also wrote the initial manuscript. During drafting and revision, the authors used AI tools, including GPT-5.6 and GPT-6, to varying extents. 
	AI assistance was used mainly to improve grammar and exposition, clarify and supplement proof details, and identify possible errors. 

	\bibliographystyle{plain}
	\bibliography{grushinRef}

@book{Agmon82,
	URL = {http://www.jstor.org/stable/j.ctt13x1d8z},
	author = {Agmon, S.},
	publisher = {Princeton University Press},
	title = {Lectures on Exponential Decay of Solutions of Second-Order Elliptic Equations: Bounds on Eigenfunctions of N-Body Schrodinger Operations. (MN-29)},
	urldate = {2026-05-06},
	year = {1982}
}

@article {AS25,
	AUTHOR = {Alphonse, P. and Seelmann, A.},
	TITLE = {Unique continuation estimates for {B}aouendi-{G}rushin equations on cylinders},
	JOURNAL = {Pure Appl. Anal.},
	FJOURNAL = {Pure and Applied Analysis},
	VOLUME = {7},
	YEAR = {2025},
	NUMBER = {3},
	PAGES = {733--770},
	ISSN = {2578-5893,2578-5885},
	MRCLASS = {35P05 (35P10 93B05 93B07)},
	MRNUMBER = {4947067},
	DOI = {10.2140/paa.2025.7.733},
	URL = {https://doi.org/10.2140/paa.2025.7.733},
}

@article {BCG14,
	AUTHOR = {Beauchard, K. and Cannarsa, P. and Guglielmi, R.},
	TITLE = {Null controllability of {G}rushin-type operators in dimension two},
	JOURNAL = {J. Eur. Math. Soc. (JEMS)},
	FJOURNAL = {Journal of the European Mathematical Society (JEMS)},
	VOLUME = {16},
	YEAR = {2014},
	NUMBER = {1},
	PAGES = {67--101},
	ISSN = {1435-9855,1435-9863},
	MRCLASS = {93B05 (35K65 93B07 93C20)},
	MRNUMBER = {3141729},
	MRREVIEWER = {Angelo\ Favini},
	DOI = {10.4171/JEMS/428},
	URL = {https://doi.org/10.4171/JEMS/428},
}

@article {BPS18,
	AUTHOR = {Beauchard, K. and Pravda-Starov, K.},
	TITLE = {Null-controllability of hypoelliptic quadratic differential equations},
	JOURNAL = {J. \'Ec. polytech. Math.},
	FJOURNAL = {Journal de l'\'Ecole polytechnique. Math\'ematiques},
	VOLUME = {5},
	YEAR = {2018},
	PAGES = {1--43},
	ISSN = {2429-7100,2270-518X},
	MRCLASS = {93B05 (35H10 35K65 93C20)},
	MRNUMBER = {3732691},
	MRREVIEWER = {El\ Mustapha\ Ait Ben Hassi},
	DOI = {10.5802/jep.62},
	URL = {https://doi.org/10.5802/jep.62},
}

@article {DSV24,
	AUTHOR = {Dicke, A. and Seelmann, A. and Veseli\'c, I.},
	TITLE = {Spectral inequality with sensor sets of decaying density for {S}chr\"odinger operators with power growth potentials},
	JOURNAL = {Partial Differ. Equ. Appl.},
	FJOURNAL = {Partial Differential Equations and Applications},
	VOLUME = {5},
	YEAR = {2024},
	NUMBER = {2},
	PAGES = {Paper No. 7, 18},
	ISSN = {2662-2963,2662-2971},
	MRCLASS = {35Pxx (35B40 35J10)},
	MRNUMBER = {4725290},
	DOI = {10.1007/s42985-024-00276-0},
	URL = {https://doi.org/10.1007/s42985-024-00276-0},
}

@article {HWW21,
	AUTHOR = {Huang, S. and Wang, G. and Wang, M.},
	TITLE = {Characterizations of stabilizable sets for some parabolic equations in {$\Bbb{R}^n$}},
	JOURNAL = {J. Differential Equations},
	FJOURNAL = {Journal of Differential Equations},
	VOLUME = {272},
	YEAR = {2021},
	PAGES = {255--288},
	ISSN = {0022-0396,1090-2732},
	MRCLASS = {93D15 (35K15 35Q41 35R11 93C20)},
	MRNUMBER = {4159677},
	MRREVIEWER = {Michela\ Egidi},
	DOI = {10.1016/j.jde.2020.09.038},
	URL = {https://doi.org/10.1016/j.jde.2020.09.038},
}

@article {Koe20,
	AUTHOR = {Koenig, A.},
	TITLE = {Lack of null-controllability for the fractional heat equation and related equations},
	JOURNAL = {SIAM J. Control Optim.},
	FJOURNAL = {SIAM Journal on Control and Optimization},
	VOLUME = {58},
	YEAR = {2020},
	NUMBER = {6},
	PAGES = {3130--3160},
	ISSN = {0363-0129,1095-7138},
	MRCLASS = {93B05 (35K65 93B07 93C20)},
	MRNUMBER = {4169745},
	MRREVIEWER = {El\ Mustapha\ Ait Ben Hassi},
	DOI = {10.1137/19M1256610},
	URL = {https://doi.org/10.1137/19M1256610},
}

@article {Kov01,
	AUTHOR = {Kovrijkine, O.},
	TITLE = {Some results related to the {L}ogvinenko-{S}ereda theorem},
	JOURNAL = {Proc. Amer. Math. Soc.},
	FJOURNAL = {Proceedings of the American Mathematical Society},
	VOLUME = {129},
	YEAR = {2001},
	NUMBER = {10},
	PAGES = {3037--3047},
	ISSN = {0002-9939,1088-6826},
	MRCLASS = {46E10 (42A16)},
	MRNUMBER = {1840110},
	DOI = {10.1090/S0002-9939-01-05926-3},
	URL = {https://doi.org/10.1090/S0002-9939-01-05926-3},
}

@article {Lis25,
	AUTHOR = {Lissy, P.},
	TITLE = {A non-controllability result for the half-heat equation on the whole line based on the prolate spheroidal wave functions and its application to the {G}rushin equation},
	JOURNAL = {J. Differential Equations},
	FJOURNAL = {Journal of Differential Equations},
	VOLUME = {433},
	YEAR = {2025},
	PAGES = {Paper No. 113306, 27},
	ISSN = {0022-0396,1090-2732},
	MRCLASS = {93B05 (35K15 35K65 35S10 93B07)},
	MRNUMBER = {4891568},
	MRREVIEWER = {Zhenbin\ Fan},
	DOI = {10.1016/j.jde.2025.113306},
	URL = {https://doi.org/10.1016/j.jde.2025.113306},
}

@article {LWXY22,
	AUTHOR = {Liu, H. and Wang, G. and Xu, Y. and Yu,
	H.},
	TITLE = {Characterizations of complete stabilizability},
	JOURNAL = {SIAM J. Control Optim.},
	FJOURNAL = {SIAM Journal on Control and Optimization},
	VOLUME = {60},
	YEAR = {2022},
	NUMBER = {4},
	PAGES = {2040--2069},
	ISSN = {0363-0129,1095-7138},
	MRCLASS = {93D15 (93B05 93C25 93D23)},
	MRNUMBER = {4450289},
	MRREVIEWER = {Hongyinping\ Feng},
	DOI = {10.1137/20M1386761},
	URL = {https://doi.org/10.1137/20M1386761},
}

@article {BJP21,
	AUTHOR = {Beauchard, K. and Jaming, P. and Pravda-Starov, K.},
	TITLE = {Spectral estimates for finite combinations of {H}ermite functions and null-controllability of hypoelliptic quadratic equations},
	JOURNAL = {Studia Math.},
	FJOURNAL = {Studia Mathematica},
	VOLUME = {260},
	YEAR = {2021},
	NUMBER = {1},
	PAGES = {1--43},
	DOI = {10.4064/sm191205-12-10},
	URL = {https://doi.org/10.4064/sm191205-12-10},
}

@article {TWX20,
	AUTHOR = {Tr\'elat, E. and Wang, G. and Xu, Y.},
	TITLE = {Characterization by observability inequalities of controllability and stabilization properties},
	JOURNAL = {Pure Appl. Anal.},
	FJOURNAL = {Pure and Applied Analysis},
	VOLUME = {2},
	YEAR = {2020},
	NUMBER = {1},
	PAGES = {93--122},
	ISSN = {2578-5893,2578-5885},
	MRCLASS = {93B05 (93B07 93C25)},
	MRNUMBER = {4041279},
	MRREVIEWER = {Yacine\ Chitour},
	DOI = {10.2140/paa.2020.2.93},
	URL = {https://doi.org/10.2140/paa.2020.2.93},
}

@article {WWZZ19,
	AUTHOR = {Wang, G. and Wang, M. and Zhang, C. and Zhang, Y.},
	TITLE = {Observable set, observability, interpolation inequality and spectral inequality for the heat equation in {$\Bbb{R}^n$}},
	JOURNAL = {J. Math. Pures Appl. (9)},
	FJOURNAL = {Journal de Math\'ematiques Pures et Appliqu\'ees. Neuvi\`eme
	S\'erie},
	VOLUME = {126},
	YEAR = {2019},
	PAGES = {144--194},
	ISSN = {0021-7824,1776-3371},
	MRCLASS = {93B05 (35K05 35K20 93B07)},
	MRNUMBER = {3950015},
	MRREVIEWER = {Larissa\ V.\ Fardigola},
	DOI = {10.1016/j.matpur.2019.04.009},
	URL = {https://doi.org/10.1016/j.matpur.2019.04.009},
}

@article {ZZG25,
	AUTHOR = {Zhu, J. and Zhuge, J.},
	TITLE = {Spectral inequality for {S}chr\"odinger equations with power growth potentials},
	JOURNAL = {Indiana Univ. Math. J.},
	FJOURNAL = {Indiana University Mathematics Journal},
	VOLUME = {74},
	YEAR = {2025},
	NUMBER = {5},
	PAGES = {1467--1494},
	ISSN = {0022-2518,1943-5258},
	MRCLASS = {35J10 (35P05)},
	MRNUMBER = {4995025},
	MRREVIEWER = {Weijia\ Wu},
	DOI = {10.1512/iumj.2025.74.60509},
	URL = {https://doi.org/10.1512/iumj.2025.74.60509},
}

@article {JW25,
	AUTHOR = {Jaming, P. and Wang, Y.},
	TITLE = {Null-controllability of the generalized {B}aouendi-{G}rushin
	heat like equations},
	JOURNAL = {J. Evol. Equ.},
	FJOURNAL = {Journal of Evolution Equations},
	VOLUME = {25},
	YEAR = {2025},
	NUMBER = {2},
	PAGES = {Paper No. 43, 39},
	ISSN = {1424-3199,1424-3202},
	MRCLASS = {93B05 (35J10 35K41 35K58 35K90)},
	MRNUMBER = {4894225},
	MRREVIEWER = {Reginaldo\ Demarque},
	DOI = {10.1007/s00028-025-01069-7},
	URL = {https://doi.org/10.1007/s00028-025-01069-7},
}

@incollection {JL99,
	AUTHOR = {Jerison, D. and Lebeau, G.},
	TITLE = {Nodal sets of sums of eigenfunctions},
	BOOKTITLE = {Harmonic analysis and partial differential equations
	({C}hicago, {IL}, 1996)},
	SERIES = {Chicago Lectures in Math.},
	PAGES = {223--239},
	PUBLISHER = {Univ. Chicago Press, Chicago, IL},
	YEAR = {1999},
	ISBN = {0-226-10456-7},
	MRCLASS = {58J05 (35B60)},
	MRNUMBER = {1743865},
	MRREVIEWER = {Philip\ W.\ Schaefer},
}

@article {Zhu24,
	AUTHOR = {Zhu, J.},
	TITLE = {Spectral inequalities for {S}chr\"odinger equations with various potentials},
	JOURNAL = {Trans. Amer. Math. Soc.},
	FJOURNAL = {Transactions of the American Mathematical Society},
	VOLUME = {379},
	YEAR = {2026},
	NUMBER = {8},
	PAGES = {5683--5712},
	ISSN = {0002-9947,1088-6850},
	MRCLASS = {35J10 (35P99 47A11 93B05)},
	MRNUMBER = {5093987},
	DOI = {10.1090/tran/9668},
	URL = {https://doi.org/10.1090/tran/9668},
}

@article{LM18,
	author={Logunov, A. and Malinnikova, E.},
	title={Quantitative propagation of smallness for solutions of elliptic equations}, 
	JOURNAL = {Proc. International Congress of Mathematicians (Rio de Janeiro, 2018)}, 
	year={2018},
	chapter = {},
	pages = {2391-2411},
	doi = {10.1142/9789813272880_0143},
	URL = {https://www.worldscientific.com/doi/abs/10.1142/9789813272880_0143},
	eprint = {https://www.worldscientific.com/doi/pdf/10.1142/9789813272880_0143},
}

@book {RS72,
	AUTHOR = {Reed, M. and Simon, B.},
	TITLE = {Methods of modern mathematical physics. {I}. {F}unctional
	analysis},
	PUBLISHER = {Academic Press, New York-London},
	YEAR = {1972},
	PAGES = {xvii+325},
	MRCLASS = {47-02 (81.47)},
	MRNUMBER = {493419},
	MRREVIEWER = {P.\ R.\ Chernoff},
}

@article {NTTV20,
	AUTHOR = {Naki\'c, I. and T\"aufer, M. and Tautenhahn, M.
	and Veseli\'c, I.},
	TITLE = {Sharp estimates and homogenization of the control cost of the
	heat equation on large domains},
	JOURNAL = {ESAIM Control Optim. Calc. Var.},
	FJOURNAL = {ESAIM. Control, Optimisation and Calculus of Variations},
	VOLUME = {26},
	YEAR = {2020},
	PAGES = {Paper No. 54, 26},
	ISSN = {1292-8119,1262-3377},
	MRCLASS = {93C20 (35K10 35Q93 93B05)},
	MRNUMBER = {4145245},
	MRREVIEWER = {Mohamed\ Ouzahra},
	DOI = {10.1051/cocv/2019058},
	URL = {https://doi.org/10.1051/cocv/2019058},
}

@article {AM22,
	AUTHOR = {Alphonse, P. and Martin, J.},
	TITLE = {Stabilization and approximate null-controllability for a large
	class of diffusive equations from thick control supports},
	JOURNAL = {ESAIM Control Optim. Calc. Var.},
	FJOURNAL = {ESAIM. Control, Optimisation and Calculus of Variations},
	VOLUME = {28},
	YEAR = {2022},
	PAGES = {Paper No. 16, 30},
	ISSN = {1292-8119,1262-3377},
	MRCLASS = {93D15 (26A33 35R11 93B05 93C20)},
	MRNUMBER = {4385099},
	MRREVIEWER = {Debayan\ Maity},
	DOI = {10.1051/cocv/2022009},
	URL = {https://doi.org/10.1051/cocv/2022009},
}

@book {OLBC10,
	TITLE = {N{IST} handbook of mathematical functions},
	EDITOR = {Olver, F. W. J. and Lozier, D. W. and Boisvert, R.
	F. and Clark, C. W.},
	NOTE = {With 1 CD-ROM (Windows, Macintosh and UNIX)},
	PUBLISHER = {U.S. Department of Commerce, National Institute of Standards
	and Technology, Washington, DC; Cambridge University Press,
	Cambridge},
	YEAR = {2010},
	PAGES = {xvi+951},
	ISBN = {978-0-521-14063-8},
	MRCLASS = {33-00 (00A20 65-00)},
	MRNUMBER = {2723248},
}

@article {EG18,
	AUTHOR = {Eremenko, A. and Gabrielov, A.},
	TITLE = {P{T}-symmetric eigenvalues for homogeneous potentials},
	JOURNAL = {J. Math. Phys.},
	FJOURNAL = {Journal of Mathematical Physics},
	VOLUME = {59},
	YEAR = {2018},
	NUMBER = {5},
	PAGES = {053503, 17},
	ISSN = {0022-2488,1089-7658},
	MRCLASS = {34L40 (30B10 81Q05 81Q12)},
	MRNUMBER = {3800598},
	MRREVIEWER = {G\'eza\ L\'evai},
	DOI = {10.1063/1.5016390},
	URL = {https://doi.org/10.1063/1.5016390},
}

@article{MZ25,
	title={Spectral inequalities for {S}chr\"odinger equations and quantitative propagation of smallness in the plane}, 
	author={Malinnikova, E. and Zhu, J.},
	year={2025},
	eprint={2505.03996},
	archivePrefix={arXiv},
	primaryClass={math.AP},
    note = {arXiv:2505.03996v1 [math.AP]},
	url={https://arxiv.org/abs/2505.03996}, 
}

@article {BMM15,
	AUTHOR = {Beauchard, K. and Miller, L. and Morancey, M.},
	TITLE = {2{D} {G}rushin-type equations: minimal time and null controllable data},
	JOURNAL = {J. Differential Equations},
	FJOURNAL = {Journal of Differential Equations},
	VOLUME = {259},
	YEAR = {2015},
	NUMBER = {11},
	PAGES = {5813--5845},
	ISSN = {0022-0396,1090-2732},
	MRCLASS = {93B05 (35K20 35K65 93B07 93C20)},
	MRNUMBER = {3397310},
	MRREVIEWER = {Dan\ Tiba},
	DOI = {10.1016/j.jde.2015.07.007},
	URL = {https://doi.org/10.1016/j.jde.2015.07.007},
}

@article {ABM21,
	AUTHOR = {Allonsius, D. and Boyer, F. and Morancey, M.},
	TITLE = {Analysis of the null controllability of degenerate parabolic systems of {G}rushin type via the moments method},
	JOURNAL = {J. Evol. Equ.},
	FJOURNAL = {Journal of Evolution Equations},
	VOLUME = {21},
	YEAR = {2021},
	NUMBER = {4},
	PAGES = {4799--4843},
	ISSN = {1424-3199,1424-3202},
	MRCLASS = {93B05 (30E05 35P20 93C20)},
	MRNUMBER = {4350587},
	MRREVIEWER = {Reginaldo\ Demarque},
	DOI = {10.1007/s00028-021-00733-y},
	URL = {https://doi.org/10.1007/s00028-021-00733-y},
}

@article{Tam26,
	title={Minimal time for null controllability of the parabolic spherical {B}aouendi-{G}rushin equation}, 
	author={Tamekue, C.},
	year={2026},
	eprint={2604.19700},
	archivePrefix={arXiv},
	primaryClass={math.AP},
    note={arXiv:2604.19700 [math.AP]},
	url={https://arxiv.org/abs/2604.19700}, 
}

@article {BDE20,
	AUTHOR = {Beauchard, K. and Dard\'e, J. and Ervedoza, S.},
	TITLE = {Minimal time issues for the observability of {G}rushin-type equations},
	JOURNAL = {Ann. Inst. Fourier (Grenoble)},
	FJOURNAL = {Universit\'e{} de Grenoble. Annales de l'Institut Fourier},
	VOLUME = {70},
	YEAR = {2020},
	NUMBER = {1},
	PAGES = {247--312},
	ISSN = {0373-0956,1777-5310},
	MRCLASS = {93B07 (35K65)},
	MRNUMBER = {4105940},
	MRREVIEWER = {Can\ Zhang},
	DOI = {10.5802/aif.3313},
	URL = {https://doi.org/10.5802/aif.3313},
}

@article {Koe17,
	AUTHOR = {Koenig, A.},
	TITLE = {Non-null-controllability of the {G}rushin operator in 2{D}},
	JOURNAL = {C. R. Math. Acad. Sci. Paris},
	FJOURNAL = {Comptes Rendus Math\'ematique. Acad\'emie des Sciences. Paris},
	VOLUME = {355},
	YEAR = {2017},
	NUMBER = {12},
	PAGES = {1215--1235},
	ISSN = {1631-073X,1778-3569},
	MRCLASS = {93B05 (35H10 35K10 93C20)},
	MRNUMBER = {3730500},
	MRREVIEWER = {Sergio\ Polidoro},
	DOI = {10.1016/j.crma.2017.10.021},
	URL = {https://doi.org/10.1016/j.crma.2017.10.021},
}

@article {DKR23,
	AUTHOR = {Dard\'e, J. and Koenig, A. and Royer, J.},
	TITLE = {Null-controllability properties of the generalized two-dimensional {B}aouendi-{G}rushin equation with non-rectangular control sets},
	JOURNAL = {Ann. H. Lebesgue},
	FJOURNAL = {Annales Henri Lebesgue},
	VOLUME = {6},
	YEAR = {2023},
	PAGES = {1479--1522},
	ISSN = {2644-9463},
	MRCLASS = {93B05 (34C20 35K20 35K65 47A10)},
	MRNUMBER = {4717395},
	MRREVIEWER = {Dany\ Nina-Huaman},
	DOI = {10.5802/ahl.193},
	URL = {https://doi.org/10.5802/ahl.193},
}

@article {LZ17,
	AUTHOR = {Liu, H. and Zhang, C.},
	TITLE = {Observability from measurable sets for a parabolic equation
	involving the {G}rushin operator and applications},
	JOURNAL = {Math. Methods Appl. Sci.},
	FJOURNAL = {Mathematical Methods in the Applied Sciences},
	VOLUME = {40},
	YEAR = {2017},
	NUMBER = {10},
	PAGES = {3821--3832},
	ISSN = {0170-4214,1099-1476},
	MRCLASS = {93B07 (35K20 35K67)},
	MRNUMBER = {3664612},
	MRREVIEWER = {Peng\ Gao},
	DOI = {10.1002/mma.4266},
	URL = {https://doi.org/10.1002/mma.4266},
}

@article {AS24,
	AUTHOR = {Alphonse, P. and Seelmann, A.},
	TITLE = {Quantitative spectral inequalities for the anisotropic {S}hubin operators and applications to null-controllability},
	JOURNAL = {C. R. Math. Acad. Sci. Paris},
	FJOURNAL = {Comptes Rendus Math\'ematique. Acad\'emie des Sciences. Paris},
	VOLUME = {362},
	YEAR = {2024},
	PAGES = {1635--1659},
	ISSN = {1631-073X,1778-3569},
	MRCLASS = {35P05 (93B05 93C20)},
	MRNUMBER = {4834571},
	DOI = {10.5802/crmath.670},
	URL = {https://doi.org/10.5802/crmath.670},
}

@article {AT13,
	AUTHOR = {Anh, C. T. and Toi, V. M.},
	TITLE = {Null controllability of a parabolic equation involving the {G}rushin operator in some multi-dimensional domains},
	JOURNAL = {Nonlinear Anal.},
	FJOURNAL = {Nonlinear Analysis. Theory, Methods \& Applications. An
	International Multidisciplinary Journal},
	VOLUME = {93},
	YEAR = {2013},
	PAGES = {181--196},
	ISSN = {0362-546X,1873-5215},
	MRCLASS = {93B05 (35K65 93C20)},
	MRNUMBER = {3117159},
	DOI = {10.1016/j.na.2013.08.003},
	URL = {https://doi.org/10.1016/j.na.2013.08.003},
}

@article {EV18,
	AUTHOR = {Egidi, M. and Veseli\'c, I.},
	TITLE = {Sharp geometric condition for null-controllability of the heat equation on {$\Bbb {R}^d$} and consistent estimates on the control cost},
	JOURNAL = {Arch. Math. (Basel)},
	FJOURNAL = {Archiv der Mathematik},
	VOLUME = {111},
	YEAR = {2018},
	NUMBER = {1},
	PAGES = {85--99},
	ISSN = {0003-889X,1420-8938},
	MRCLASS = {93B05 (35L20 35Q93 93C20)},
	MRNUMBER = {3816981},
	MRREVIEWER = {Cristina\ Pignotti},
	DOI = {10.1007/s00013-018-1185-x},
	URL = {https://doi.org/10.1007/s00013-018-1185-x},
}

@misc{DT26,
      title={On the critical time of observability of the multi-dimensional {B}aouendi-{G}rushin equation}, 
      author={Dardé, J. and Trabut, M.},
      year={2026},
      eprint={2604.01419},
      archivePrefix={arXiv},
      primaryClass={math.OC},
      note={arXiv:2604.01419 [math.OC]},
      url={https://arxiv.org/abs/2604.01419}, 
}

	
\end{document}